\documentclass{article}
\usepackage[utf8]{inputenc}
\usepackage{amssymb}
\usepackage{latexsym}
\usepackage{amsmath}
\usepackage{array}
\usepackage{multirow}
\usepackage{graphicx}
\graphicspath{Images/}
\usepackage{caption}
\usepackage{subcaption}
\usepackage{algorithmic}
\usepackage{float}
\usepackage{hyperref}
\usepackage{url} 
\usepackage{cleveref}
\usepackage{algorithm}
\usepackage{amsthm}
\usepackage{verbatim}
\usepackage{xcolor}
\usepackage{geometry}  
\newtheorem{theorem}{Theorem}

\newcommand{\avg}[1]{\left\{\!\!\left\{#1\right\}\!\!\right\}}
\newcommand{\diff}[1]{\left[\!\left[#1\right]\!\right]}

\newcommand{\modL}[1]{\left|#1\right|}
\newcommand{\paraL}[1]{\left(#1\right)}
\newcommand{\sumintK}[1]{\sum_{\sigma \in \mathcal{E}(K) \bigcap \mathcal{E}_{\text{int}}} \frac{|\sigma|}{|K|}}
\newcommand{\sumintKK}[1]{\sum_{\sigma \in \mathcal{E}(K)} \frac{|\sigma|}{|K|}}
\newcommand{\sumK}[1]{\sum_{K \in \mathcal{T}} |K|}
\newcommand{\sumintall}[1]{\sum_{\sigma \in \mathcal{E}} |\sigma|}

\newcommand{\dx}[1]{\ \mathrm{d}x}
\newcommand{\dt}[1]{\ \mathrm{d}t}
\newcommand{\divh}[1]{\text{div}_h}

\title{An asymptotic preserving scheme satisfying entropy stability for the barotropic Euler system\footnote{ \textbf{Funding:} This work was funded by the DAAD-DST (German-India) Project based personnel exchange programme: Development and analysis of higher-order structure-preserving numerical methods for hyperbolic balance laws. \\
M. A. was funded by the Ministry of Education, Government of India, under the Prime Minister's Research Fellowship/grant PM/MHRD-19-17567.03, and by the Gutenberg Research College University Mainz. \\
M.L.-M. gratefully acknowledges the support of DFG Project 5258 3336 funded within Focused Programme SPP 2410 "Hyperbolic Balance Laws: Complexity, Scales and Randomness" and of the Mainz Institute of Multiscale Modeling. }}
\author{Megala Anandan\footnote{Institute for Mathematics, Johannes Gutenberg University of Mainz, 55128 Mainz, Germany. e-mail: manandan@uni-mainz.de}, M\'{a}ria Luk\'{a}\v{c}ov\'{a}-Medvid'ov\'{a}\footnote{Institute for Mathematics, Johannes Gutenberg University of Mainz, 55128 Mainz, Germany. e-mail: lukacova@uni-mainz.de} , S. V. Raghurama Rao\footnote{Indian Institute of Science, Bangalore-560012, India. e-mail: raghu@iisc.ac.in}}

\date{}

\begin{document}

\maketitle

%\begin{abstract} 
%The dimensionless barotropic Euler system has an associated entropy inequality corresponding
%to a convex entropy function, for all non-zero values of Mach number. Further, this system
%reduces to incompressible system of equations when Mach number approaches zero. Hence, in
%this paper, we attempt to develop a numerical scheme that preserves both the structures, i.e.,
%it is asymptotic preserving at the low Mach number limit, and entropy stable for different values
%of Mach number. A time semi-discrete implicit-explicit (IMEX) treatment is utilised to achieve
%the asymptotic preserving property. The entropy stabilities of different space discretisation
%strategies have been studied, and one of these strategies yield long-time entropy stability for
%different values of Mach number as validated in the numerical experiments.
%\end{abstract}

\begin{abstract}
In this paper we study structure-preserving numerical methods for low Mach number barotropic Euler equations. 
Besides their asymptotic preserving properties that are crucial in order to obtain uniformly consistent and stable approximations of the Euler equations in their singular limit as the Mach number approaches zero, our aim is also to preserve discrete entropy stability.
Suitable acoustic/advection splitting approach combined with time implicit-explicit approximations are used to achieve
the asymptotic preserving property. The entropy stability of different space discretisation
strategies is studied for
different values of Mach number and is validated by the numerical experiments.
\end{abstract}

\section{Introduction}
Many problems arising in science and engineering often contain dimensionless parameters that appear when suitable non-dimensionalisation is employed. For barotropic/full Euler systems, the Mach number appears as a parameter (denoted by $\epsilon$), which dictates whether the flow is compressible ( $\epsilon \sim \mathcal{O}(1)$) or incompressible ($\epsilon<<1$).  The appearance of such a small parameter in the denominator also leads to stiffness in the system, reflected in wide disparity of propagating wave speeds. It has been rigorously proved in \cite{Basic_cont_Majda1,Basic_cont_Majda2,Basic_cont_Schochet} that the solutions of such hyperbolic systems converge to those of mixed hyperbolic-elliptic incompressible system when $\epsilon$ approaches zero. Explicit numerical methods for these systems require restrictive $\epsilon-$dependent stability condition on time step, and hence they become computationally very expensive when $\epsilon$ becomes small. Further, Godunov-type compressible flow solvers suffer from loss of accuracy as numerical dissipation is inversely proportional to $\epsilon$ \cite{Basic_exp_Dellacherie}. Fully implicit numerical methods, on the other hand, are complicated to implement due to nonlinearity of the Euler systems. Hence, attempts were made to efficiently extend the compressible flow solvers to incompressible/low Mach number limit consisting of divergence-free constraint on velocity field. In particular, semi-implicit time stepping techniques allow for the compressible solver to transform into an incompressible solver as $\epsilon$ becomes small, and its stability requirements are independent of $\epsilon$. Such schemes are called asymptotic preserving (AP) schemes, first introduced by Jin \cite{ap_jin} for kinetic equations and later extended to hyperbolic systems (see \cite{ap_jin2} for review). Semi-implicit time stepping is often achieved by implicit-explicit (IMEX) approach involving implicit treatment of stiff terms and explicit treatment of non-stiff terms. Several IMEX-AP schemes have been formulated by using different strategies to split the flux into stiff and non-stiff parts, and we refer the interested reader to \cite{Basic_AP_Degond1,Basic_AP_Tang,Basic_AP_Degond2,Basic_AP_Jin,Basic_AP_Lukacova1,Basic_AP_Lukacova2,Basic_AP_Lukacova3,Basic_AP_Dimarco1,Basic_AP_Boscarino1,Basic_AP_Dimarco2,Basic_AP_Noelle,Basic_AP_Boscarino2,Basic_AP_Lukacova4,Basic_AP_Arun,AP_Saurav}. \\
In addition to asymptotic preserving properties, another crucially required property for a numerical method is its stability. For the Euler equations, this also means nonlinear stability dictated by the second law of thermodynamics, the entropy inequality. Consequently, entropy stability has emerged as a nonlinear stability criterion for numerical schemes since the seminal work of Tadmor \cite{10.2307/2008251,tadmor_2003,TADMOR2016}. Several entropy stable numerical methods for different hyperbolic systems have been developed. These include developments specific to shallow water equations \cite{GASSNER2016291,WINTERMEYER2017200,https://doi.org/10.1002/fld.4766,VK_Megala}, the Euler equations \cite{Barth,ISMAIL20095410,puppo_semplice_2011-A,puppo_semplice_2011,chandrashekar_2013,Ray2013335,ray_chandrashekar_fjordholm_mishra_2016,GASSNER201639,CREAN2018410,Chizari2021,Yan2023}, 
%Navier-Stokes equations \cite{YAMALEEV2019108897,MANZANERO2020109241,RDP2020}, 
and magnetohydrodynamics equations \cite{doi:10.1137/15M1013626}. However, these entropy stable schemes were proposed for fixed Mach number ($\epsilon$ being order one). On the other hand, the governing system is associated with an entropy inequality for all non-zero values of $\epsilon$. Hence, our aim in this paper is to develop a numerical scheme that is entropy stable for different values of $\epsilon$, and also AP as $\epsilon$ approaches zero. As fas as we are aware, this question has not yet been studied in literature. The present paper makes the first step in this research direction, discusses possible discretisation strategies, and validates them by a series of numerical experiments.  \\
The paper is organised as follows: Section \ref{AP_ES_Sec:MM} presents the barotropic Euler system, its entropy stability property for different values of Mach number $\epsilon$, and its asymptotic limit as $\epsilon$ approaches zero. Section \ref{AP_ES_Sec:NM} presents the numerical method that employs an IMEX-AP time discretisation in the spirit of \cite{Basic_AP_Dimarco1,Basic_AP_Boscarino2}, and three different space discretisation strategies. The asymptotic preserving property of fully discrete scheme is also presented. Section \ref{AP_ES_Sec:NR} presents the numerical validation of our scheme by depicting the AP and entropy stability properties. Section \ref{AP_ES_Sec:Con} concludes the paper.

\section{Mathematical model}
\label{AP_ES_Sec:MM}
In this section, we present the barotropic Euler system, its entropy stability property, and its asymptotic limit as Mach number approaches $0$.
\subsection{The barotropic Euler system}
Consider the barotropic Euler system,
\begin{gather}
\label{Bar Euler mass}
    \partial_t \rho + \nabla \cdot (\rho \mathbf{u}) = 0 \\
\label{Bar Euler mom}
    \partial_t (\rho \mathbf{u}) + \nabla \cdot (\rho \mathbf{u} \otimes \mathbf{u} ) + \nabla p(\rho) = \mathbf{0},
\end{gather}
where $\mathbf{x} \in \Omega \subset \mathbb{R}^d$, $t\in \mathbb{R}^+ \cup \{0\}$, $\rho (\mathbf{x},t) : \Omega \times \mathbb{R}^+ \cup \{0\}  \rightarrow \mathbb{R}^+$ is the fluid density, $\mathbf{u}(\mathbf{x},t) : \Omega \times \mathbb{R}^+ \cup \{0\} \rightarrow \mathbb{R}^d$ is the fluid velocity, and $p(\rho(\mathbf{x},t))=\kappa \rho^{\gamma} \in \mathbb{R}^+$ is the pressure. Here, $d$ is the dimension in space, and $\kappa$, $\gamma > 1$ are constants. This system is hyperbolic with eigenvalues (in direction $\mathbf{n}$) $\mathbf{u}\cdot \mathbf{n} - c$ and $\mathbf{u}\cdot \mathbf{n} + c$, where $c=\sqrt{\gamma p/\rho}$ is the sound speed, and the conserved quantities are density, $\rho$ and momentum, $\rho \mathbf{u}$. The initial conditions required for the system are $\rho (\mathbf{x},0)=\rho^0(\mathbf{x})$ and $\mathbf{u}(\mathbf{x},0)=\mathbf{u}^0(\mathbf{x})$, and the boundary is considered to have periodic or zero flux conditions.\\
We perform non-dimensionalization of the above barotropic Euler system in \eqref{Bar Euler mass} and \eqref{Bar Euler mom} by using the reference values $x_r,t_r,\rho_r,u_r,p_r$. The dimensionless variables are given as, 
\begin{equation}
    \hat{\mathbf{x}}=\frac{\mathbf{x}}{x_r}, \hat{t}=\frac{t}{t_r}, \hat{\rho}=\frac{\rho}{\rho_r}, \hat{\mathbf{u}}=\frac{\mathbf{u}}{u_r}, \hat{p}=\frac{p}{p_r}.
\end{equation}
Inserting these into \eqref{Bar Euler mass} and \eqref{Bar Euler mom} and omitting the hat symbol, we obtain the dimensionless barotropic Euler system,
\begin{gather}
    \label{ND Bar Euler mass}
    \partial_t \rho + \nabla \cdot (\rho \mathbf{u}) = 0 \\
\label{ND Bar Euler mom}
    \partial_t (\rho \mathbf{u}) + \nabla \cdot (\rho \mathbf{u} \otimes \mathbf{u} ) + \frac{1}{\epsilon^2} \nabla p(\rho) = \mathbf{0},
\end{gather}
where $\epsilon = u_r \sqrt{\rho_r/p_r}$ is proportional to the Mach number. This system is also hyperbolic, and its eigenvalues (in direction $\mathbf{n}$) are $\mathbf{u}\cdot \mathbf{n} - c/\epsilon$ and $\mathbf{u}\cdot \mathbf{n} + c/\epsilon$. Hereafter, we consider the dimensionless form of barotropic Euler system in \eqref{ND Bar Euler mass} and \eqref{ND Bar Euler mom} for the presentation of analysis and numerical methods.  

\subsection{Entropy stability property}
Most hyperbolic systems in general have entropy inequality associated with them. In this section, we present the entropy inequality corresponding to the system in \eqref{ND Bar Euler mass} and \eqref{ND Bar Euler mom}. As we will see in what follows, the physical energy plays the role of (mathematical) entropy. Consequently, the entropy inequality reduces to the energy dissipation property.  \\
Let $\mathbf{U}=[\rho, \rho u_1, \dots , \rho u_d]^T$ be the vector of conserved variables and its flux vector in $k^{th}$ direction be denoted by 
$\mathbf{G}^k(\mathbf{U})=[\rho u_k, p\delta_{k1}/\epsilon^2 +\rho u_1 u_k, \dots , p\delta_{kd}/\epsilon^2 +\rho u_d u_k]^T$. Here $u_i$ is the $i^{th}$ component of fluid velocity $\mathbf{u}$. In this notation, the barotropic Euler system in \eqref{ND Bar Euler mass} and \eqref{ND Bar Euler mom} can be recast as:
\begin{equation}
\label{ND sys}
    \partial_t \mathbf{U} + \partial_{x_k} \mathbf{G}^k (\mathbf{U}) = \mathbf{0}.
\end{equation}
%Let us further assume that $\Omega$ is a convex subset of $\mathbb{R}^d$. 
The function, 
\begin{equation}
\label{Ent fn}
    \eta (\mathbf{U})= \frac{1}{2} \rho \| \mathbf{u} \|_2^2 + \frac{1}{\epsilon^2} \frac{p(\rho)}{\gamma -1}
\end{equation}
which is convex with respect to $\mathbf{U}$ is an entropy for the system \eqref{ND sys}. Further $ \frac{\text{d}^2\eta}{\text{d}\mathbf{U}^2}$ is a symmetric positive definite matrix and,
\begin{equation}
\label{Ent condn}
    \frac{\text{d}^2\eta}{\text{d}\mathbf{U}^2} \cdot \frac{\text{d}\mathbf{G}^k}{\text{d}\mathbf{U}} \text{ is symmetric} \iff \frac{\text{d}\omega_k}{\text{d}\mathbf{U}}= \frac{\text{d}\eta}{\text{d}\mathbf{U}} \cdot \frac{\text{d}\mathbf{G}^k}{\text{d}\mathbf{U}}.
\end{equation}
Here $\omega_k$ is the $k^{th}$ component of the entropy flux function $\boldsymbol{\omega}(\mathbf{U})=\mathbf{u} \left( \eta (\mathbf{U}) + p(\rho)/\epsilon^2 \right)$ corresponding to $\eta (\mathbf{U})$. For sufficiently smooth solutions, the inner product of \eqref{ND sys} with $\frac{\text{d} \eta}{\text{d} \mathbf{U}}$ gives entropy equality
\begin{equation}
    \partial_t \eta (\mathbf{U}) + \partial_{x_k} \omega_k (\mathbf{U}) = 0.
\end{equation}
For weak (non-smooth) solutions, we only get
\begin{equation}
\label{Ent ineq}
    \partial_t \eta (\mathbf{U}) + \partial_{x_k} \omega_k (\mathbf{U}) \leq 0
\end{equation}
due to the convexity of $\eta (\mathbf{U})$. Note that \eqref{Ent ineq} is understood in the distributional sense.  

\subsection{Asymptotic limit}
Our aim in this section is to define a limiting system of \eqref{ND Bar Euler mass}, \eqref{ND Bar Euler mom} as $\epsilon \to 0$. We point out that all calculations presented below are formal assuming enough regularity of the corresponding solutions. We assume that solutions can be expanded with respect to $\epsilon-$powers as follows:
\begin{gather}
\label{AP ansatz 1}
    \rho = \rho_0 + \epsilon \rho_1 + \epsilon^2 \rho_2 + \dots, \\
\label{AP ansatz 2}
    \mathbf{u} = \mathbf{u}_0 + \epsilon \mathbf{u}_1 + \epsilon^2 \mathbf{u}_2 + \dots, \\
\label{AP ansatz 3}
    p = p_0 + \epsilon p_1 + \epsilon^2 p_2 + \dots
\end{gather}
The asymptotic behavior as $\epsilon \to 0$ is determined by inserting \eqref{AP ansatz 1}, \eqref{AP ansatz 2} and \eqref{AP ansatz 3} into the system in \eqref{ND Bar Euler mass} and \eqref{ND Bar Euler mom}. Balancing $\mathcal{O}(\epsilon^{-2})$ terms in the momentum conservation equation, we obtain,
\begin{equation*}
    \nabla p_0 = 0.
\end{equation*}
Hence, $p_0$ is spatially constant and is function of time alone. Since $p_0=\kappa \rho_0^{\gamma}$, $\rho_0$ is also spatially constant and is function of time alone. Similarly balancing $\mathcal{O}(\epsilon^{-1})$ terms in the momentum conservation equation, we infer that $p_1$ and $\rho_1$ are also spatial constants and are functions of only time. Now, balancing $\mathcal{O}(1)$ terms in both mass and momentum conservation equations, we get,
\begin{gather}
\label{Order 1 Mass}
    \partial_t \rho_0 + \rho_0 \nabla \cdot \mathbf{u}_0 = 0, \\
\label{Order 1 Mom}
    \partial_t (\rho_0 \mathbf{u}_0) + \rho_0 \nabla \cdot (\mathbf{u}_0 \otimes \mathbf{u}_0 ) + \nabla p_2 = \mathbf{0}.
\end{gather}
Here $p_2$ is interpreted as the hydrostatic pressure. Integrating the $\mathcal{O}(1)$ mass balance in \eqref{Order 1 Mass} on $\Omega$, we get,
\begin{equation}
    |\Omega| \partial_t \rho_0 = - \rho_0 \int_{\Omega} \nabla \cdot \mathbf{u}_0 d \Omega = - \rho_0 \int_{\partial \Omega} \mathbf{u}_0 \cdot \mathbf{n} ds. 
\end{equation}
Taking $\mathbf{u} \cdot \mathbf{n} = 0$ on $\partial \Omega$ or considering periodic boundary conditions, we get $\int_{\partial \Omega} \mathbf{u}_0 \cdot \mathbf{n} ds = 0$. Thus, $\partial_t \rho_0 =0$ and $\rho_0$ is constant in both space and time, resulting in $\nabla \cdot \mathbf{u}_0 = 0$ according to \eqref{Order 1 Mass}. The $\mathcal{O}(1)$ momentum balance in \eqref{Order 1 Mom} therefore becomes,
\begin{equation}
    \partial_t \mathbf{u}_0 + \nabla \cdot (\mathbf{u}_0 \otimes \mathbf{u}_0 ) + \frac{\nabla p_2}{\rho_0} = \mathbf{0}.
\end{equation}
Similarly, integration of $\mathcal{O}(\epsilon)$ mass balance equation and usage of $\mathbf{u} \cdot \mathbf{n} = 0$ on $\partial \Omega$ or periodic boundary conditions result in $\partial_t \rho_1 =0$, and hence $\rho_1$ is constant in both space and time. \\
Further, the initial conditions are assumed to be compatible with the equations of different orders of $\epsilon$ (such as, $\epsilon^{-2}$, $\epsilon^{-1}$, $\epsilon^{0}$). In this paper, we consider the well-prepared initial conditions,
\begin{gather}
\label{wp rho ic}
    \rho(\mathbf{x},0)=\rho^0(\mathbf{x})=\rho_0^0 +  \epsilon^2 \rho_2^0 (\mathbf{x}) \\
\label{wp u ic}
    u(\mathbf{x},0)=u^0(\mathbf{x})=u_0^0 (\mathbf{x}) + \epsilon u_1^0 (\mathbf{x})
\end{gather}
such that $\rho_0^0$ is constant and $\nabla \cdot \mathbf{u}_0^0 = 0$.

\section{Numerical method}
\label{AP_ES_Sec:NM}
In this section, we want to construct a numerical method that is both asymptotic preserving and entropy stable. That is, we expect the method to satisfy the asymptotic limits of dimensionless barotropic Euler system in \eqref{ND Bar Euler mass} and \eqref{ND Bar Euler mom} as $\epsilon \to 0$, and also satisfy discrete entropy inequality in different regimes of $\epsilon$. To achieve this goal, we use implicit-explicit (IMEX) time discretisation required for attaining asymptotic consistency, and compare the entropy stability property of three different types of space discretisation in different regimes of $\epsilon$. We also present the asymptotic preserving property of considered numerical methods. 
\subsection{Semi-discrete IMEX time discretisation}
We begin with the presentation of first order IMEX time discretisation of the barotropic Euler system in \eqref{ND Bar Euler mass} and \eqref{ND Bar Euler mom} for clarity. 
\begin{gather}
\label{sd time_mass}
    \rho^{n+1}=\rho^n - \Delta t_n \nabla \cdot (\rho \mathbf{u})^{n+1} \\
\label{sd time_mom}
    (\rho \mathbf{u})^{n+1} = (\rho \mathbf{u})^{n} - \Delta t_n \nabla \cdot (\rho \mathbf{u} \otimes \mathbf{u} )^{n} - \frac{\Delta t_n}{\epsilon^2} \nabla p(\rho)^{n+1} 
\end{gather}
Here, $\Delta t_n=t_{n+1}-t_n$. The mass flux $\nabla \cdot (\rho \mathbf{u})$ and the pressure term $\frac{1}{\epsilon^2} \nabla p(\rho)$ are treated implicitly, while $\nabla \cdot (\rho \mathbf{u} \otimes \mathbf{u} )$ in the momentum flux is treated explicitly. It is important to treat the mass flux implicitly in order to get $\nabla \cdot \mathbf{u}^{n+1}=0$ as $\mathcal{O}(1)$ constraint. Indeed, if the mass flux is treated explicitly, then the whole method would become explicit and require severe $\epsilon$ dependent time step restriction enforced by stability. \\
Substituting the momentum equation \eqref{sd time_mom} in $\nabla \cdot (\rho \mathbf{u})^{n+1}$ of \eqref{sd time_mass}, we get,
\begin{equation}
\label{sd time_mass2}
    \rho^{n+1}=\rho^n - \Delta t_n \nabla \cdot (\rho \mathbf{u})^{n} + \Delta t_n^2 \nabla^2 : (\rho \mathbf{u} \otimes \mathbf{u} )^{n} + \frac{\Delta t_n^2}{\epsilon^2} \Delta p(\rho)^{n+1}.
\end{equation}
Since $p(\rho)=\kappa \rho^{\gamma}$, the presence of $\Delta p(\rho)^{n+1}$ in the above equation calls for a need to use the nonlinear iterative solver to find $\rho^{n+1}$. To avoid the computational effort, we perform linearisation of $p(\rho)^{n+1}$ around the incompressible constant density $\rho_0$ as:
\begin{equation}
    p(\rho)^{n+1}=p(\rho_0) + (\rho^{n+1}-\rho_0) p'(\rho)|_{\rho=\rho_0} + \mathcal{O}(\epsilon^4).
\end{equation}
The above linearisation is true if the higher derivatives of $p$ are $\mathcal{O}(1)$ and the method is asymptotic preserving (that is, $(\rho^{n+1}-\rho_0) \simeq \mathcal{O}(\epsilon^2)$). We intend to construct our method such that it is asymptotic preserving, and we have used this information \textit{a priori} in the linearisation of $p(\rho)^{n+1}$. Using this linearisation in \eqref{sd time_mass2}, we get,
\begin{equation}
\label{sd time_mass3}
    \rho^{n+1}=\rho^n - \Delta t_n \nabla \cdot (\rho \mathbf{u})^{n} + \Delta t_n^2 \nabla^2 : (\rho \mathbf{u} \otimes \mathbf{u} )^{n} + \frac{\Delta t_n^2}{\epsilon^2} p'(\rho)|_{\rho=\rho_0} \Delta \rho^{n+1} + \mathcal{O}(\Delta t_n^2 \epsilon^2).
\end{equation}
In a crude sense (by formally Taylor expanding $\rho^n$ about $\rho^{n+1}$), the modified or equivalent partial differential equation of the above time discrete equation is obtained as,
\begin{eqnarray*}
    \rho^{n+1} &=& \rho^n - \Delta t_n \nabla \cdot (\rho \mathbf{u})^{n} + \Delta t_n^2 \nabla^2 : (\rho \mathbf{u} \otimes \mathbf{u} )^{n} + \frac{\Delta t_n^2}{\epsilon^2} p'(\rho)|_{\rho=\rho_0} \Delta \rho^{n+1}  \\
    &=& \rho^n - \Delta t_n \nabla \cdot (\rho \mathbf{u})^{n+1} + \mathcal{O}(\Delta t_n^2 \epsilon^2)  \text{ (from \eqref{sd time_mom})}\\
    \implies \partial_t \rho^{n+1} &=& - \nabla \cdot (\rho \mathbf{u})^{n+1} + \mathcal{O}(\Delta t_n) +  \mathcal{O}(\Delta t_n \epsilon^2).
\end{eqnarray*}
Thus, the first order temporal accuracy of the method remains unaffected due to the linearisation as long as $\mathcal{O}(\epsilon^2) \leq \mathcal{O}(1)$. The higher derivatives are considered to be $\mathcal{O}(1)$ in this argument. \\
As indicated by \eqref{sd time_mass3}, \eqref{sd time_mom}, we have split the part governed by the acoustic waves from the rest non-stiff part. The latter models the nonlinear advection waves. \\
From the algorithmic viewpoint, \eqref{sd time_mass3} can be solved easily by inversion of a matrix as follows,
\begin{equation}
\label{sd time_1 mass}
    \rho^{n+1}=  \left( I - \left(\frac{\Delta t_n}{\epsilon}\right)^2  p'(\rho)|_{\rho=\rho_0} \Delta \right)^{-1}  \left(  \rho^n - \Delta t_n \nabla \cdot (\rho \mathbf{u})^{n} + \Delta t_n^2 \nabla^2 : (\rho \mathbf{u} \otimes \mathbf{u} )^{n} \right).
\end{equation}
Then, $\rho^{n+1}$ evaluated as above is used to find $p(\rho)^{n+1}$. Inserting this into \eqref{sd time_mom}, we get $(\rho \mathbf{u})^{n+1}$ and thus the algorithm is complete. \eqref{sd time_1 mass} and \eqref{sd time_mom} together form the update equations for first order time semi-discrete scheme. 

Next, we present the higher order IMEX Runge Kutta (IMEX-RK) time discretisation of the barotropic Euler system in \eqref{ND Bar Euler mass} and \eqref{ND Bar Euler mom}. An IMEX-RK time discretisation is represented by the following double Butcher tableau:
\begin{equation}
\label{AP_ES_gen Butcher tableau}
\begin{tabular}{p{0.25cm}|p{0.5cm}}
\centering $\Tilde{c}$ & \centering $\Tilde{A}$ \cr 
\hline
 & \centering \vspace{-0.25em}$\Tilde{b}^T$
\end{tabular} \quad \quad
\begin{tabular}{p{0.25cm}|p{0.5cm}}
\centering $c$ & \centering $A$ \cr 
\hline
 &  \vspace{-0.25em}\centering $b^T$ 
\end{tabular}
\end{equation}
where $\Tilde{A}=(\Tilde{a}_{ij}), A=(a_{ij}) \in \mathbb{R}^{s \times s}$; $c,\Tilde{c},b,\Tilde{b} \in \mathbb{R}^{s}$. The matrices $\Tilde{A},A$ correspond to explicit (strictly lower triangular matrix with diagonal elements as $0$) and implicit (lower triangular with non-zero diagonal elements) parts of the scheme. Such $A$ are known as \textit{diagonally implicit} matrices. The coefficients $\Tilde{c}$ and $c$ are given by  
\begin{equation}
    \Tilde{c}_i=\sum_{j=1}^{i-1} \Tilde{a}_{ij}, c_i=\sum_{j=1}^{i} a_{ij},
\end{equation} and the vectors $\Tilde{b}=(\Tilde{b}_j)$ and $b=(b_j)$ give quadrature weights that combine the stages. For AP schemes, it turns to be important to work with globally stiffly accurate (GSA) IMEX-RK scheme that satisfies the following property: 
\begin{equation}
\label{AP_ES_GSA}
    c_s=\Tilde{c}_s=1 \text{ and } a_{sj}=b_j, \Tilde{a}_{sj}=\Tilde{b}_j, \;\; \forall j \in \{1,2,\dots,s\}.  
\end{equation}
The GSA property ensures that the update at $t_{n+1}$ is same as the update at $s^{th}$ stage. \\
The $i^{th}$ stage update (for $i \in \{1,2, \dots,s \}$) of the barotropic Euler system in \eqref{ND Bar Euler mass} and \eqref{ND Bar Euler mom} is given by,
\begin{gather}
\label{sd stage time_mass}
    \rho^{i}=\rho^n - \Delta t_n \sum_{j=1}^i a_{ij} \nabla \cdot (\rho \mathbf{u})^{j} \\
\label{sd stage time_mom}
    (\rho \mathbf{u})^{i} = (\rho \mathbf{u})^{n} - \Delta t_n \sum_{j=1}^{i-1} \Tilde{a}_{ij} \nabla \cdot (\rho \mathbf{u} \otimes \mathbf{u} )^{j} - \frac{\Delta t_n}{\epsilon^2} \sum_{j=1}^i a_{ij} \nabla p(\rho)^{j} 
\end{gather}
where $\Delta t_n=t_{n+1}-t_n$. Substituting the momentum equation \eqref{sd stage time_mom} in $\nabla \cdot (\rho \mathbf{u})^{i}$ of \eqref{sd stage time_mass}, we get,
\begin{equation}
\label{sd stage time_mass2}
    \rho^{i}=\rho^n - \Delta t_n \sum_{j=1}^{i-1} a_{ij} \nabla \cdot (\rho \mathbf{u})^{j} - \Delta t_n a_{ii} \nabla \cdot (\rho \mathbf{u})^{n} + \Delta t_n^2 a_{ii} \sum_{j=1}^{i-1} \Tilde{a}_{ij} \nabla^2 : (\rho \mathbf{u} \otimes \mathbf{u} )^{j} + \frac{\Delta t_n^2}{\epsilon^2} a_{ii} \sum_{j=1}^i a_{ij} \Delta p(\rho)^{j}.
\end{equation}
The above equation requires a nonlinear solver to find $\rho^{i}$. Similar to first order method, we perform linearisation around incompressible constant density $\rho_0$: 
\begin{equation}
    p(\rho)^{i}=p(\rho_0) + (\rho^{i}-\rho_0) p'(\rho)|_{\rho=\rho_0} + \mathcal{O}(\epsilon^4).
\end{equation}
The asymptotic preserving property ($(\rho^{i}-\rho_0) \simeq \mathcal{O}(\epsilon^2)$) of the method is used a priori in the linearisation. Plugging in \eqref{sd stage time_mass2} yields,
\begin{multline}
\label{sd stage time_high mass}
    \rho^{i}= \left( I - \left(\frac{\Delta t_n}{\epsilon}\right)^2 a_{ii}^2 p'(\rho)|_{\rho=\rho_0} \Delta \right)^{-1} \\ \left(  \rho^n - \Delta t_n \sum_{j=1}^{i-1} a_{ij} \nabla \cdot (\rho \mathbf{u})^{j} - \Delta t_n a_{ii} \nabla \cdot (\rho \mathbf{u})^{n} + \Delta t_n^2 a_{ii} \sum_{j=1}^{i-1} \Tilde{a}_{ij}  \nabla^2 : (\rho \mathbf{u} \otimes \mathbf{u} )^{j} + \frac{\Delta t_n^2}{\epsilon^2} a_{ii} \sum_{j=1}^{i-1} a_{ij} p'(\rho)|_{\rho=\rho_0} \Delta \rho^{j} \right).
\end{multline}
Then, $\rho^{i}$ evaluated as above is used to find $p(\rho)^{i}$. Inserting this into \eqref{sd stage time_mom}, we get $(\rho \mathbf{u})^{i}$ and thus the evaluation of stage values is complete. \eqref{sd stage time_high mass} and \eqref{sd stage time_mom} together form the stage update equations for higher order IMEX-RK time semi-discrete scheme. Further, $\rho^{n+1}=\rho^s$ and $(\rho \mathbf{u})^{n+1}=(\rho \mathbf{u})^{s}$ due to the GSA property and therefore the algorithm is complete. 

\subsection{Asymptotic preserving property of the time semi-discrete scheme}
In this section, we show that the higher order GSA IMEX-RK time semi-discrete scheme \eqref{sd stage time_mass2} and \eqref{sd stage time_mom} is asymptotic preserving. 
\begin{theorem}
    Assume well-prepared initial conditions in \eqref{wp rho ic} and \eqref{wp u ic}, the asymptotic expansion in \eqref{AP ansatz 1}-\eqref{AP ansatz 3}, and periodic boundary conditions on $\rho$ and $\mathbf{u}$. Then the time semi-discrete GSA IMEX-RK scheme given by \eqref{sd stage time_mass2} and \eqref{sd stage time_mom} satisfies for $\epsilon \to 0$
    \begin{gather}
        \rho_0^i \equiv \text{ constant }, \  \rho_1^i \equiv \text{ constant }, \ \rho_0^i+\epsilon \rho_1^i = \rho_0, \\
        \nabla \cdot \mathbf{u}_0^i = 0, \\
        \mathbf{u}_0^i = \mathbf{u}_0^n - \Delta t_n \sum_{j=1}^{i-1} \Tilde{a}_{ij} \nabla \cdot (\mathbf{u}_0 \otimes \mathbf{u}_0)^{j} - \frac{\Delta t_n}{\rho_0} \sum_{j=1}^i a_{ij} \nabla p_2^{j},
    \end{gather}
for all $i \in \{1,2,\dots,s \}$, which is a consistent approximation of the incompressible Euler system \eqref{Order 1 Mass}, \eqref{Order 1 Mom}.
\end{theorem}
\begin{proof}
    Inserting the asymptotic expansion \eqref{AP ansatz 1}-\eqref{AP ansatz 3} into the momentum update equation \eqref{sd stage time_mom} and equating $\mathcal{O}\left( \frac{1}{\epsilon^2} \right)$ terms, we obtain
    \begin{equation*}
        \Delta t_n \sum_{j=1}^i a_{ij} \nabla p_0^{j} =0, \ \text{ for all } i \in \{1,2,\dots,s \} \implies \nabla p_0^{i} =0, \ \text{ for all } i \in \{1,2,\dots,s \}.
    \end{equation*}
    Since $p_0^i=\kappa \rho_0^{i^{\gamma}}$, $\rho_0^i$ is spatially constant for all $i \in \{1,2,\dots,s \}$. Similarly equating $\mathcal{O}\left( \frac{1}{\epsilon} \right)$ terms in the momentum update equation \eqref{sd stage time_mom}, we infer that $\rho_1^i$ is spatially constant for all $i \in \{1,2,\dots,s \}$. \\
    Inserting the asymptotic expansion \eqref{AP ansatz 1}-\eqref{AP ansatz 3} into the mass update equation \eqref{sd stage time_mass2} and equating $\mathcal{O}\left( 1 \right)$ terms, we obtain
    \begin{equation}
        \rho_0^{i}=\rho_0^n - \Delta t_n \sum_{j=1}^{i-1} a_{ij} \rho_0^j \nabla \cdot (\mathbf{u}_0)^{j} - \Delta t_n a_{ii} \rho_0^n \nabla \cdot (\mathbf{u}_0)^{n} + \Delta t_n^2 a_{ii} \sum_{j=1}^{i-1} \Tilde{a}_{ij} \rho_0^j  \nabla^2 : (\mathbf{u}_0 \otimes \mathbf{u}_0 )^{j} + \Delta t_n^2 a_{ii} \sum_{j=1}^i a_{ij}  \Delta p_2^{j}.
    \end{equation}
    Integrating the above equation on $\Omega$ and using periodic boundary conditions on $\rho_2$ and $\mathbf{u}_0$, we obtain 
    \begin{equation}
        \rho_0^i = \rho_0^n, \  \text{ for all } i \in \{1,2,\dots,s \}.
    \end{equation}
    Repeating the similar procedure for $\mathcal{O}\left( \epsilon \right)$ terms of the mass update equation \eqref{sd stage time_mass2}, we obtain,
    \begin{equation}
        \rho_1^i = \rho_1^n, \  \text{ for all } i \in \{1,2,\dots,s \}.
    \end{equation}
    Since $\rho_{0,1}^{n+1}=\rho_{0,1}^s$ due to the GSA property of IMEX-RK time discretisation, we have $\rho_{0,1}^{n+1}=\rho_{0,1}^s=\rho_{0,1}^n$. Therefore, $\rho_{0,1}^n=\rho_{0,1}^0 \equiv \text{constant}, \text{ for all } n=1,2,\dots$. \\
    Inserting this into the $\mathcal{O}(1)$ mass and momentum update equations, we get for all $i \in \{1,2,\dots,s \}$
    \begin{gather}
    \label{AP time proof mass}
        \sum_{j=1}^{i-1} a_{ij} \nabla \cdot (\mathbf{u}_0)^{j} + a_{ii} \nabla \cdot (\mathbf{u}_0)^{n} - \Delta t_n a_{ii} \sum_{j=1}^{i-1} \Tilde{a}_{ij}  \nabla^2 : (\mathbf{u}_0 \otimes \mathbf{u}_0 )^{j} - \frac{\Delta t_n}{\rho_0} a_{ii} \sum_{j=1}^i a_{ij}  \Delta p_2^{j} = 0, \\
    \label{AP time proof mom}
        \mathbf{u}_0^{i} = \mathbf{u}_0^{n} - \Delta t_n \sum_{j=1}^{i-1} \Tilde{a}_{ij} \nabla \cdot ( \mathbf{u}_0 \otimes \mathbf{u}_0)^{j} - \frac{\Delta t_n}{\rho_0} \sum_{j=1}^i a_{ij} \nabla p_2^{j}, 
    \end{gather}
    where $\rho_0 = \rho_0^0 + \epsilon \rho_1^0 $. Taking divergence of \eqref{AP time proof mom} and inserting it into \eqref{AP time proof mass}, we obtain
    \begin{equation}
        \sum_{j=1}^{i} a_{ij} \nabla \cdot (\mathbf{u}_0)^{j} = 0, \ \text{ for all } i \in \{1,2,\dots,s \} \implies \nabla \cdot (\mathbf{u}_0)^{i} =0, \ \text{ for all } i \in \{1,2,\dots,s \}.
    \end{equation}
    \begin{comment}
        \begin{gather}
        \rho_0^i &=& \left(- \Delta t_n^2 a_{ii}^2 p'(\rho)|_{\rho=\rho_0} \Delta \right)^{-1} \left( \Delta t_n^2 a_{ii} \sum_{j=1}^{i-1} a_{ij} p'(\rho)|_{\rho=\rho_0} \Delta \rho_0^{j}  \right) \\
        \implies & \ & \sum_{j=1}^i \frac{a_{ij}}{a_{ii}} \rho_0^j = 0
    \end{gather}
    $\forall i \in \{1,2,\dots,s \}$. 
    \end{comment}
    \end{proof}
The above theorem shows the asymptotic consistency of the IMEX-RK time semi-discrete scheme. Due to the GSA property, the expressions for $\rho_0^s, \rho_1^s, u_0^s$ follow for $\rho_0^{n+1}, \rho_1^{n+1}, u_0^{n+1}$. In the next section, we explain the discretisation techniques for spatial derivatives present in the scheme.  

\subsection{Space discretisation}
In this section, we discuss various consistent spatial discretisations for the time semi-discrete scheme proposed above. It is important to keep the numerical diffusion coefficients free of the small parameter $\epsilon$, in-order to avoid an uncontrollable growth in numerical diffusion term as $\epsilon \to 0$. In what follows, we present three different types of space discretisation.
\par We consider the first order time semi-discrete scheme given by \eqref{sd time_1 mass} and \eqref{sd time_mom} for presentation of the spatial discretisation. The discretisation of all the spatial derivatives present in this scheme will be explained. The corresponding spatial derivatives in higher order time semi-discrete scheme given by \eqref{sd stage time_high mass} and \eqref{sd stage time_mom} will be approximated analogously. The additional terms $\sum_{j=1}^{i-1} a_{ij} \nabla \cdot (\rho \mathbf{u})^{j}$ and $\sum_{j=1}^{i-1} a_{ij} \Delta p(\rho)^{j}$ present in the mass update equation (given by \eqref{sd stage time_high mass}) will also follow the same discretisation as $\nabla \cdot (\rho \mathbf{u})^n$ and $\Delta p(\rho)^{n+1}$, respectively. %For clarity, we explain the ideas in a multi-dimensional setting. 
\par Let us consider the mesh $\mathcal{T}_h$ of closed polygons or polyhedra on the physical domain, such that $\overline{\Omega} = \bigcup_{K\in \mathcal{T}_h} K$, $h>0$ is the mesh parameter. Let $\mathcal{E}(K)$ denote the set of faces $\sigma$ of the cell $K$, $\sigma \in \partial K$. For each face $\sigma = K|L$, $\mathbf{n}_{\sigma,K}$ denotes the normal vector of $\sigma$, that is pointed outwards from the element $K$ to the element $L$. Let $|K|$ and $|\sigma|$ denote the measure of cell $K$ and face $\sigma$, respectively. Let $\phi \in L^1(\Omega;\mathbb{R})$, then $\phi_h$ denotes a piecewise constant projection, \textit{i.e.}
\begin{equation}
    \phi_h|_K \equiv \phi_K = \frac{1}{|K|} \int_K \phi(x) \text{ d}x. 
\end{equation}
For each face $\sigma=K|L$, $\avg{\phi_h}_{\sigma} = \frac{\phi_K+\phi_L}{2}$ and $\diff{\phi_h}_{\sigma}=\phi_L-\phi_K$ stands for the average and jump over $\sigma$, respectively. Let $d_{\sigma}$ be the length of a line connecting the centres of $K$ and $L$ that is orthogonal to $\sigma$.

\subsubsection{Type 1}
We apply an upwind discretisation for $\nabla \cdot (\rho \mathbf{u} \otimes \mathbf{u})^n$ in the momentum equation \eqref{sd time_mom}, while all other first and second derivatives present in the scheme \eqref{sd time_1 mass} and \eqref{sd time_mom} are treated in a central fashion. Since our goal is to achieve entropy stability, we do not add numerical diffusion to implicit terms as they are entropy stable with central discretisation. 
\begin{enumerate}
    \item Spatial discretisation of $\nabla \cdot (\rho \mathbf{u} \otimes \mathbf{u})^n$ in momentum equation \eqref{sd time_mom} is given by upwind discretisation 
    \begin{equation}
        D_{upw}(\rho \mathbf{u} \otimes \mathbf{u})_K := \sumintKK{} \paraL{ \rho_K \mathbf{u}_K [\avg{\mathbf{u}_h}_{\sigma} \cdot \mathbf{n}_{\sigma,K}]^+ + \rho_L \mathbf{u}_L [\avg{\mathbf{u}_h}_{\sigma} \cdot \mathbf{n}_{\sigma,K}]^- }, 
    \end{equation}
    where $ [a]^{\pm} = \frac{(a)\pm \modL{a}}{2}$ stands for the positive and negative part of $a \in \mathbb{R}$, respectively.
    %The subscript $I$ denotes the interface, and the subscripts $L$ and $R$ denote the cell values to the left and right of interface, $I$. $\mathbf{n}$ denotes the normal to the interface $I$, and it is outward corresponding to the cell $L$ at which update takes place. 
    \item Spatial discretisation of $\nabla \cdot (\rho \mathbf{u})^{n}$ in mass update equation \eqref{sd time_1 mass} and $\nabla p(\rho)^{n+1}$ in momentum equation \eqref{sd time_mom} are given by central discretisations
    \begin{gather}
        D_{cen} (\rho \mathbf{u})_K := \sumintKK{} \avg{\rho_h \mathbf{u}_h}_{\sigma} \cdot \mathbf{n}_{\sigma,K} \\
        D_{cen} (p)_K :=  \sumintKK{} \avg{p(\rho_h)}_{\sigma} \mathbf{n}_{\sigma,K}.
    \end{gather} 
    \item Spatial discretisation of $\nabla^2 : (\rho \mathbf{u} \otimes \mathbf{u} )^{n}$ and $\Delta p^{n+1}$ in mass update equation \eqref{sd time_1 mass} are given by central discretisations 
    \begin{gather}
        D^2 (\rho \mathbf{u} \otimes \mathbf{u})_K :=  D_{cen} D_{cen} (\rho \mathbf{u} \mathbf{u})_K; \ D_{cen} (\cdot)_K := \sumintKK{} \avg{(\cdot)_h}_{\sigma} \cdot \mathbf{n}_{\sigma,K}  \\
        D^2 (p)_K := \sumintKK{} \frac{\diff{p(\rho_h)}_{\sigma}}{d_{\sigma}}. 
    \end{gather} 
\end{enumerate}
The corresponding terms in higher order time semi-discrete scheme given by \eqref{sd stage time_high mass} and \eqref{sd stage time_mom} also follow the same discretisation as described above. 

\subsubsection{Type 2}
All the terms in \eqref{sd time_1 mass} and \eqref{sd time_mom} follow the discretisation in type 1, except that the term $\nabla \cdot (\rho \mathbf{u})^n$ in the mass update equation \eqref{sd time_1 mass} is treated in an upwind fashion. Spatial discretisation for $\nabla \cdot (\rho \mathbf{u})^n$ is given by, 
\begin{equation}
    D_{upw}(\rho \mathbf{u})_K := \sumintKK{} \paraL{ \rho_K [\avg{\mathbf{u}_h}_{\sigma} \cdot \mathbf{n}_{\sigma,K}]^+ + \rho_L [\avg{\mathbf{u}_h}_{\sigma} \cdot \mathbf{n}_{\sigma,K}]^- }.
\end{equation}
The corresponding terms in higher order time semi-discrete scheme given by \eqref{sd stage time_high mass} and \eqref{sd stage time_mom} also follow the same discretisation as described above. 

\subsubsection{Type 3}
Here, all the terms in \eqref{sd time_1 mass} and \eqref{sd time_mom} follow the discretisation in type 1, except that the term $\nabla \cdot (\rho \mathbf{u} \otimes \mathbf{u})^n$ in momentum equation \eqref{sd time_mom} is discretised by using an entropy stable flux. We apply the entropy conservative flux discretisation, and add numerical diffusion for attaining entropy stability. We first derive the entropy conserving and stable fluxes for the full barotropic Euler system, and use only the part of these fluxes corresponding to $\nabla \cdot (\rho \mathbf{u} \otimes \mathbf{u})^n$ in momentum equation \eqref{sd time_mom}. All other first and second derivatives in the scheme \eqref{sd time_1 mass} and \eqref{sd time_mom} are treated in central fashion as in type 1. \\
The entropy variable corresponding to the convex entropy function \eqref{Ent fn} is: 
\begin{equation}
    \mathbf{V}=\frac{\text{d} \eta}{\text{d} \mathbf{U}} = \begin{bmatrix}-\frac{1}{2}\modL{\mathbf{u}}^2+\frac{1}{\epsilon^2}\frac{\kappa \gamma}{\gamma -1} \rho^{\gamma -1} & \mathbf{u} \end{bmatrix}^T
\end{equation} 
For entropy conservation, we require (cf. \cite{10.2307/2008251,tadmor_2003,TADMOR2016}) 
\begin{gather}
\label{Ent cons condn}
    \diff{\mathbf{V} \cdot \mathbf{G}^k}_{\sigma} - \diff{\mathbf{V}}_{\sigma} \cdot \mathbf{G}_{\sigma}^{k*} = \diff{\omega^k}_{\sigma},
\end{gather}
where $k \in \{1,2,\dots,d\}$ is an index for dimension.  The following interface flux function
\begin{gather}
    \mathbf{G}_{\sigma}^{k*} = \begin{bmatrix}
        \paraL{\rho u^{(k)}}_{\sigma}^* \\ \left( \frac{1}{\epsilon^2}p \delta_{kj}+ \rho u^{(k)} u^{(j)} \right)_{\sigma}^* \end{bmatrix} = \begin{bmatrix}
            \avg{\rho}^{\gamma}_{\sigma} \avg{u_h^{(k)}}_{\sigma} \\  \frac{1}{\epsilon^2} \avg{p(\rho_h)}_{\sigma} \delta_{kj} + \avg{\rho}^{\gamma}_{\sigma} \avg{u_h^{(k)}}_{\sigma} \avg{u_h^{(j)}}_{\sigma}
        \end{bmatrix}, \text{ for } j \in \{1,2,\dots,d\}
\end{gather}
with $\delta_{kj}$ denoting the Kronecker delta function, and
\begin{eqnarray}
    \avg{\rho}^{\gamma}_{\sigma} = \frac{\gamma -1}{\gamma} \frac{\diff{ \rho_h^{\gamma} }_{\sigma}}{ \diff{\rho^{\gamma -1}}_{\sigma}},
\end{eqnarray}
satisfies the entropy conserving condition in \eqref{Ent cons condn}. Hence, the entropy conserving spatial discretisation of $\nabla \cdot (\rho \mathbf{u} \otimes \mathbf{u})^n$ in momentum equation \eqref{sd time_mom} is given by 
\begin{gather}
\label{Space_t3_ec0}
    D_{EC}(\rho \mathbf{u} \otimes \mathbf{u})_K = \sumintKK{} \avg{\rho}^{\gamma}_{\sigma} \avg{\mathbf{u}_h}_{\sigma} \otimes \avg{\mathbf{u}_h}_{\sigma} \cdot \mathbf{n}_{\sigma,K}.
\end{gather} 
Note that \eqref{Space_t3_ec0} yields second order accurate approximation. \\
To achieve entropy stability, we consider a scalar dissipation $\Lambda^k = \modL{\avg{u_h^{(k)}}_{\sigma} }$ that is independent of $\epsilon$. 
Thus, the entropy stable flux for full barotropic Euler system reads, 
\begin{equation}
    \mathbf{G}^k_{\sigma} = \mathbf{G}_{\sigma}^{k*} - \frac{q}{2} \Lambda^k \diff{ \mathbf{V} }_{\sigma},
\end{equation} 
where $q>0$ is a suitable constant. This entropy stable flux results in first order spatial accuracy. For second order accuracy, we use (cf. \cite{doi:10.1137/110836961}) 
\begin{equation}
    \mathbf{G}^k_{\sigma} = \mathbf{G}_{\sigma}^{k*} - \frac{q}{2} \Lambda^k \left<\!\left< \mathbf{V} \right>\!\right>_{\sigma},
\end{equation}
where min-mod reconstruction is used to evaluate an average of a piecewise linear reconstruction on $\sigma$ denoted by $\left<\!\left< \mathbf{V} \right>\!\right>_{\sigma}$:
\begin{eqnarray}
\label{ES space disc type 3}
     D_{ES} (\rho \mathbf{u} \otimes \mathbf{u})_K  &=& \sumintKK{} \left\{ \begin{matrix} \modL{\avg{\mathbf{u}}_{\sigma}} \left( \avg{\rho}^{\gamma}_{\sigma} \modL{\avg{\mathbf{u}}_{\sigma}} - \frac{q}{2}  \diff{\mathbf{u} }_{\sigma} \right), & 1^{st} \text{ order} \\  \modL{\avg{\mathbf{u}}_{\sigma}} \left( \avg{\rho}^{\gamma}_{\sigma} \modL{\avg{\mathbf{u}}_{\sigma}}  - \frac{q}{2}  \left<\!\left<\mathbf{u} \right>\!\right>_{\sigma} \right),  & 2^{nd} \text{ order} . \end{matrix} \right.
\end{eqnarray}
Hence, in this type, $\nabla \cdot (\rho \mathbf{u} \otimes \mathbf{u})^n$ in momentum equation \eqref{sd time_mom} is discretised as shown in \eqref{ES space disc type 3}. All the other first and second derivatives present in \eqref{sd time_1 mass} and \eqref{sd time_mom} are discretised in central fashion as shown in type 1. \\
The corresponding terms in higher order time semi-discrete scheme given by \eqref{sd stage time_high mass} and \eqref{sd stage time_mom} also follow the same discretisation as described above. 

\subsection{Asymptotic preserving property of the fully discrete scheme}
In this section, we show the asymptotic consistency of our fully discrete scheme as $\epsilon \to 0$. For this, we present a general theorem that considers all the three types of spatial discretisation. For convenience and clarity, we prove the asymptotic preserving property in a one-dimensional setting. The results presented below directly generalize to a multi-dimensional case
\begin{theorem}
    Assume well-prepared initial conditions in \eqref{wp rho ic} and \eqref{wp u ic}, the asymptotic expansion in \eqref{AP ansatz 1}-\eqref{AP ansatz 3}, and periodic boundary conditions on $\rho$ and $u_1$. Consider the fully discrete scheme
    \begin{gather}  
    \label{AP thm full mass}
        \rho_k^{i}=\rho_k^n - \Delta t_n \sum_{j=1}^{i-1} a_{ij} D (\rho u_1)_k^{j} - \Delta t_n a_{ii} D (\rho u_1)_k^{n} + \Delta t_n^2 a_{ii} \sum_{j=1}^{i-1} \Tilde{a}_{ij} D^2 (\rho u_1^2 )_k^{j} + \frac{\Delta t_n^2}{\epsilon^2} a_{ii} \sum_{j=1}^i a_{ij} D^2 (p)_k^{j} \\
    \label{AP thm full mom}
        (\rho u_1)_k^{i} = (\rho u_1)_k^{n} - \Delta t_n \sum_{j=1}^{i-1} \Tilde{a}_{ij} D_1 (\rho u_1^2 )_k^{j} - \frac{\Delta t_n}{\epsilon^2} \sum_{j=1}^i a_{ij} D_{cen} (p)_k^{j} 
    \end{gather}
    with $D$ as $D_{cen}$ or $D_{upw}$, and $D_1$ as $D_{upw}$ or $D_{ES}$. Then for $\epsilon \to 0$, a solution of \eqref{AP thm full mass}, \eqref{AP thm full mom} satisfies
    \begin{gather}
        \rho_{0_k}^i \equiv \text{ constant }, \ \rho_{1_k}^i \equiv \text{ constant }, \ \rho_{0_k}^i + \epsilon \rho_{1_k}^i = \rho_0,  \\
        %\nabla \cdot u_{1_0}^i = 0, \\
        (u_{1_0})_k^{i} = (u_{1_0})_k^{n} - \Delta t_n \sum_{j=1}^{i-1} \Tilde{a}_{ij} D_1 ( u_{1_0}^2)_k^{j} - \frac{\Delta t_n}{\rho_0} \sum_{j=1}^i a_{ij} D_{cen} (p_2)_k^{j},
    \end{gather}
for all $k=1,2,\dots$, and all $i \in \{1,2,\dots,s \}$. 
\end{theorem}
\begin{proof}
    Substituting the asymptotic expansion into the momentum update equation in \eqref{AP thm full mom} and equating $\mathcal{O}\left(\frac{1}{\epsilon^2}\right)$ terms, we get
    \begin{equation}
        \sum_{j=1}^i a_{ij} D_{cen} (p_0)_k^j=0,  \text{ for all } i \in \{ 1,2,\dots,s\} \implies D_{cen} (p_0)_k^i=0, \text{ for all } i \in \{ 1,2,\dots,s\}.  
    \end{equation}
    Note that this property does not allow us to conclude that $(p_0)_k^i$ is spatially constant. Depending on boundary conditions, checkerboard modes could occur. In order to conclude that $(p_0)_k^i$ is spatially constant, it is required to consider a ghost point on the left and impose $p_{0_{ghost}}=(p_0)_{k=0}$. Then, $(p_0)_k^i$ is spatially constant and hence $(\rho_0)_k^i$ is also spatially constant since $(p_0)_k^i=\kappa \left( (\rho_0)_k^i \right)^{\gamma}$. \\
    Similarly equating $\mathcal{O}\left(\frac{1}{\epsilon}\right)$ terms in the momentum balance \eqref{AP thm full mom}, we infer that $(\rho_1)_k^i$ is spatially constant. \\
    Inserting the asymptotic expansion into the mass update equation in \eqref{AP thm full mass} and equating $\mathcal{O}(1)$ terms, we obtain,
    \begin{equation}
        \rho_{0_k}^{i}=\rho_{0_k}^n - \Delta t_n \sum_{j=1}^{i-1} a_{ij} \rho_{0_k}^j D(u_{1_0})_k^{j} - \Delta t_n a_{ii} \rho_{0_k}^n D ( u_{1_0})_k^{n} + \Delta t_n^2 a_{ii} \sum_{j=1}^{i-1} \Tilde{a}_{ij} \rho_{0_k}^j D^2 ( u_{1_0}^2 )_k^{j} + \Delta t_n^2 a_{ii} \sum_{j=1}^i a_{ij} D^2 (p_2)_k^{j}.
    \end{equation}
    Summing over all the points in the domain and using periodic boundary conditions on $\rho_2$ and $u_{1_0}$, we obtain,
    \begin{equation}
        \rho_{0_k}^{i}=\rho_{0_k}^n, \ \text{ for all } k, \ \text{ for all } i \in \{ 1,2,\dots,s\}.
    \end{equation}
    Repeating the similar procedure for $\mathcal{O}\left( \epsilon \right)$ terms of the mass update equation \eqref{AP thm full mass}, we obtain,
    \begin{equation}
        \rho_{1_k}^{i}=\rho_{1_k}^n, \ \text{ for all } k, \ \text{ for all } i \in \{ 1,2,\dots,s\}.
    \end{equation}
    Since $\rho_{{0,1}_k}^{n+1}=\rho_{{0,1}_k}^s$ due to the GSA property of IMEX-RK time discretisation, we have $\rho_{{0,1}_k}^{n+1}=\rho_{{0,1}_k}^s=\rho_{{0,1}_k}^n$. Therefore, $\rho_{{0,1}_k}^n=\rho_{{0,1}_k}^0 \equiv \text{constant}, \text{ for all } n=1,2,\dots$. \\
    Inserting this into the $\mathcal{O}(1)$ mass and momentum update equations, we get for all $i \in \{1,2,\dots,s \}$,
    \begin{gather}
    \label{AP full proof mass}
        \sum_{j=1}^{i-1} a_{ij} D (u_{1_0})_k^{j} + a_{ii} D(u_{1_0})_k^{n} - \Delta t_n a_{ii} \sum_{j=1}^{i-1} \Tilde{a}_{ij}  D^2 (u_{1_0}^2)_k^{j} - \frac{\Delta t_n}{\rho_0} a_{ii} \sum_{j=1}^i a_{ij}  D^2 (p_2)_k^{j} = 0, \\
    \label{AP full proof mom}
        (u_{1_0})_k^{i} = (u_{1_0})_k^{n} - \Delta t_n \sum_{j=1}^{i-1} \Tilde{a}_{ij} D_1 ( u_{1_0}^2)_k^{j} - \frac{\Delta t_n}{\rho_0} \sum_{j=1}^i a_{ij} D_{cen} (p_2)_k^{j}, 
    \end{gather}
    where $\rho_0 = \rho_{0_k}^0 + \epsilon \rho_{1_k}^0, \text{for any } k=1,2,\dots$.
\end{proof}
Due to the GSA property, the expressions for $\rho_{0_k}^s, \rho_{1_k}^s, (u_{1_0})_k^s$ follow for $\rho_{0_k}^{n+1}, \rho_{1_k}^{n+1}, (u_{1_0})_k^{n+1}$, for all $k=1,2,\dots$. Thus, we have devised an asymptotic preserving IMEX-RK scheme with three different types of space discretisation techniques. 

\section{Numerical results}
\label{AP_ES_Sec:NR}
In this section, we present the numerical results obtained by our asymptotic preserving IMEX-RK scheme with three different types of spatial discretisation. The numerical results include: entropy, potential energy (PE), kinetic energy (KE) plots and accuracy tables of a standard periodic problem for different values of $\epsilon$; entropy, density and momentum plots for colliding acoustic waves and Riemann problems; entropy, PE and KE plots for Gresho and travelling vortex problems.  

\subsection{Standard periodic problem}
The domain of the problem is $\Omega:=[0,1]$, and the initial conditions are: 
\begin{gather}
    \rho^0(x)=1+\epsilon^2 \sin \left(2\pi x\right) \\
    u_{1}^0(x)=1+\epsilon \sin \left( 2\pi x\right).
\end{gather}
The parametric values are: $\kappa=1$ and $\gamma=2$. The entropy plots and accuracy tables of this problem will be presented for different values of $\epsilon$.  
\subsubsection{Entropy, kinetic energy (KE) and potential energy (PE)}
\label{subsubsec:Num_Ent}
\begin{figure}[h!]
\centering
\begin{subfigure}[b]{0.32\textwidth}
\centering
\includegraphics[width=\textwidth]{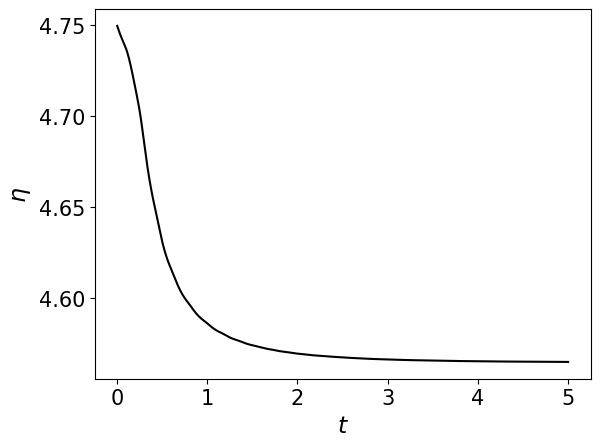}
\caption{Entropy - Type $1$}
\label{Fig:Entropytype1eps05}
\end{subfigure}
\hspace{-0.2cm}
\begin{subfigure}[b]{0.32\textwidth}
\centering
\includegraphics[width=\textwidth]{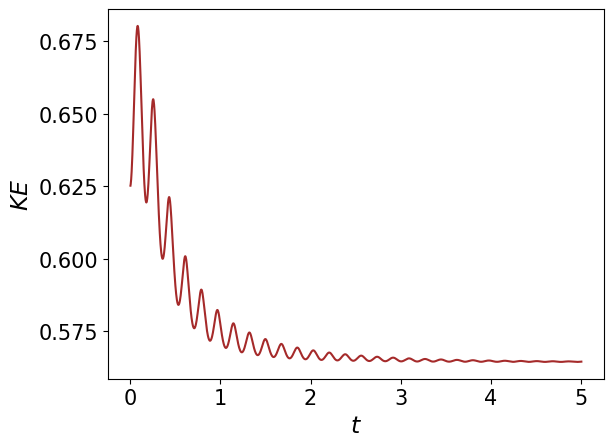}
\caption{KE - Type $1$}
\label{Fig:KEtype1eps05}
\end{subfigure}
\hspace{-0.2cm}
\begin{subfigure}[b]{0.32\textwidth}
\centering
\includegraphics[width=\textwidth]{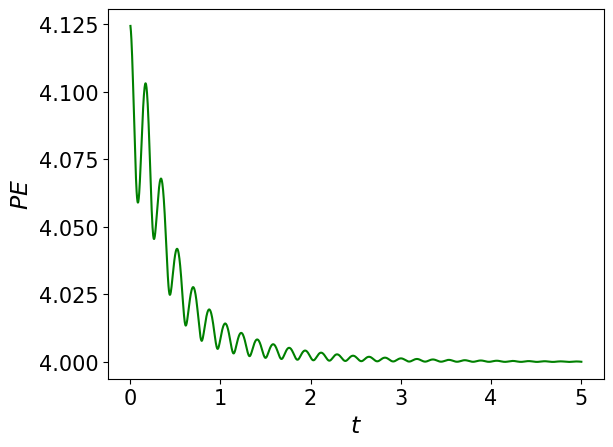}
\caption{PE - Type $1$}
\label{Fig:PEtype1eps05}
\end{subfigure}
\vfill
\begin{subfigure}[b]{0.32\textwidth}
\centering
\includegraphics[width=\textwidth]{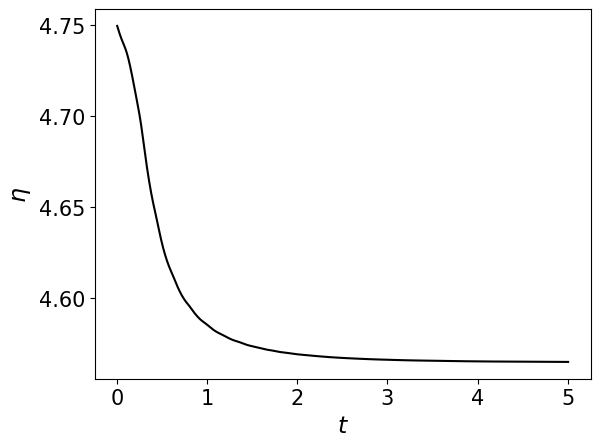}
\caption{Entropy - Type $2$}
\label{Fig:Entropytype2eps05}
\end{subfigure}
\hspace{-0.2cm}
\begin{subfigure}[b]{0.32\textwidth}
\centering
\includegraphics[width=\textwidth]{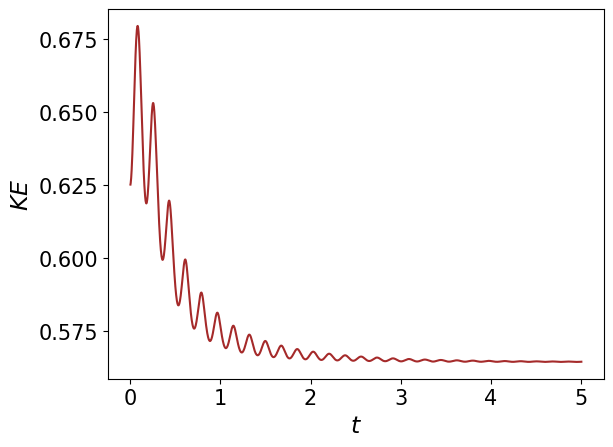}
\caption{KE - Type $2$}
\label{Fig:KEtype2eps05}
\end{subfigure}
\hspace{-0.2cm}
\begin{subfigure}[b]{0.32\textwidth}
\centering
\includegraphics[width=\textwidth]{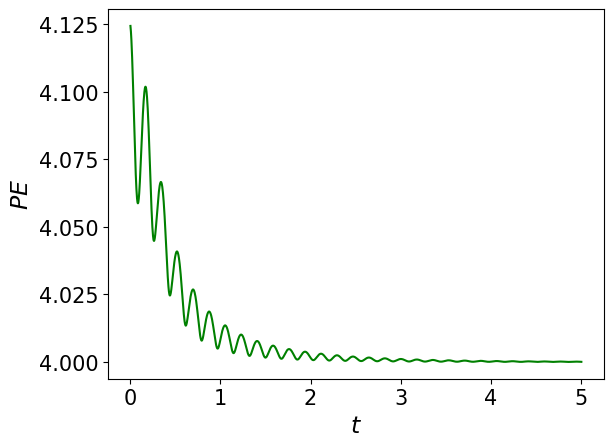}
\caption{PE - Type $2$}
\label{Fig:PEtype2eps05}
\end{subfigure}
\vfill
\begin{subfigure}[b]{0.32\textwidth}
\centering
\includegraphics[width=\textwidth]{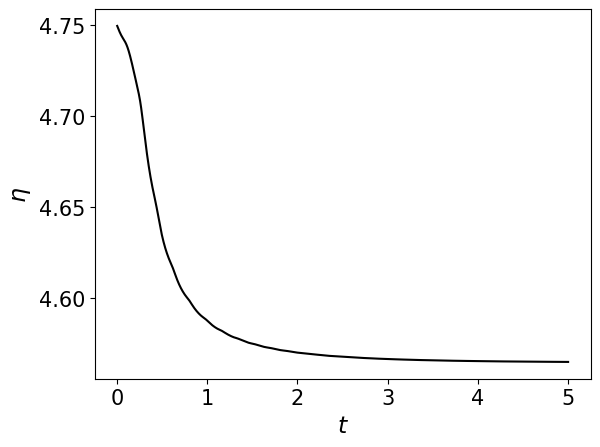}
\caption{Entropy - Type $3$}
\label{Fig:Entropytype3eps05}
\end{subfigure}
\hspace{-0.2cm}
\begin{subfigure}[b]{0.32\textwidth}
\centering
\includegraphics[width=\textwidth]{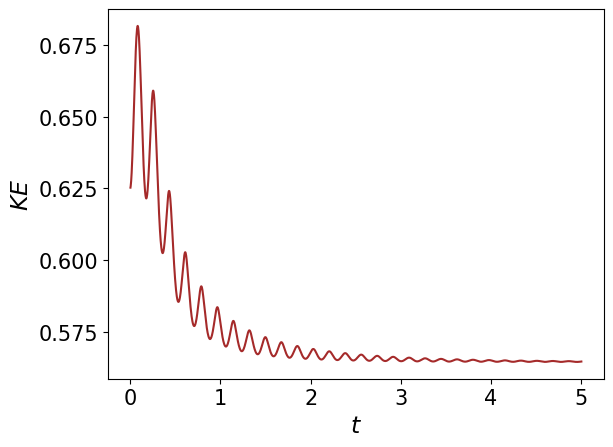}
\caption{KE - Type $3$}
\label{Fig:KEtype3eps05}
\end{subfigure}
\hspace{-0.2cm}
\begin{subfigure}[b]{0.32\textwidth}
\centering
\includegraphics[width=\textwidth]{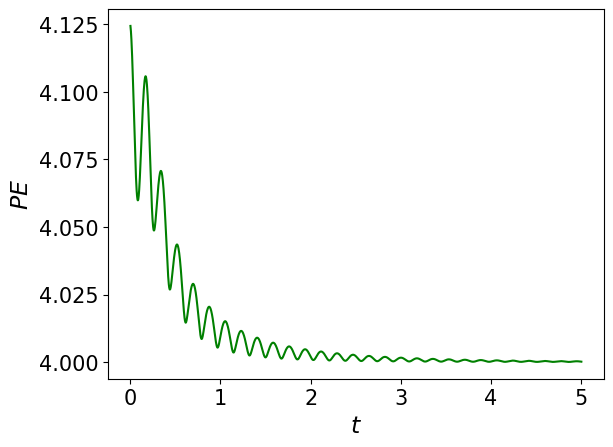}
\caption{PE - Type $3$}
\label{Fig:PEtype3eps05}
\end{subfigure}
\caption{\centering \textbf{Standard periodic problem:} Entropy, KE and PE plots for $\epsilon=0.5$ using space discretisation types 1, 2 and 3.}
\label{Fig:eps05}
\end{figure}

\begin{figure}[h!]
\centering
\begin{subfigure}[b]{0.32\textwidth}
\centering
\includegraphics[width=\textwidth]{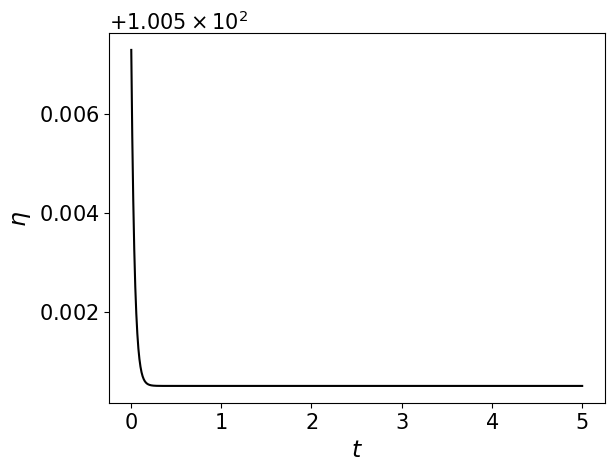}
\caption{Entropy - Type $1$}
\label{Fig:Entropytype1eps01}
\end{subfigure}
\hspace{-0.2cm}
\begin{subfigure}[b]{0.32\textwidth}
\centering
\includegraphics[width=\textwidth]{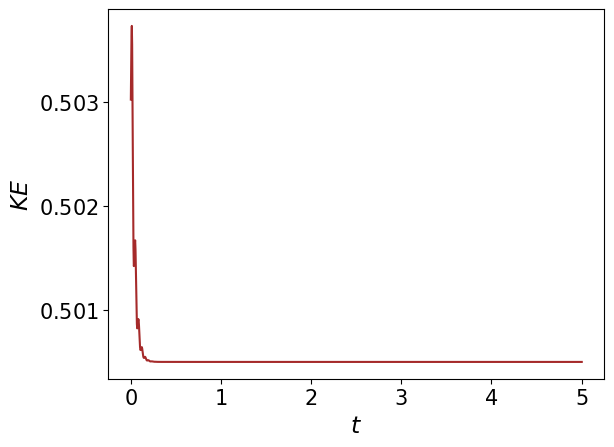}
\caption{KE - Type $1$}
\label{Fig:KEtype1eps01}
\end{subfigure}
\hspace{-0.2cm}
\begin{subfigure}[b]{0.32\textwidth}
\centering
\includegraphics[width=\textwidth]{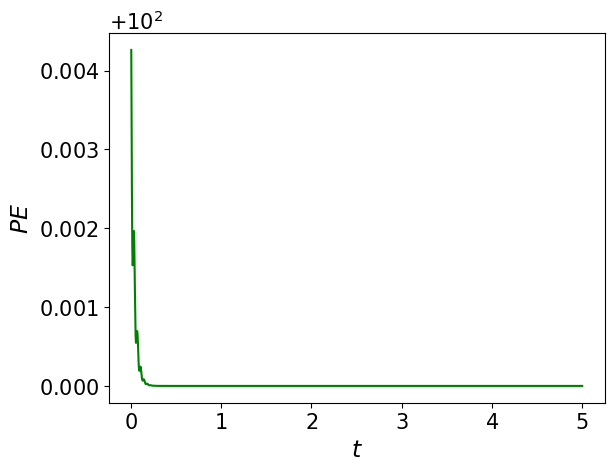}
\caption{PE - Type $1$}
\label{Fig:PEtype1eps01}
\end{subfigure}
\vfill
\begin{subfigure}[b]{0.32\textwidth}
\centering
\includegraphics[width=\textwidth]{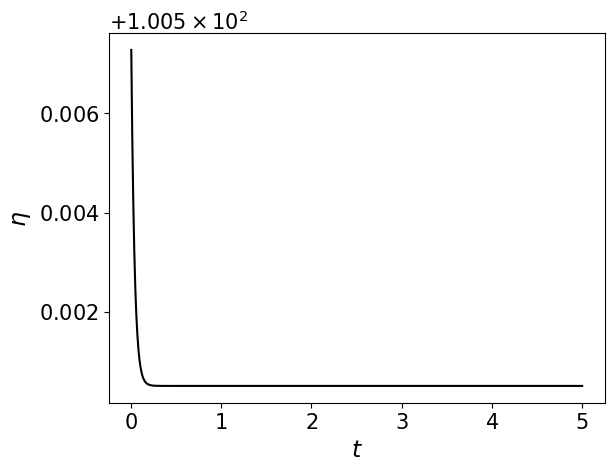}
\caption{Entropy - Type $2$}
\label{Fig:Entropytype2eps01}
\end{subfigure}
\hspace{-0.2cm}
\begin{subfigure}[b]{0.32\textwidth}
\centering
\includegraphics[width=\textwidth]{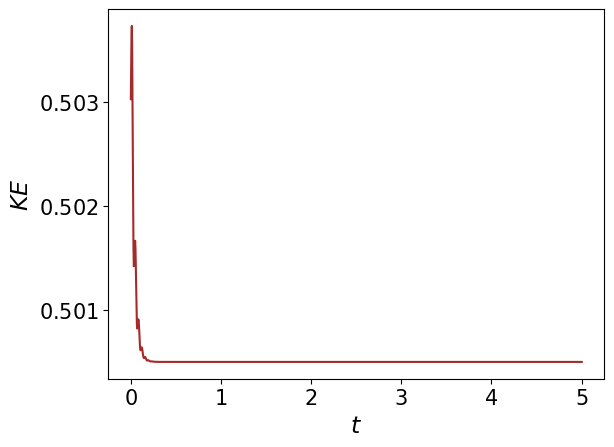}
\caption{KE - Type $2$}
\label{Fig:KEtype2eps01}
\end{subfigure}
\hspace{-0.2cm}
\begin{subfigure}[b]{0.32\textwidth}
\centering
\includegraphics[width=\textwidth]{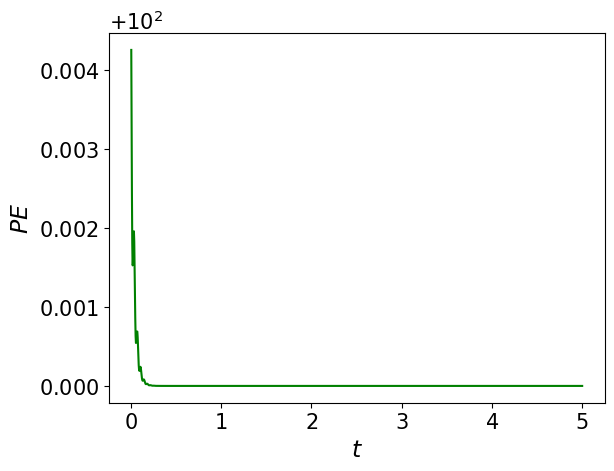}
\caption{PE - Type $2$}
\label{Fig:PEtype2eps01}
\end{subfigure}
\vfill
\begin{subfigure}[b]{0.32\textwidth}
\centering
\includegraphics[width=\textwidth]{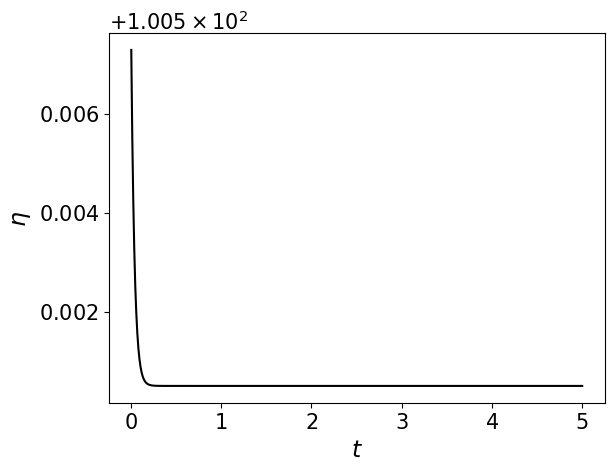}
\caption{Entropy - Type $3$}
\label{Fig:Entropytype3eps01}
\end{subfigure}
\hspace{-0.2cm}
\begin{subfigure}[b]{0.32\textwidth}
\centering
\includegraphics[width=\textwidth]{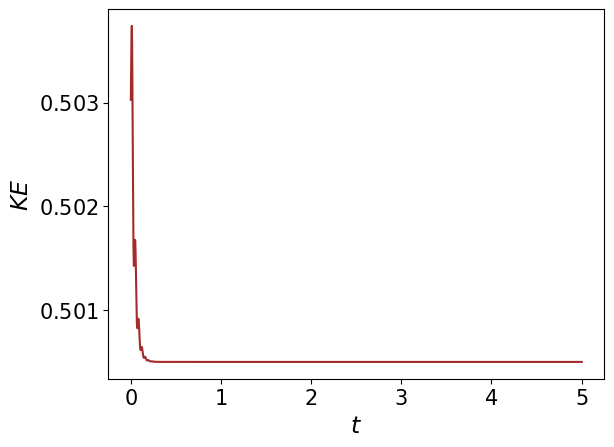}
\caption{KE - Type $3$}
\label{Fig:KEtype3eps01}
\end{subfigure}
\hspace{-0.2cm}
\begin{subfigure}[b]{0.32\textwidth}
\centering
\includegraphics[width=\textwidth]{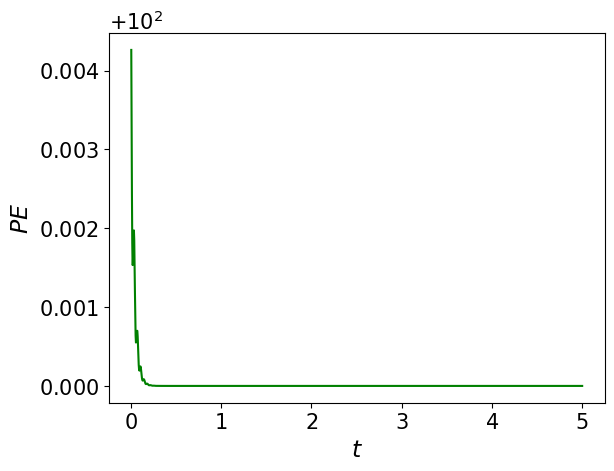}
\caption{PE - Type $3$}
\label{Fig:PEtype3eps01}
\end{subfigure}
\caption{\centering \textbf{Standard periodic problem:} Entropy, KE and PE plots for $\epsilon=0.1$ using space discretisation types 1, 2 and 3.}
\label{Fig:eps01}
\end{figure}

\begin{figure}[h!]
\centering
\begin{subfigure}[b]{0.32\textwidth}
\centering
\includegraphics[width=\textwidth]{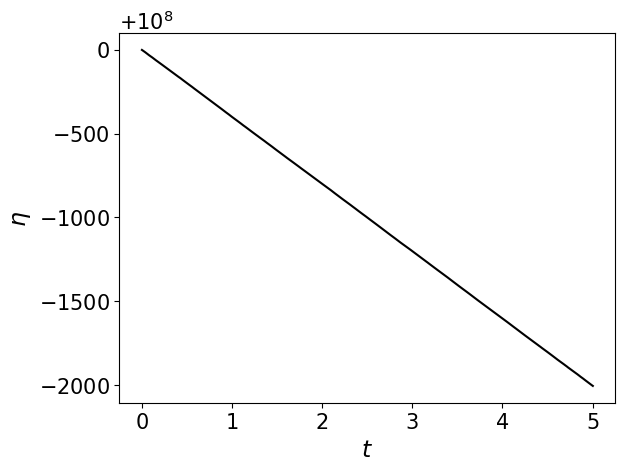}
\caption{Entropy - Type $1$}
\label{Fig:Entropytype1eps00001}
\end{subfigure}
\hspace{-0.2cm}
\begin{subfigure}[b]{0.29\textwidth}
\centering
\includegraphics[width=\textwidth]{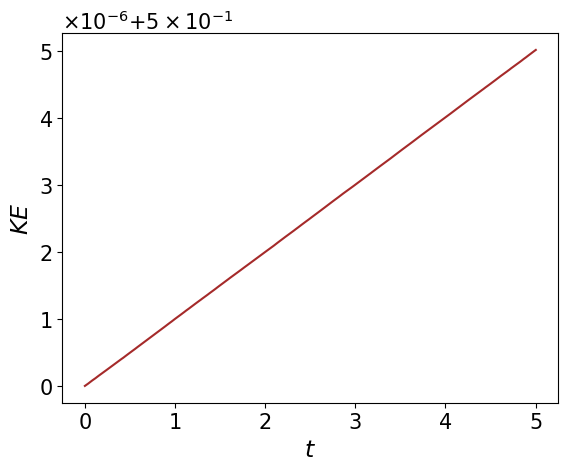}
\caption{KE - Type $1$}
\label{Fig:KEtype1eps00001}
\end{subfigure}
\hspace{-0.2cm}
\begin{subfigure}[b]{0.32\textwidth}
\centering
\includegraphics[width=\textwidth]{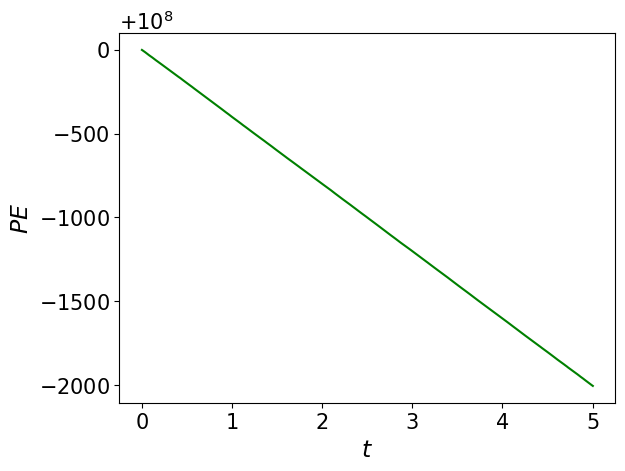}
\caption{PE - Type $1$}
\label{Fig:PEtype1eps00001}
\end{subfigure}
\vfill
\begin{subfigure}[b]{0.32\textwidth}
\centering
\includegraphics[width=\textwidth]{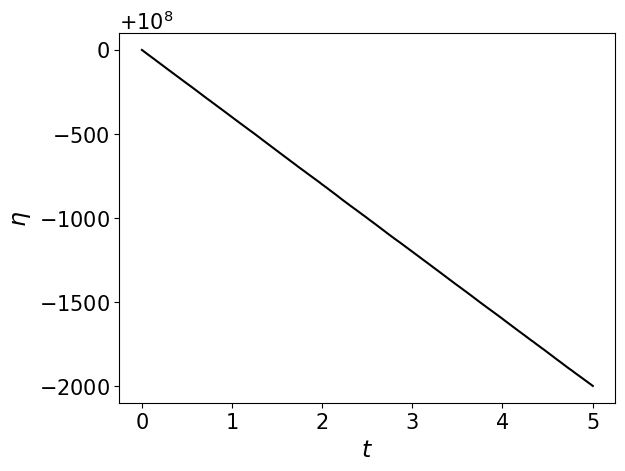}
\caption{Entropy - Type $2$}
\label{Fig:Entropytype2eps00001}
\end{subfigure}
\hspace{-0.2cm}
\begin{subfigure}[b]{0.29\textwidth}
\centering
\includegraphics[width=\textwidth]{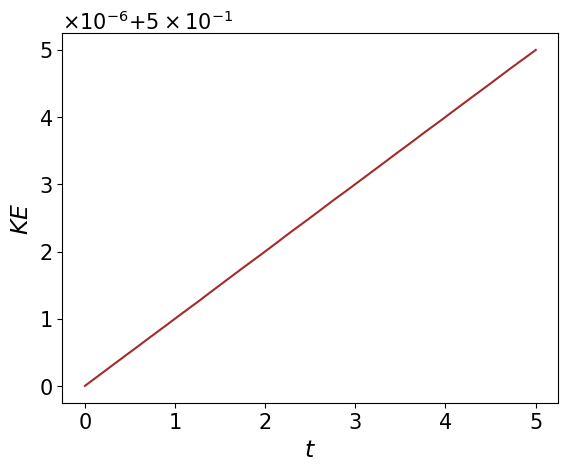}
\caption{KE - Type $2$}
\label{Fig:KEtype2eps00001}
\end{subfigure}
\hspace{-0.2cm}
\begin{subfigure}[b]{0.32\textwidth}
\centering
\includegraphics[width=\textwidth]{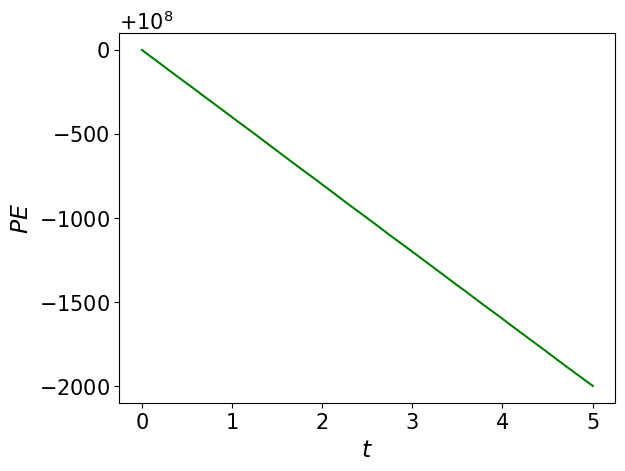}
\caption{PE - Type $2$}
\label{Fig:PEtype2eps00001}
\end{subfigure}
\vfill
\begin{subfigure}[b]{0.32\textwidth}
\centering
\includegraphics[width=\textwidth]{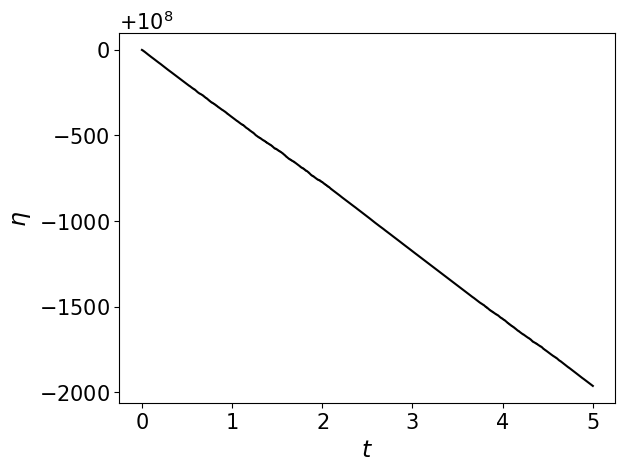}
\caption{Entropy - Type $3$}
\label{Fig:Entropytype3eps00001}
\end{subfigure}
\hspace{-0.2cm}
\begin{subfigure}[b]{0.31\textwidth}
\centering
\includegraphics[width=\textwidth]{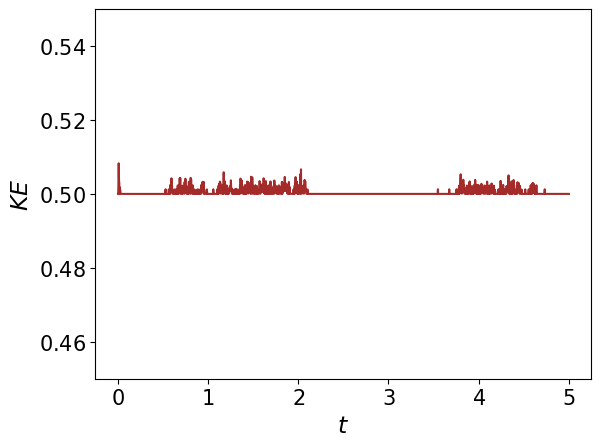}
\caption{KE - Type $3$}
\label{Fig:KEtype3eps00001}
\end{subfigure}
\hspace{-0.2cm}
\begin{subfigure}[b]{0.32\textwidth}
\centering
\includegraphics[width=\textwidth]{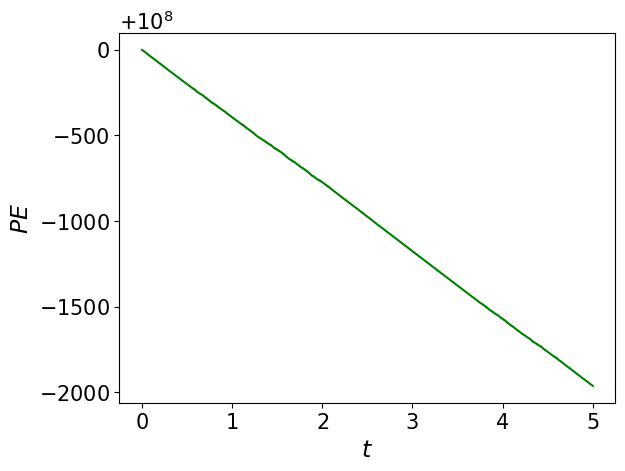}
\caption{PE - Type $3$}
\label{Fig:PEtype3eps00001}
\end{subfigure}
\caption{\centering \textbf{Standard periodic problem:} Entropy, KE and PE plots for $\epsilon=0.0001$ using space discretisation types 1, 2 and 3.}
\label{Fig:eps00001}
\end{figure}

The domain $\Omega$ is discretised into $N=200$ grid points. The first order IMEX scheme $ARS(1,1,1)$ is used along with three types of spatial discretisation techniques to obtain the plots on entropy, KE and PE. The time step is chosen as: 
\begin{equation}
\label{Num res time step}
    \Delta t_n = C \frac{\Delta x}{\max\limits_{i \in \Omega_N} \left(\left|u_{1_i}^n\right|\right)},
\end{equation}
where $C$ is the CFL number, and $\Omega_N$ is the discretised set of the domain $\Omega$ (that is, $\Omega_N=\{1,2,\dots,N\}$). \\
Figures \ref{Fig:eps05}, \ref{Fig:eps01} and \ref{Fig:eps00001} show the entropy, KE and PE plots obtained for $\epsilon=0.5,0.1$ and $10^{-4}$ respectively by using type 1, type 2 and type 3 spatial discretisation techniques. The global value of the convex entropy function at time $t_n$ is given by
\begin{equation}
    \eta^n = KE^n + PE^n, \text{ where } KE^n = \sum_{k=1}^N  \frac{1}{2}\rho_k^n \left(u_{1_k}^n\right)^2 \Delta x_k, \ PE^n = \sum_{k=1}^N \frac{1}{\epsilon^2} \frac{p_k^n}{\gamma-1} \Delta x_k
\end{equation}
with $\Delta x_k = \Delta x = L/N$ for $k \in \{ 1,2,\dots,N\}$, and $L$ is the length of the domain. The plots in Figures \ref{Fig:eps05}, \ref{Fig:eps01} and \ref{Fig:eps00001} are obtained at time $T=5$ to depict the long time behaviours of entropy, KE and PE. \\
$\mathbf{\epsilon=0.5:}$ From Figure \ref{Fig:eps05}, we observe that entropy which is the sum of kinetic and potential energies decays while kinetic and potential energies separately do not decay. It is expected as we only maintain entropy stability in the space discretisation. Such an entropy decay is observed for $C=0.1,0.2,\dots,0.9$ (results corresponding to $C=0.8$ are shown in Figure \ref{Fig:eps05}). The results for type 3 use $q=0$, meaning that entropy conserving spacial discretisation is used for $\nabla \cdot (\rho \mathbf{u} \otimes \mathbf{u})$. \\
$\mathbf{\epsilon=0.1:}$ Similar to $\epsilon=0.5$, we observe for $\epsilon=0.1$ (from Figure \ref{Fig:eps01}) that entropy decays while potential and kinetic energies do not. Such entropy decay is observed for $C=0.1,0.2,\dots,0.9$ (results corresponding to $C=0.8$ are shown in Figure \ref{Fig:eps01}). The results for type 3 are obtained with $q=0$. \\
$\mathbf{\epsilon=0.0001:}$ Unlike the previous cases, for $\epsilon=0.0001$, the entropy decay is very sensitive to the value of $C$. This may be due to the entropy production from the explicit time discretisation of the nonlinear terms. The results in Figure \ref{Fig:eps00001} are obtained with $C=0.5$. Further, the results for type 3 use $q=2$, meaning that entropy stable spatial discretisation is used for $\nabla \cdot (\rho \mathbf{u} \otimes \mathbf{u})$.  

\subsubsection{Order of accuracy}
In this subsection, we show the order of accuracy of type 2 spatial discretisation paired with $ARS(1,1,1)$ IMEX time discretisation that is first order accurate in time. In \cref{tab:EOC_rho,tab:EOC_u}, we observe the order of accuracy of density $\rho$ and velocity $u_1$ by using different number of grid points: $N=20,25,100,200,500$. The reference solution is obtained with $N=1000$. The time step is chosen according to \eqref{Num res time step} with CFL number, $C=0.5$. \\ %For $\epsilon=0.5,0.1$, a type 3 spatial discretisation with $q=0$ is used (meaning an entropy conserving spatial discretisation for $\nabla \cdot (\rho \mathbf{u} \otimes \mathbf{u})$), whereas for $\epsilon=0.0001$, a type 3 spatial discretisation with $q=2$ is used (meaning an entropy stable spatial discretisation for $\nabla \cdot (\rho \mathbf{u} \otimes \mathbf{u})$).
For $\epsilon=0.5, 0.1$, both density and velocity converge with a rate of about $1$. For $\epsilon=10^{-4}$, density shows very small errors of $\mathcal{O}(10^{-7})$ for all tested values of $N$. This indicates that $\rho_0$ and $\rho_1$ in asymptotic expansion \eqref{AP ansatz 1} are constant, and thus our scheme is asymptotic preserving. Further $\rho_2$ in asymptotic expansion \eqref{AP ansatz 1} and the velocity $u_1$ are converging with a rate of about $0.5$. 

\begin{table}[h!]
    \centering
    \renewcommand{\arraystretch}{1.3} 
    \setlength{\tabcolsep}{10pt}      
    
    \begin{tabular}{|c|c|c|c|c|c|c|c|}
        \hline
        \multirow{2}{*}{\textbf{N}} & \multirow{2}{*}{$\mathbf{\Delta x}$} & 
        \multicolumn{2}{c|}{$\mathbf{\epsilon=0.5}$} & 
        \multicolumn{2}{c|}{$\mathbf{\epsilon=0.1}$} & 
        \multicolumn{2}{c|}{$\mathbf{\epsilon=10^{-4}}$} \\ 
        \cline{3-8}
        & & $\mathbf{||\rho \textbf{ error}||_{L_2}}$ & \textbf{EOC} 
        & $\mathbf{||\rho \textbf{ error}||_{L_2}}$ & \textbf{EOC} 
        & $\mathbf{||\rho \textbf{ error}||_{L_2}}$ & \textbf{EOC} \\ 
        \hline
        20  & 0.05  & 0.03267  & -      & 0.00447  & -      & 4.89 $\times 10^{-7}$ & -      \\  
        50  & 0.02  & 0.01644  & 0.7497 & 0.00370  & 0.2069 & 4.77 $\times 10^{-7}$ & 0.0266 \\  
        100 & 0.01  & 0.01006  & 0.7083 & 0.00256  & 0.5315 & 4.52 $\times 10^{-7}$ & 0.0773 \\  
        250 & 0.004 & 0.00282  & 1.3874 & 0.00117  & 0.8510 & 3.76 $\times 10^{-7}$ & 0.1999 \\  
        500 & 0.002 & 0.00156  & 0.8580 & 0.00044  & 1.4177 & 2.49 $\times 10^{-7}$ & 0.5962 \\  
        \hline
    \end{tabular}
    
    \caption{\centering \textbf{Standard periodic problem:} Convergence rates of $L_2$ error in $\rho$ using $ARS(1,1,1)$ coupled with type 2 discretisation.}
    \label{tab:EOC_rho}
\end{table}

\begin{table}[h!]
    \centering
    \renewcommand{\arraystretch}{1.3} 
    \setlength{\tabcolsep}{10pt}      
    
    \begin{tabular}{|c|c|c|c|c|c|c|c|}
        \hline
        \multirow{2}{*}{\textbf{N}} & \multirow{2}{*}{$\mathbf{\Delta x}$} & 
        \multicolumn{2}{c|}{$\mathbf{\epsilon=0.5}$} & 
        \multicolumn{2}{c|}{$\mathbf{\epsilon=0.1}$} & 
        \multicolumn{2}{c|}{$\mathbf{\epsilon=10^{-4}}$} \\ 
        \cline{3-8}
        & & $\mathbf{||u_1 \textbf{ error}||_{L_2}}$ & \textbf{EOC} 
        & $\mathbf{||u_1 \textbf{ error}||_{L_2}}$ & \textbf{EOC} 
        & $\mathbf{||u_1 \textbf{ error}||_{L_2}}$ & \textbf{EOC} \\ 
        \hline
        20  & 0.05  & 0.12749  & -      & 0.08749  & -      & 5.27 $\times 10^{-7}$ & -      \\  
        50  & 0.02  & 0.02794  & 0.0021 & 0.07201  & 0.2125 & 5.04 $\times 10^{-7}$ & 0.0493 \\  
        100 & 0.01  & 0.00628  & 2.1540 & 0.05035  & 0.5163 & 4.65 $\times 10^{-7}$ & 0.1158 \\  
        250 & 0.004 & 0.00351  & 0.6346 & 0.02121  & 0.9433 & 4.08 $\times 10^{-7}$ & 0.1425 \\  
        500 & 0.002 & 0.00127  & 1.4661 & 0.00761  & 1.4784 & 2.94 $\times 10^{-7}$ & 0.4730 \\  
        \hline
    \end{tabular}
    
    \caption{\centering \textbf{Standard periodic problem:} Convergence rates of $L_2$ error in $u_1$ using $ARS(1,1,1)$ coupled with type 2 discretisation.}
    \label{tab:EOC_u}
\end{table}

\subsection{Colliding acoustic waves problem}
\begin{figure}[h!]
\centering
\begin{subfigure}[b]{0.32\textwidth}
\centering
\includegraphics[width=\textwidth]{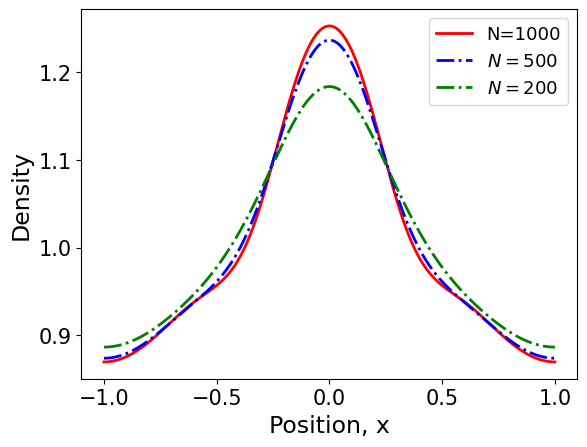}
\label{Fig:CAWtype1_den_t004}
\caption{Density at $T=0.04$}
\end{subfigure}
\hspace{-0.2cm}
\begin{subfigure}[b]{0.32\textwidth}
\centering
\includegraphics[width=\textwidth]{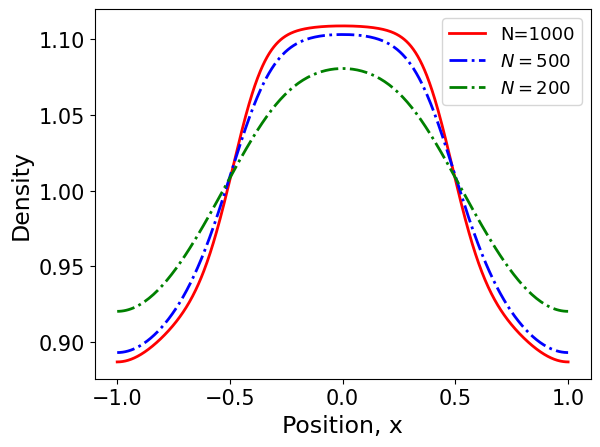}
\label{Fig:CAWtype1_den_t006}
\caption{Density at $T=0.06$}
\end{subfigure}
\hspace{-0.2cm}
\begin{subfigure}[b]{0.32\textwidth}
\centering
\includegraphics[width=\textwidth]{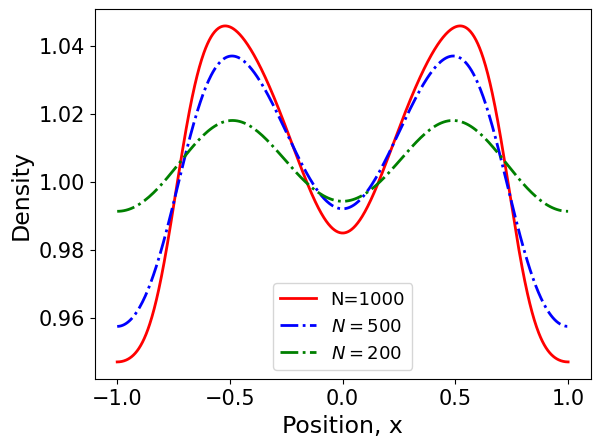}
\label{Fig:CAWtype1_den_t008}
\caption{Density at $T=0.08$}
\end{subfigure}
\vfill
\begin{subfigure}[b]{0.33\textwidth}
\centering
\includegraphics[width=\textwidth]{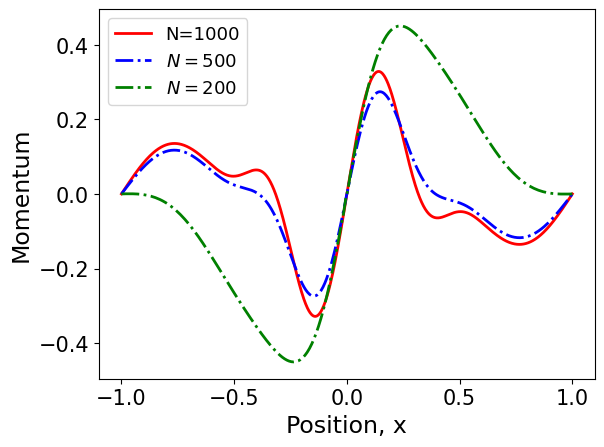}
\label{Fig:CAWtype1_mom_t004}
\caption{Momentum at $T=0.04$}
\end{subfigure}
\hspace{-0.2cm}
\begin{subfigure}[b]{0.32\textwidth}
\centering
\includegraphics[width=\textwidth]{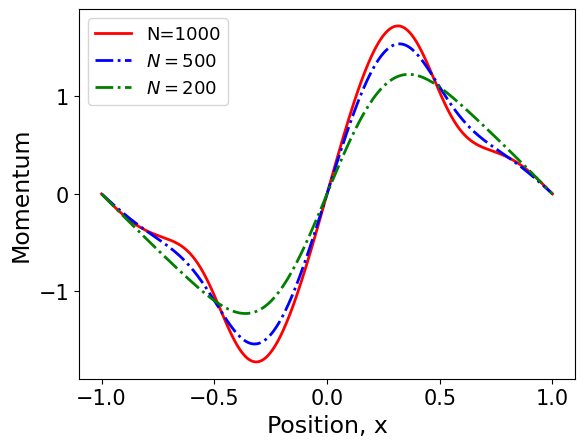}
\label{Fig:CAWtype1_mom_t006}
\caption{Momentum at $T=0.06$}
\end{subfigure}
\hfill
\begin{subfigure}[b]{0.32\textwidth}
\centering
\includegraphics[width=\textwidth]{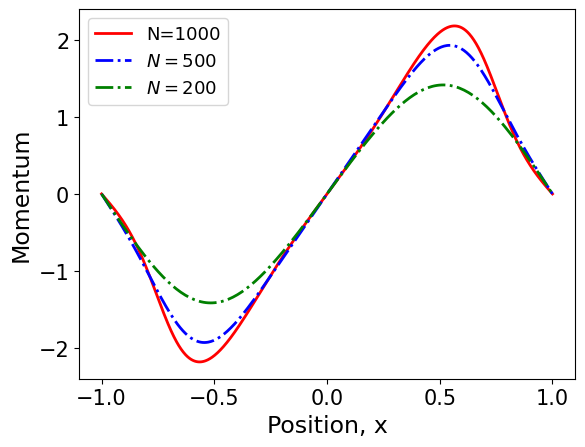}
\label{Fig:CAWtype1_mom_t008}
\caption{Momentum at $T=0.08$}
\end{subfigure}
\vfill
\begin{subfigure}[b]{0.32\textwidth}
\centering
\includegraphics[width=\textwidth]{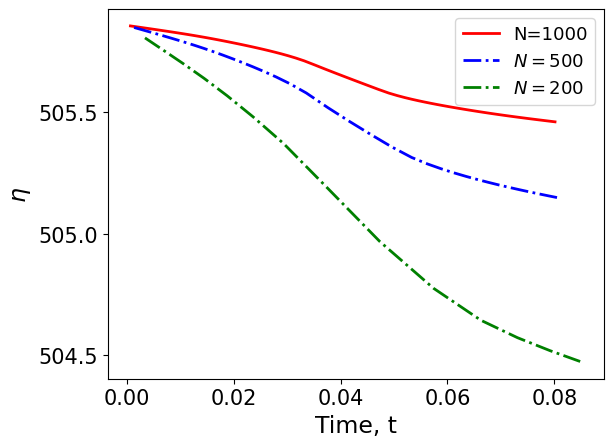}
\label{Fig:CAWtype1_ent}
\caption{Entropy}
\end{subfigure}
\hspace{-0.2cm}
\begin{subfigure}[b]{0.31\textwidth}
\centering
\includegraphics[width=\textwidth]{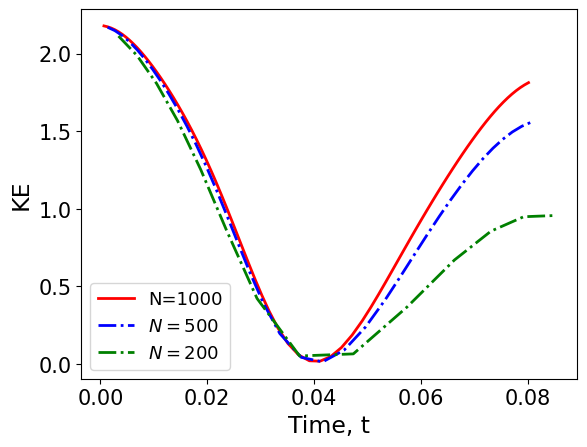}
\label{Fig:CAWtype1_KE}
\caption{Kinetic energy}
\end{subfigure}
\hspace{-0.2cm}
\begin{subfigure}[b]{0.32\textwidth}
\centering
\includegraphics[width=\textwidth]{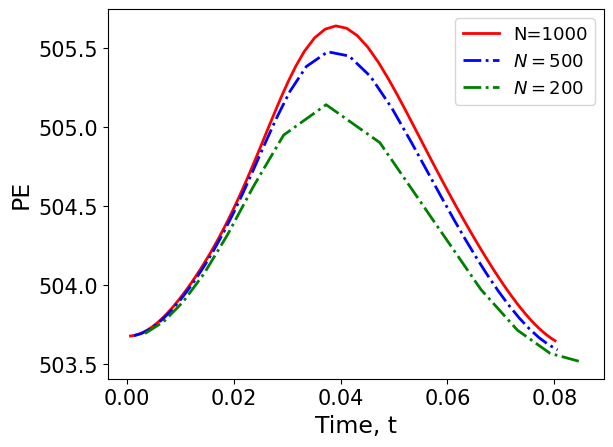}
\label{Fig:CAWtype1_PE}
\caption{Potential energy}
\end{subfigure}
\caption{\centering \textbf{Colliding acoustic waves problem} with $\epsilon=0.1$ using type 1 space discretisation.}
\label{Fig:CAWtype1}
\end{figure}

\begin{figure}[h!]
\centering
\begin{subfigure}[b]{0.32\textwidth}
\centering
\includegraphics[width=\textwidth]{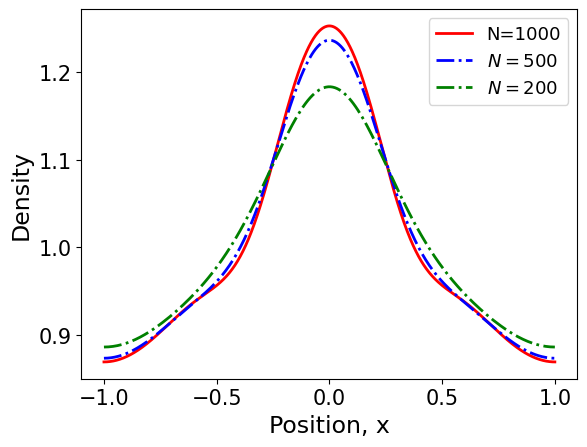}
\label{Fig:CAWtype2_den_t004}
\caption{Density at $T=0.04$}
\end{subfigure}
\hspace{-0.2cm}
\begin{subfigure}[b]{0.32\textwidth}
\centering
\includegraphics[width=\textwidth]{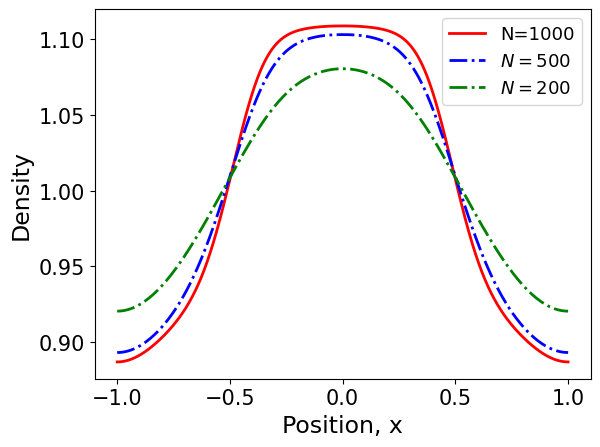}
\label{Fig:CAWtype2_den_t006}
\caption{Density at $T=0.06$}
\end{subfigure}
\hspace{-0.2cm}
\begin{subfigure}[b]{0.32\textwidth}
\centering
\includegraphics[width=\textwidth]{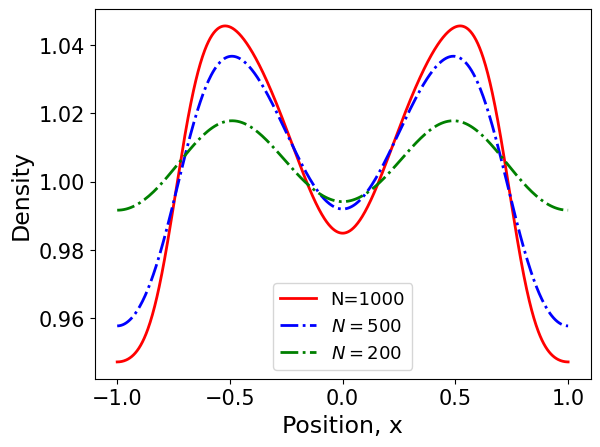}
\label{Fig:CAWtype2_den_t008}
\caption{Density at $T=0.08$}
\end{subfigure}
\vfill
\begin{subfigure}[b]{0.33\textwidth}
\centering
\includegraphics[width=\textwidth]{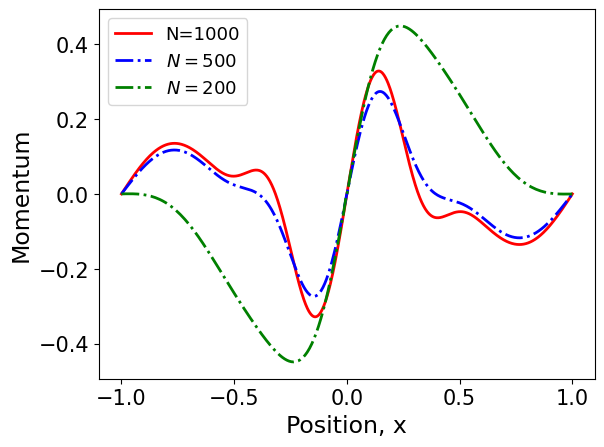}
\label{Fig:CAWtype2_mom_t004}
\caption{Momentum at $T=0.04$}
\end{subfigure}
\hspace{-0.2cm}
\begin{subfigure}[b]{0.32\textwidth}
\centering
\includegraphics[width=\textwidth]{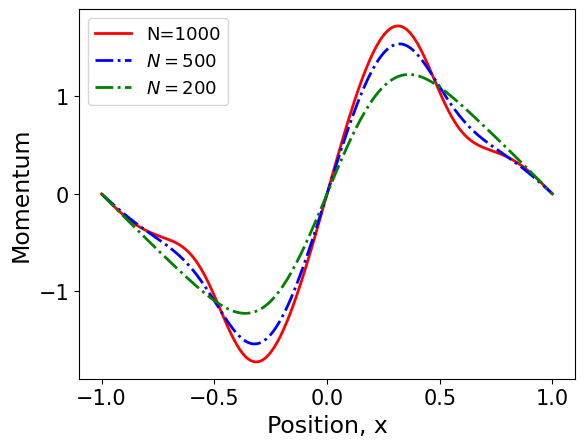}
\label{Fig:CAWtype2_mom_t006}
\caption{Momentum at $T=0.06$}
\end{subfigure}
\hfill
\begin{subfigure}[b]{0.32\textwidth}
\centering
\includegraphics[width=\textwidth]{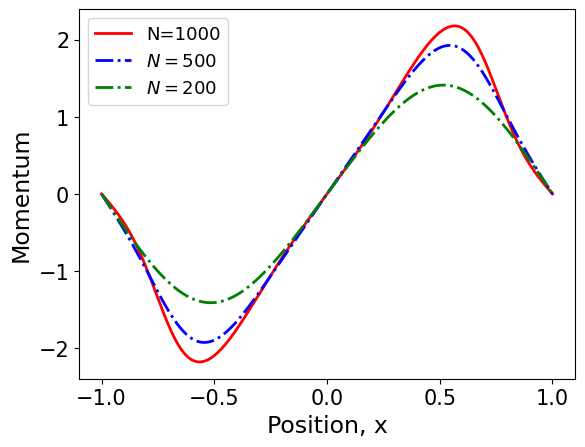}
\label{Fig:CAWtype2_mom_t008}
\caption{Momentum at $T=0.08$}
\end{subfigure}
\vfill
\begin{subfigure}[b]{0.32\textwidth}
\centering
\includegraphics[width=\textwidth]{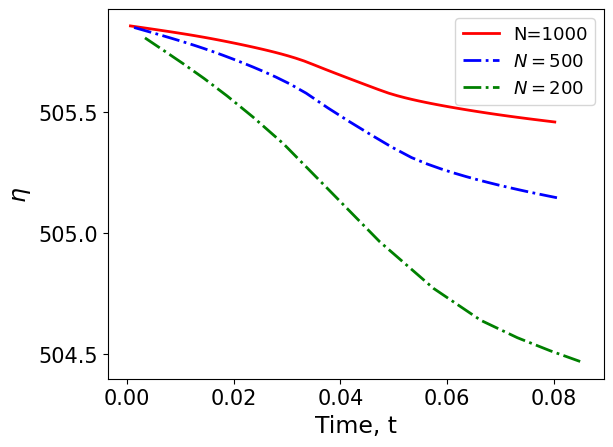}
\label{Fig:CAWtype2_ent}
\caption{Entropy}
\end{subfigure}
\hspace{-0.2cm}
\begin{subfigure}[b]{0.31\textwidth}
\centering
\includegraphics[width=\textwidth]{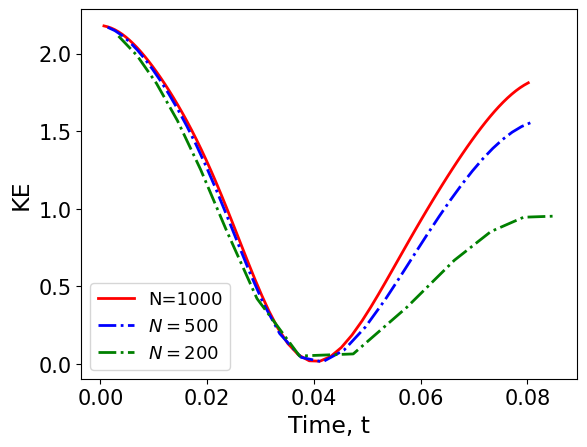}
\label{Fig:CAWtype2_KE}
\caption{Kinetic energy}
\end{subfigure}
\hspace{-0.2cm}
\begin{subfigure}[b]{0.32\textwidth}
\centering
\includegraphics[width=\textwidth]{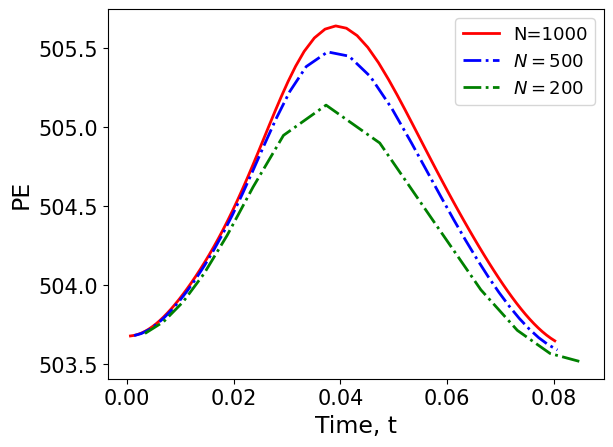}
\label{Fig:CAWtype2_PE}
\caption{Potential energy}
\end{subfigure}
\caption{\centering \textbf{Colliding acoustic waves problem} with $\epsilon=0.1$ using type 2 space discretisation.}
\label{Fig:CAWtype2}
\end{figure}

\begin{figure}[h!]
\centering
\begin{subfigure}[b]{0.32\textwidth}
\centering
\includegraphics[width=\textwidth]{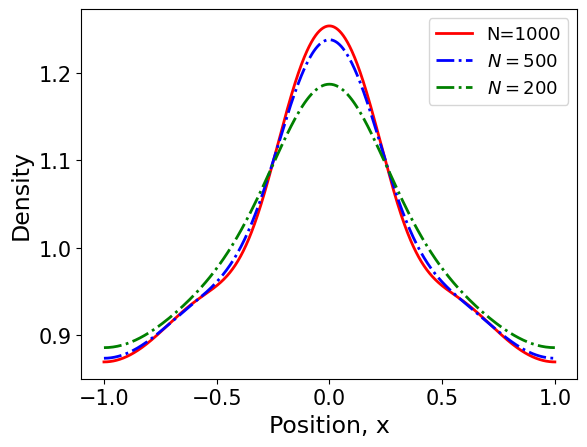}
\label{Fig:CAWtype3_den_t004}
\caption{Density at $T=0.04$}
\end{subfigure}
\hspace{-0.2cm}
\begin{subfigure}[b]{0.32\textwidth}
\centering
\includegraphics[width=\textwidth]{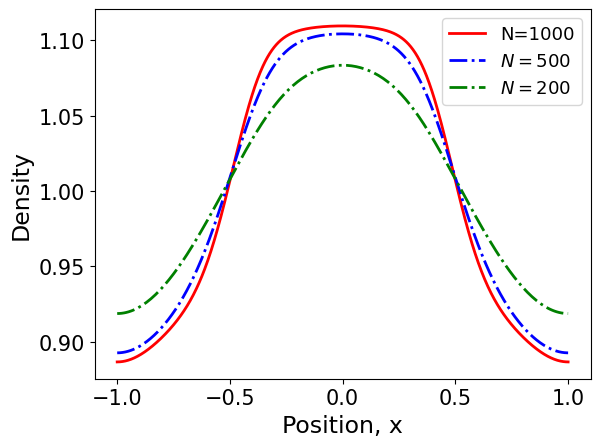}
\label{Fig:CAWtype3_den_t006}
\caption{Density at $T=0.06$}
\end{subfigure}
\hspace{-0.2cm}
\begin{subfigure}[b]{0.32\textwidth}
\centering
\includegraphics[width=\textwidth]{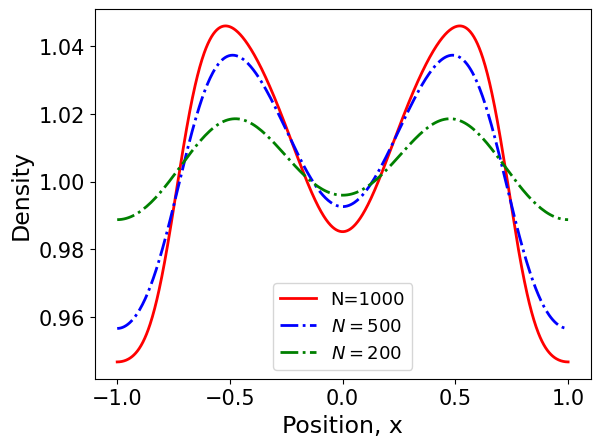}
\label{Fig:CAWtype3_den_t008}
\caption{Density at $T=0.08$}
\end{subfigure}
\vfill
\begin{subfigure}[b]{0.33\textwidth}
\centering
\includegraphics[width=\textwidth]{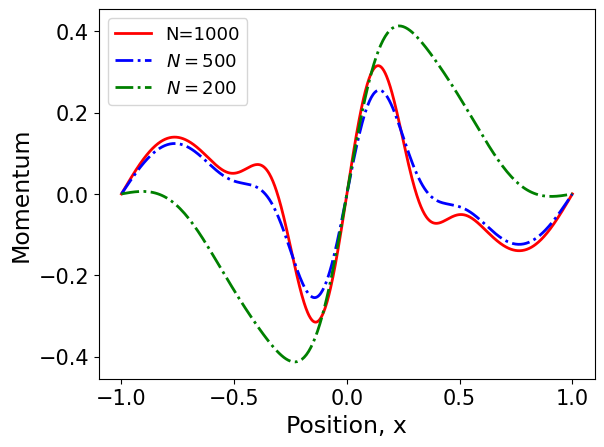}
\label{Fig:CAWtype3_mom_t004}
\caption{Momentum at $T=0.04$}
\end{subfigure}
\hspace{-0.2cm}
\begin{subfigure}[b]{0.32\textwidth}
\centering
\includegraphics[width=\textwidth]{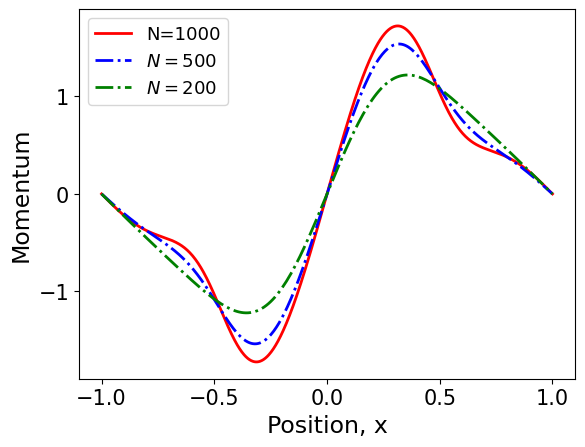}
\label{Fig:CAWtype3_mom_t006}
\caption{Momentum at $T=0.06$}
\end{subfigure}
\hfill
\begin{subfigure}[b]{0.32\textwidth}
\centering
\includegraphics[width=\textwidth]{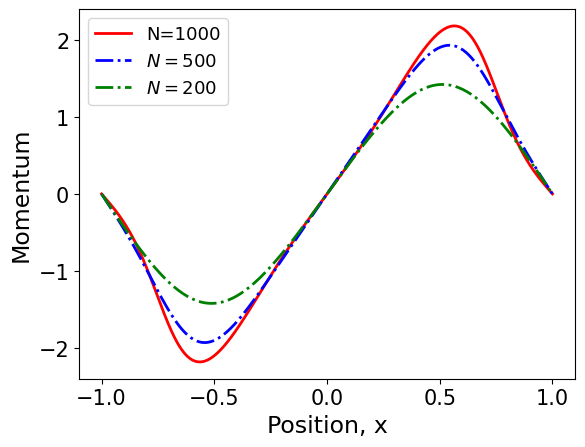}
\label{Fig:CAWtype3_mom_t008}
\caption{Momentum at $T=0.08$}
\end{subfigure}
\vfill
\begin{subfigure}[b]{0.32\textwidth}
\centering
\includegraphics[width=\textwidth]{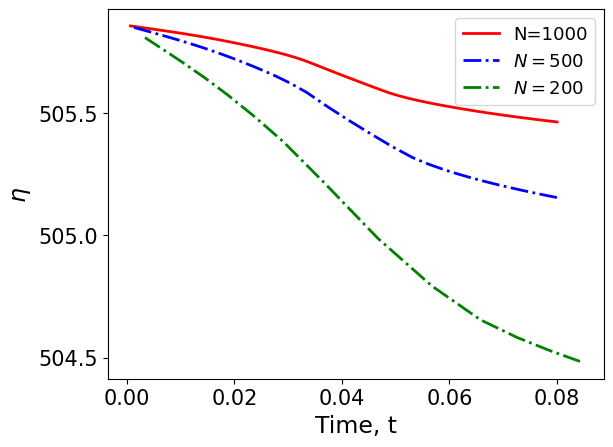}
\label{Fig:CAWtype3_ent}
\caption{Entropy}
\end{subfigure}
\hspace{-0.2cm}
\begin{subfigure}[b]{0.31\textwidth}
\centering
\includegraphics[width=\textwidth]{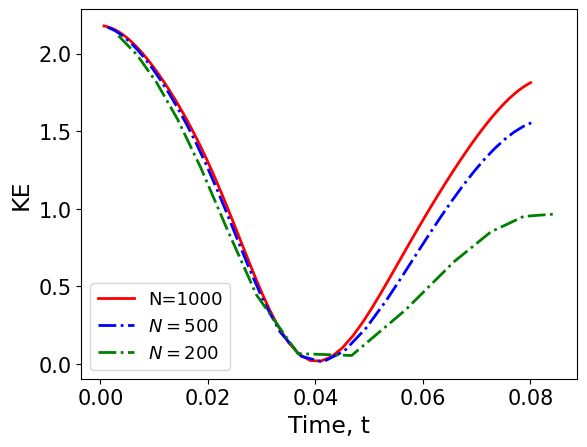}
\label{Fig:CAWtype3_KE}
\caption{Kinetic energy}
\end{subfigure}
\hspace{-0.2cm}
\begin{subfigure}[b]{0.32\textwidth}
\centering
\includegraphics[width=\textwidth]{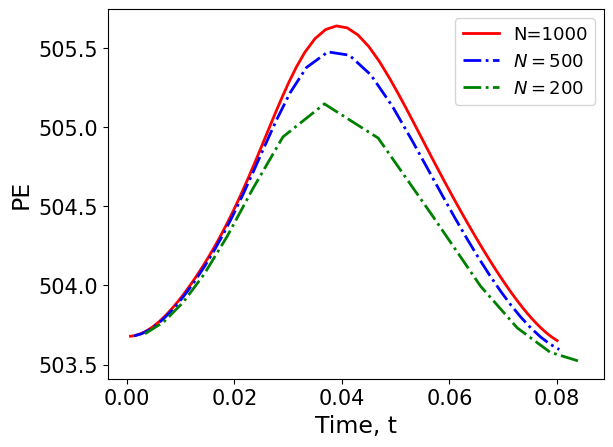}
\label{Fig:CAWtype3_PE}
\caption{Potential energy}
\end{subfigure}
\caption{\centering \textbf{Colliding acoustic waves problem} with $\epsilon=0.1$ using type 3 space discretisation.}
\label{Fig:CAWtype3}
\end{figure}

The domain of the problem is $\Omega=[-1,1]$, and the initial conditions are: 
\begin{gather}
    \rho^0(x)=0.955+0.5\epsilon \left(1-\cos \left(2\pi x\right) \right) \\
    u_1^0(x)=-\text{sign}(x)\sqrt{\gamma} \left(1- \cos \left( 2\pi x\right) \right).
\end{gather}
The parametric values are: $\kappa=1$, $\gamma=1.4$ and $\epsilon=0.1$. It is to be noted that the initial condition is not well-prepared. Periodic boundary conditions are used for this problem. Density and momentum plots are presented for different values of the final time $T=0.04,0.06,0.08$. We also present global entropy, kinetic energy, and potential energy with respect to time. $ARS(1,1,1)$ IMEX time discretisation is used. Figures \ref{Fig:CAWtype1}, \ref{Fig:CAWtype2} and \ref{Fig:CAWtype3} show the plots obtained by using type 1, type 2 and type 3 space discretisations, respectively. All three types of space discretisation result in entropy decay. While the global entropy decays, as expected, the kinetic and potential energies separately do not decay. Further, for this problem, there is no significant difference in density and momentum plots between the three types of space discretisation. \\
The results corresponding to $C=0.8$ are shown in Figures \ref{Fig:CAWtype1}, \ref{Fig:CAWtype2}, and \ref{Fig:CAWtype3}. Similar entropy decay is observed for $C=0.1,0.2,\dots,0.9$.

\subsection{Riemann problem}
\begin{figure}[h!]
\centering
\begin{subfigure}[b]{0.32\textwidth}
\centering
\includegraphics[width=\textwidth]{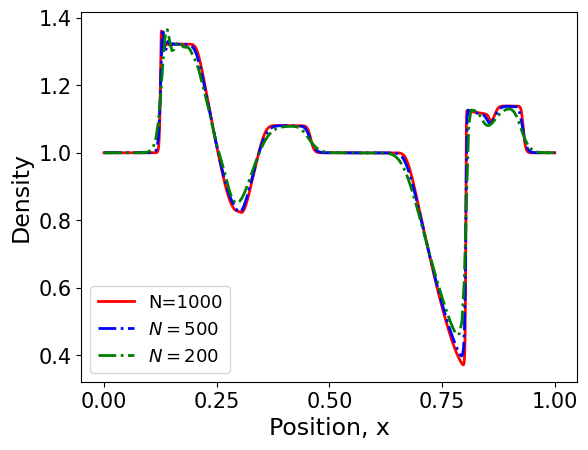}
\label{Fig:RPtype2_den_eps08}
\end{subfigure}
\hspace{-0.2cm}
\begin{subfigure}[b]{0.32\textwidth}
\centering
\includegraphics[width=\textwidth]{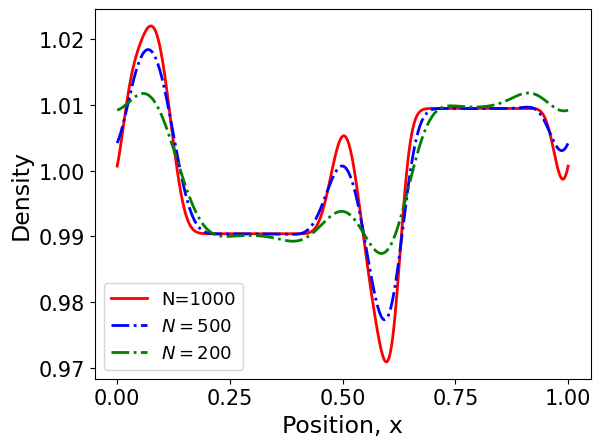}
\label{Fig:RPtype2_den_eps03}
\end{subfigure}
\hspace{-0.2cm} 
\begin{subfigure}[b]{0.34\textwidth}
\centering
\includegraphics[width=\textwidth]{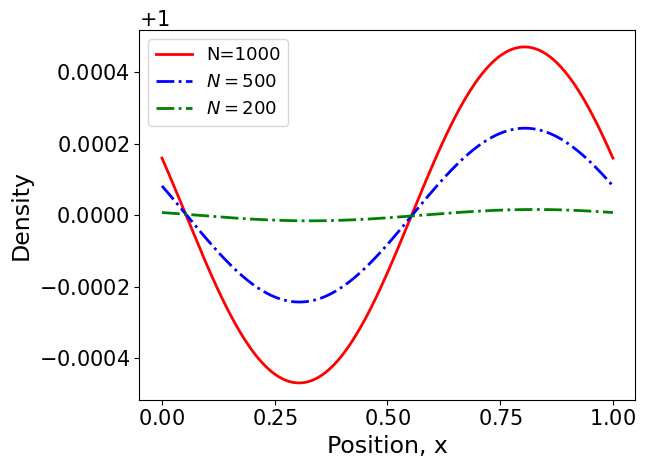}
\label{Fig:RPtype2_den_eps005}
\end{subfigure}
\vfill
\begin{subfigure}[b]{0.32\textwidth}
\centering
\includegraphics[width=\textwidth]{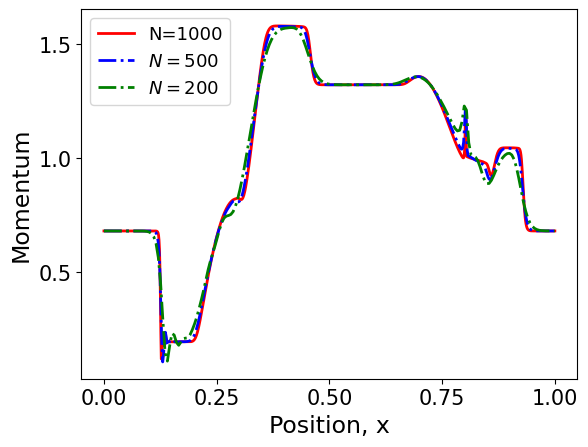}
\label{Fig:RPtype2_mom_eps08}
\end{subfigure}
\hspace{-0.2cm}
\begin{subfigure}[b]{0.32\textwidth}
\centering
\includegraphics[width=\textwidth]{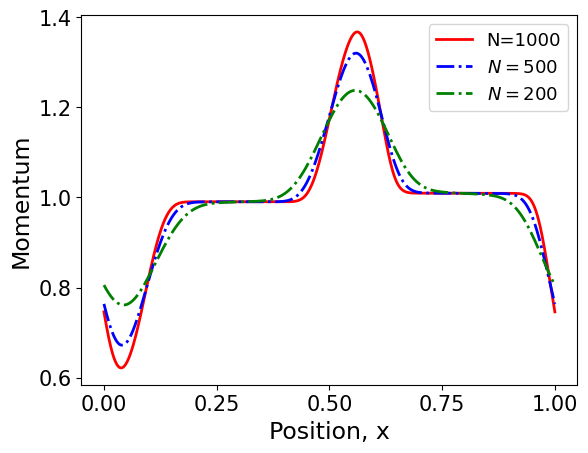}
\label{Fig:RPtype2_mom_eps03}
\end{subfigure}
\hspace{-0.2cm}
\begin{subfigure}[b]{0.34\textwidth}
\centering
\includegraphics[width=\textwidth]{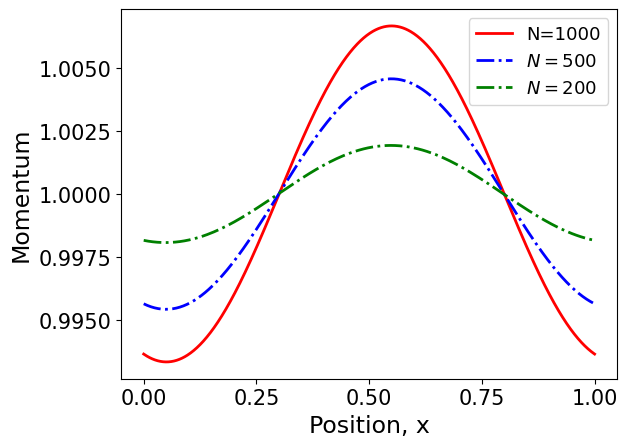}
\label{Fig:RPtype2_mom_eps005}
\end{subfigure}
\vfill
\begin{subfigure}[b]{0.32\textwidth}
\centering
\includegraphics[width=\textwidth]{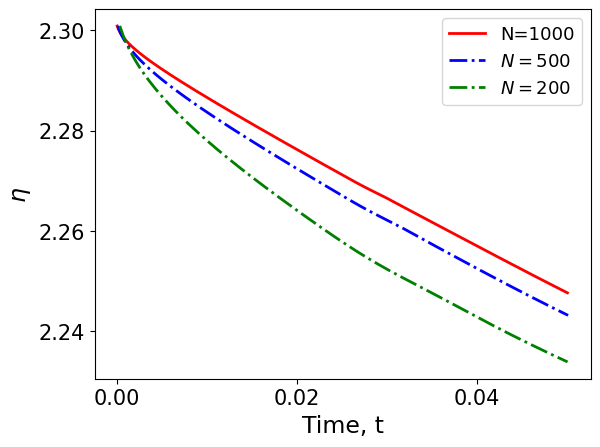}
\label{Fig:RPtype2_ent_eps08}
\caption{$\epsilon=0.8$}
\end{subfigure}
\hspace{-0.2cm}
\begin{subfigure}[b]{0.33\textwidth}
\centering
\includegraphics[width=\textwidth]{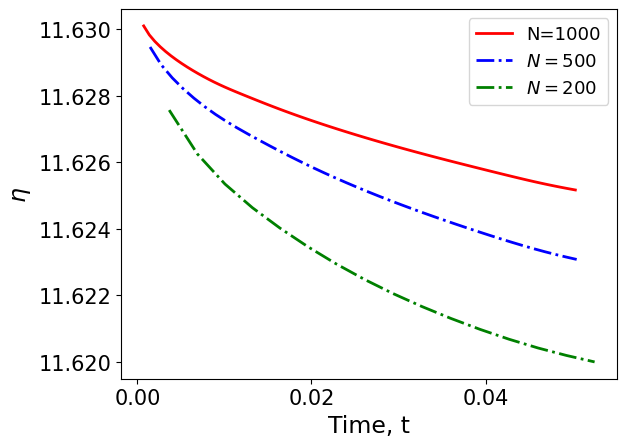}
\label{Fig:RPtype2_ent_eps03}
\caption{$\epsilon=0.3$}
\end{subfigure}
\hspace{-0.2cm}
\begin{subfigure}[b]{0.32\textwidth}
\centering
\includegraphics[width=\textwidth]{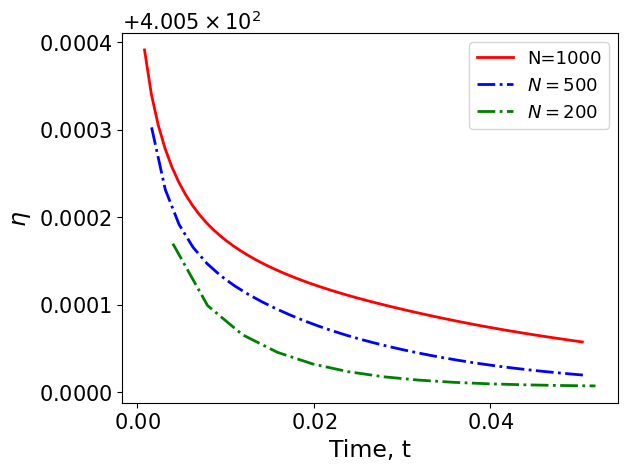}
\label{Fig:RPtype2_ent_eps005}
\caption{$\epsilon=0.05$}
\end{subfigure}
\caption{\centering \textbf{Riemann problem} at $T=0.05$ using type 2 space discretisation.}
\label{Fig:RP_type2}
\end{figure}

\begin{figure}[h!]
\centering
\begin{subfigure}[b]{0.32\textwidth}
\centering
\includegraphics[width=\textwidth]{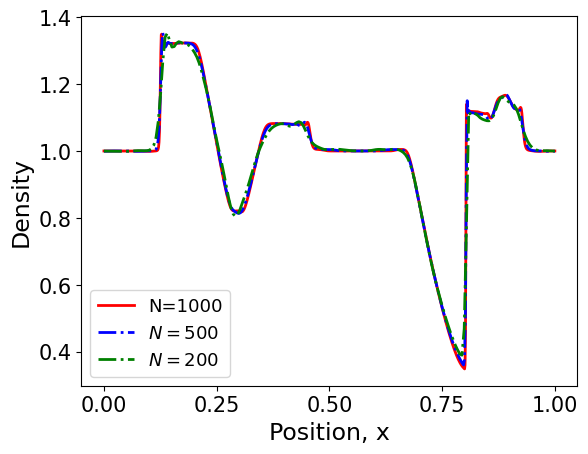}
\label{Fig:RPtype3ES1_den_eps08}
\end{subfigure}
\hspace{-0.2cm}
\begin{subfigure}[b]{0.32\textwidth}
\centering
\includegraphics[width=\textwidth]{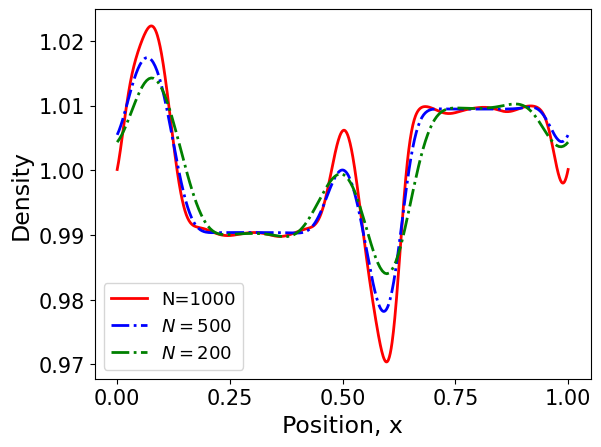}
\label{Fig:RPtype3ES1_den_eps03}
\end{subfigure}
\hspace{-0.2cm} 
\begin{subfigure}[b]{0.34\textwidth}
\centering
\includegraphics[width=\textwidth]{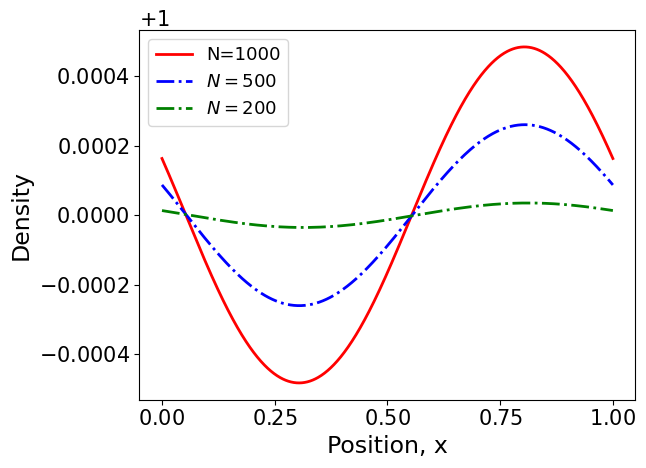}
\label{Fig:RPtype3ES1_den_eps005}
\end{subfigure}
\vfill
\begin{subfigure}[b]{0.32\textwidth}
\centering
\includegraphics[width=\textwidth]{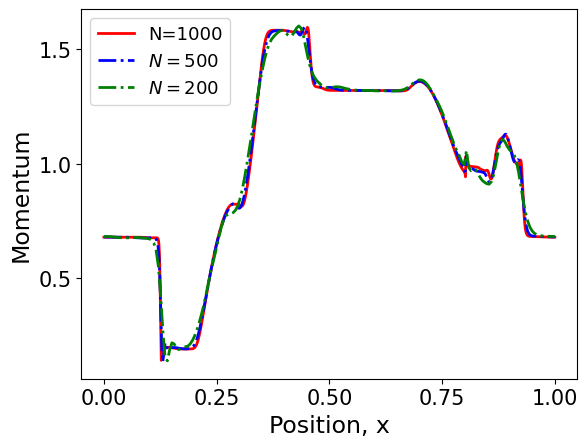}
\label{Fig:RPtype3ES1_mom_eps08}
\end{subfigure}
\hspace{-0.2cm}
\begin{subfigure}[b]{0.32\textwidth}
\centering
\includegraphics[width=\textwidth]{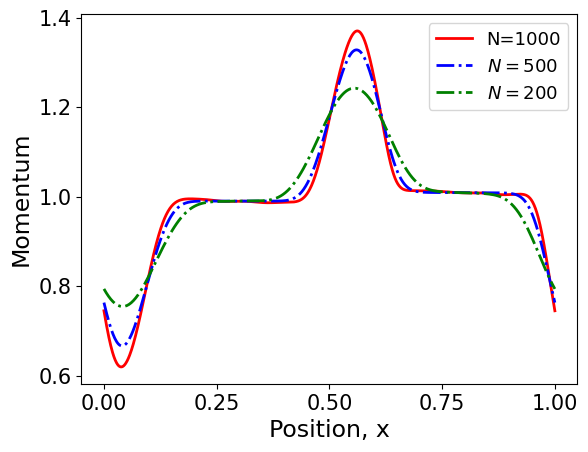}
\label{Fig:RPtype3ES1_mom_eps03}
\end{subfigure}
\hspace{-0.2cm}
\begin{subfigure}[b]{0.34\textwidth}
\centering
\includegraphics[width=\textwidth]{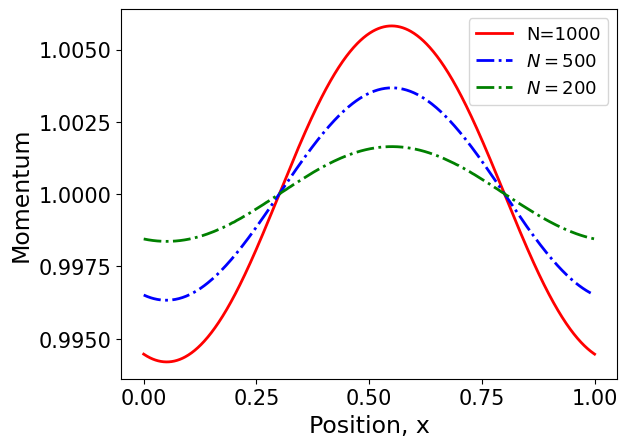}
\label{Fig:RPtype3ES1_mom_eps005}
\end{subfigure}
\vfill
\begin{subfigure}[b]{0.31\textwidth}
\centering
\includegraphics[width=\textwidth]{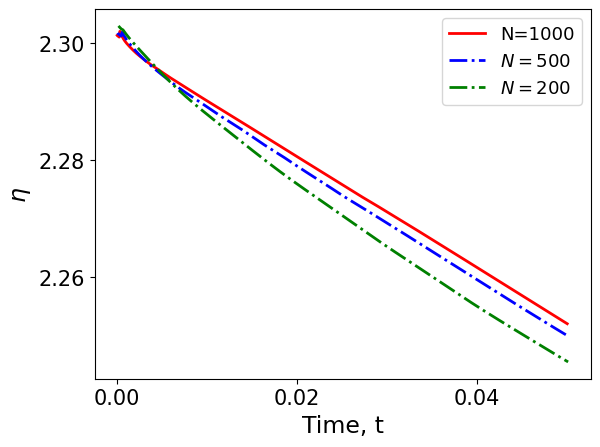}
\label{Fig:RPtype3ES1_ent_eps08}
\caption{$\epsilon=0.8$}
\end{subfigure}
\hspace{-0.2cm}
\begin{subfigure}[b]{0.32\textwidth}
\centering
\includegraphics[width=\textwidth]{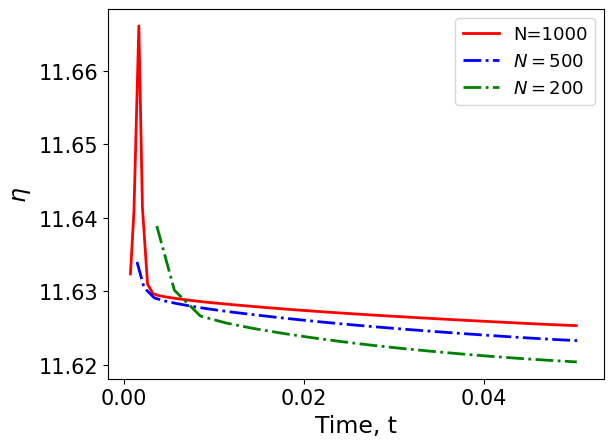}
\label{Fig:RPtype3ES1_ent_eps03}
\caption{$\epsilon=0.3$}
\end{subfigure}
\hspace{-0.2cm}
\begin{subfigure}[b]{0.34\textwidth}
\centering
\includegraphics[width=\textwidth]{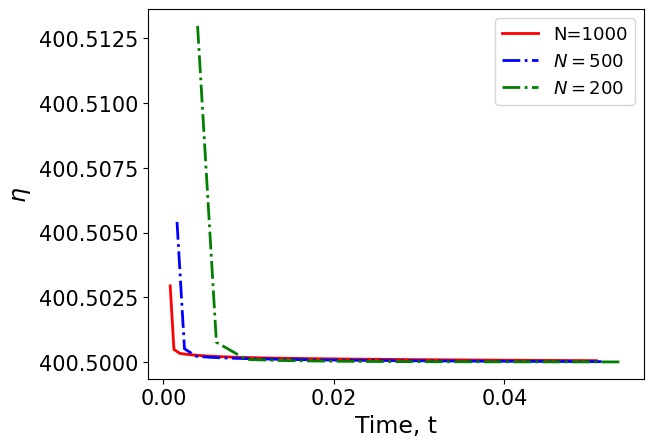}
\label{Fig:RPtype3ES1_ent_eps005}
\caption{$\epsilon=0.05$}
\end{subfigure}
\caption{\centering \textbf{Riemann problem} at $T=0.05$ using $1^{st}$ order type 3 space discretisation with $q=1$.}
\label{Fig:RP_type3ES1}
\end{figure}

\begin{figure}[h!]
\centering
\begin{subfigure}[b]{0.32\textwidth}
\centering
\includegraphics[width=\textwidth]{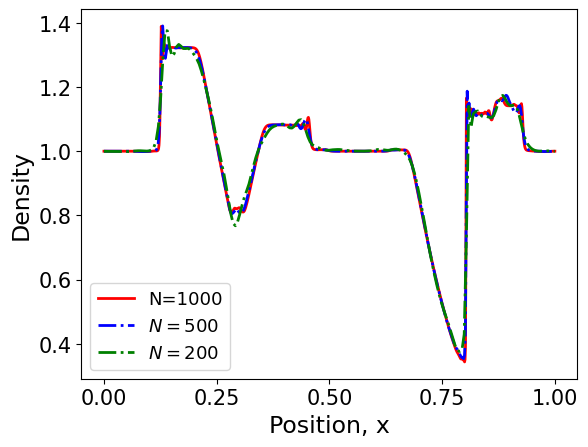}
\label{Fig:RPtype3ES2_den_eps08}
\end{subfigure}
\hspace{-0.2cm}
\begin{subfigure}[b]{0.32\textwidth}
\centering
\includegraphics[width=\textwidth]{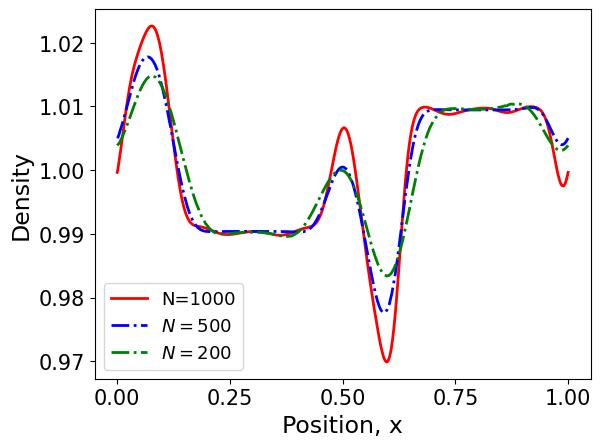}
\label{Fig:RPtype3ES2_den_eps03}
\end{subfigure}
\hspace{-0.2cm} 
\begin{subfigure}[b]{0.34\textwidth}
\centering
\includegraphics[width=\textwidth]{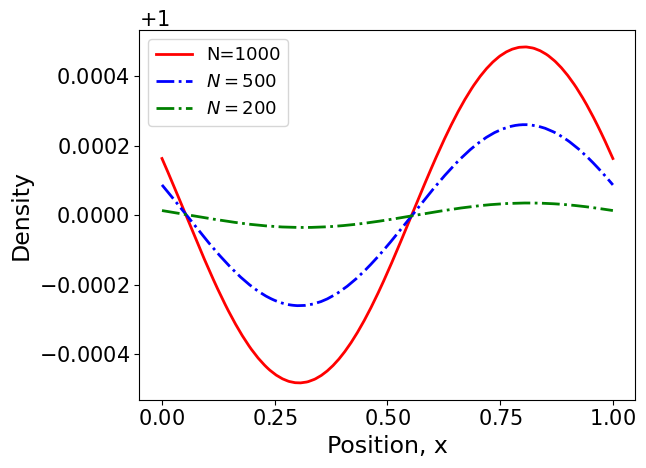}
\label{Fig:RPtype3ES2_den_eps005}
\end{subfigure}
\vfill
\begin{subfigure}[b]{0.32\textwidth}
\centering
\includegraphics[width=\textwidth]{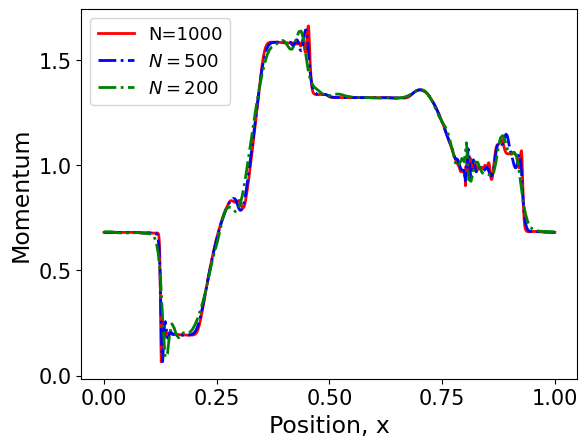}
\label{Fig:RPtype3ES2_mom_eps08}
\end{subfigure}
\hspace{-0.2cm}
\begin{subfigure}[b]{0.32\textwidth}
\centering
\includegraphics[width=\textwidth]{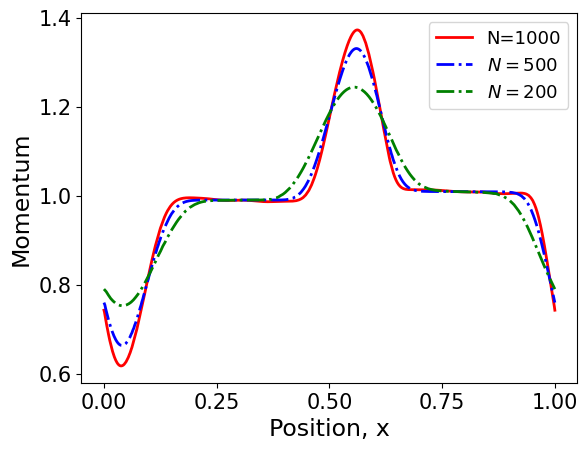}
\label{Fig:RPtype3ES2_mom_eps03}
\end{subfigure}
\hspace{-0.2cm}
\begin{subfigure}[b]{0.34\textwidth}
\centering
\includegraphics[width=\textwidth]{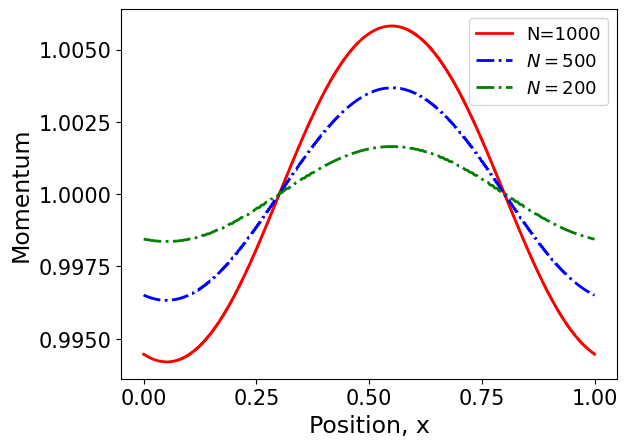}
\label{Fig:RPtype3ES2_mom_eps005}
\end{subfigure}
\vfill
\begin{subfigure}[b]{0.31\textwidth}
\centering
\includegraphics[width=\textwidth]{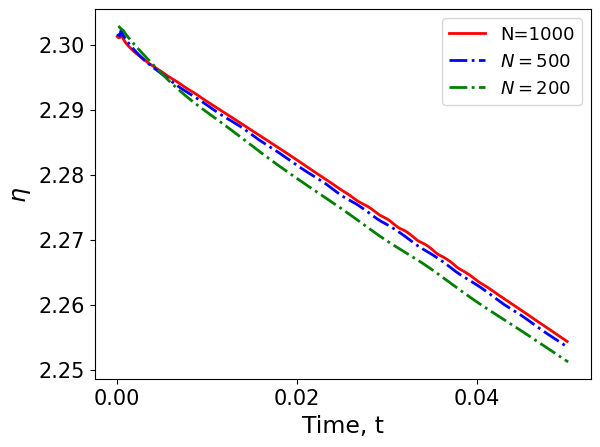}
\label{Fig:RPtype3ES2_ent_eps08}
\caption{$\epsilon=0.8$}
\end{subfigure}
\hspace{-0.2cm}
\begin{subfigure}[b]{0.32\textwidth}
\centering
\includegraphics[width=\textwidth]{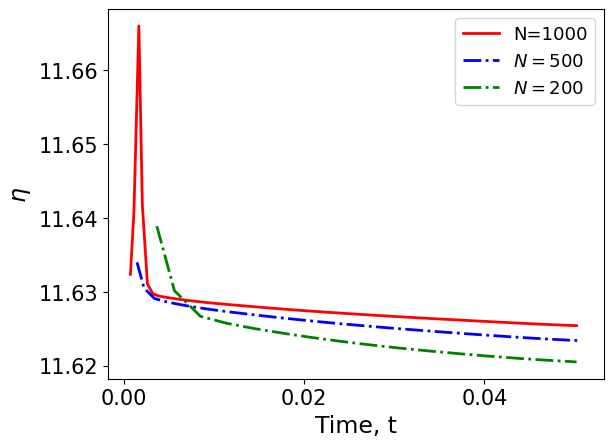}
\label{Fig:RPtype3ES2_ent_eps03}
\caption{$\epsilon=0.3$}
\end{subfigure}
\hspace{-0.2cm}
\begin{subfigure}[b]{0.34\textwidth}
\centering
\includegraphics[width=\textwidth]{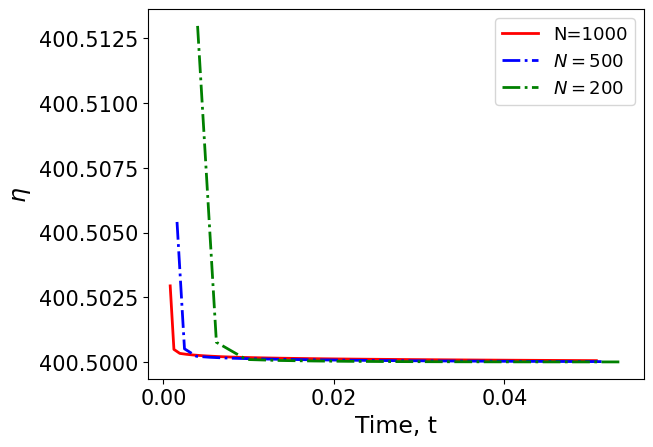}
\label{Fig:RPtype3ES2_ent_eps005}
\caption{$\epsilon=0.05$}
\end{subfigure}
\caption{\centering \textbf{Riemann problem} at $T=0.05$ using $2^{nd}$ order type 3 space discretisation with $q=1$.}
\label{Fig:RP_type3ES2}
\end{figure}

This problem is taken from \cite{Basic_AP_Degond1}. The domain is $\Omega=[0,1]$ and the initial conditions are: 
\begin{eqnarray}
    \rho^0(x)=1, & \left(\rho u_1\right)^0(x)=1-\epsilon^2/2, & \text{if } x \in [0,0.2] \cup [0.8,1] \\
    \rho^0(x)=1+\epsilon^2, & \left(\rho u_1\right)^0(x)=1, & \text{if } x \in (0.2,0.3] \\
    \rho^0(x)=1, & \left(\rho u_1\right)^0(x)=1+\epsilon^2/2, & \text{if } x \in (0.3,0.7] \\
    \rho^0(x)=1-\epsilon^2, & \left(\rho u_1\right)^0(x)=1, & \text{if } x \in (0.7,0.8]. 
\end{eqnarray}
The parametric values are: $\kappa=1$ and $\gamma=2$. Periodic boundary conditions are used for this problem. Density, momentum and global entropy (vs. time) plots are presented at $T=0.05$ for different values of $\epsilon$, such as $\epsilon=0.8,0.3,0.05$. $ARS(1,1,1)$ IMEX time discretisation is used. Figures \ref{Fig:RP_type2}, \ref{Fig:RP_type3ES1} and \ref{Fig:RP_type3ES2} show the plots obtained by using type 2, type 3 first and second order entropy stable space discretisations, respectively. Entropy decay is observed for all tested values of $\epsilon$, while using type 2 and type 3 space discretisations. However, while type 1 space discretisation results in entropy decay for $\epsilon=0.05,0.3$ (figures not shown), this discretisation results in blow-up for $\epsilon=0.8$. This is expected as type 1 space discretisation involves central discretisation of $\nabla \cdot (\rho \mathbf{u})$.\\
Results corresponding to $C=0.8$ are shown for $\epsilon=0.05,0.3$ in Figures \ref{Fig:RP_type2}, \ref{Fig:RP_type3ES1} and \ref{Fig:RP_type3ES2}. For these values of $\epsilon$, similar entropy decay is observed when using $C=0.1,0.2,\dots,0.9$. On the other hand for $\epsilon=0.8$, $C=0.2$ is used in Figure \ref{Fig:RP_type2}, and  $C=0.1$ is used in Figures \ref{Fig:RP_type3ES1} and \ref{Fig:RP_type3ES2}. 

\subsection{Gresho vortex problem}
\begin{figure}[h!]
\centering
\begin{subfigure}[b]{0.32\textwidth}
\centering
\includegraphics[width=\textwidth]{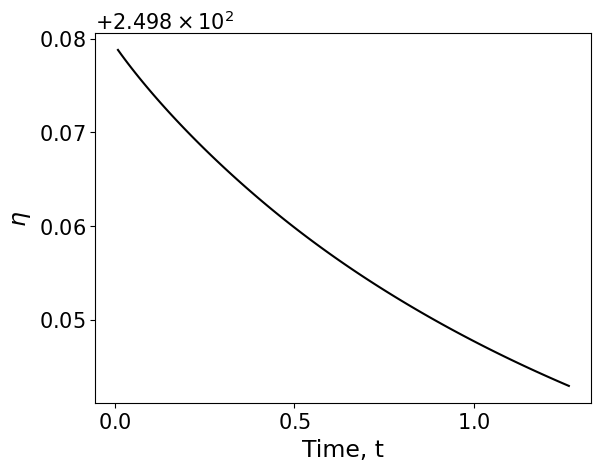}
\caption{Entropy, $\epsilon=0.1$}
\label{Fig:GV_eps01_eta}
\end{subfigure}
\hspace{-0.2cm}
\begin{subfigure}[b]{0.32\textwidth}
\centering
\includegraphics[width=\textwidth]{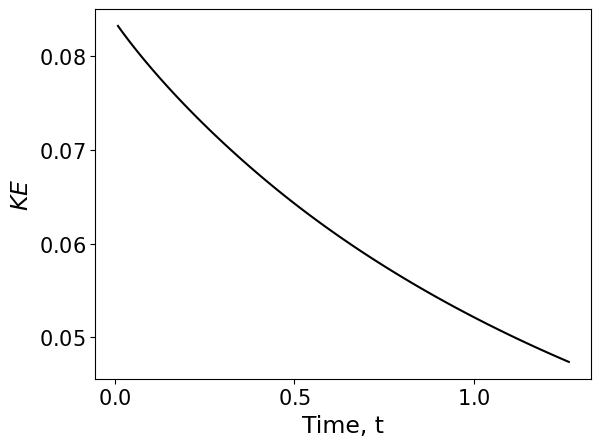}
\caption{KE, $\epsilon=0.1$}
\label{Fig:GV_eps01_KE}
\end{subfigure}
\hspace{-0.2cm}
\begin{subfigure}[b]{0.31\textwidth}
\centering
\includegraphics[width=\textwidth]{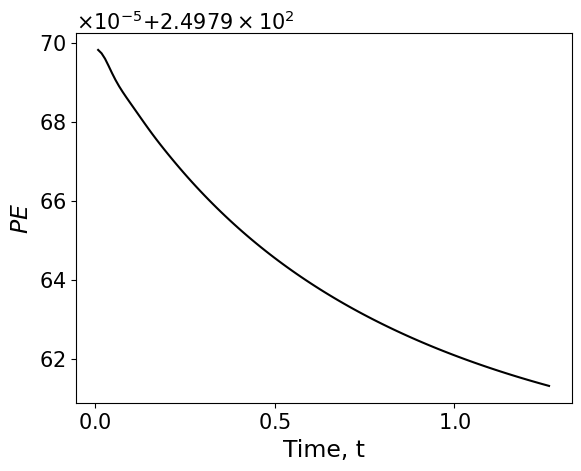}
\caption{PE, $\epsilon=0.1$}
\label{Fig:GV_eps01_PE}
\end{subfigure}
\hspace{-0.2cm}
\vfill
\begin{subfigure}[b]{0.32\textwidth}
\centering
\includegraphics[width=\textwidth]{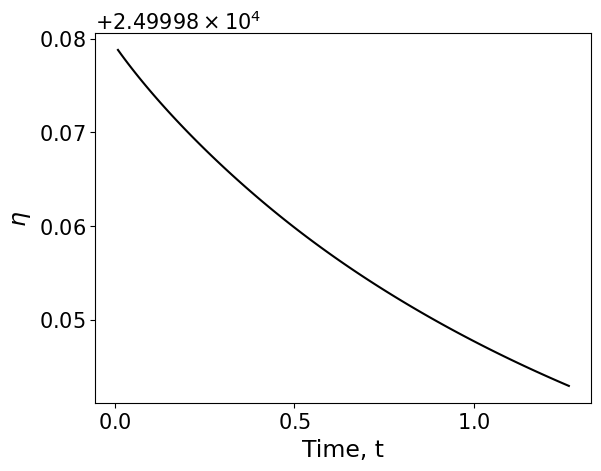}
\caption{Entropy, $\epsilon=0.01$}
\label{Fig:GV_eps001_eta}
\end{subfigure}
\hspace{-0.2cm}
\begin{subfigure}[b]{0.32\textwidth}
\centering
\includegraphics[width=\textwidth]{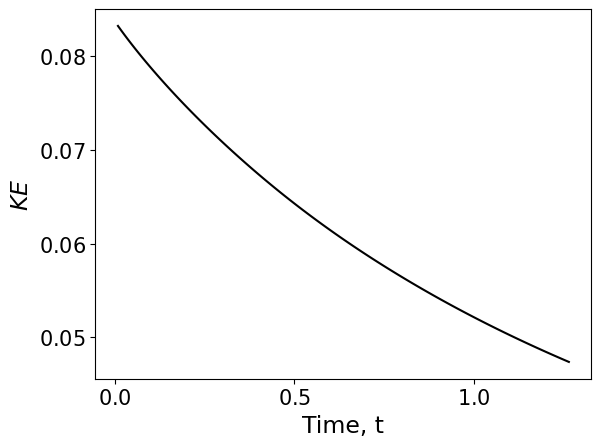}
\caption{KE, $\epsilon=0.01$}
\label{Fig:GV_eps001_KE}
\end{subfigure}
\hspace{-0.2cm}
\begin{subfigure}[b]{0.31\textwidth}
\centering
\includegraphics[width=\textwidth]{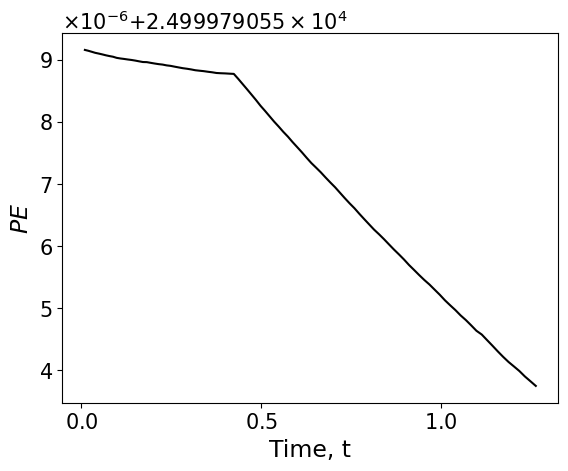}
\caption{PE, $\epsilon=0.01$}
\label{Fig:GV_eps001_PE}
\end{subfigure}
\vfill
\begin{subfigure}[b]{0.32\textwidth}
\centering
\includegraphics[width=\textwidth]{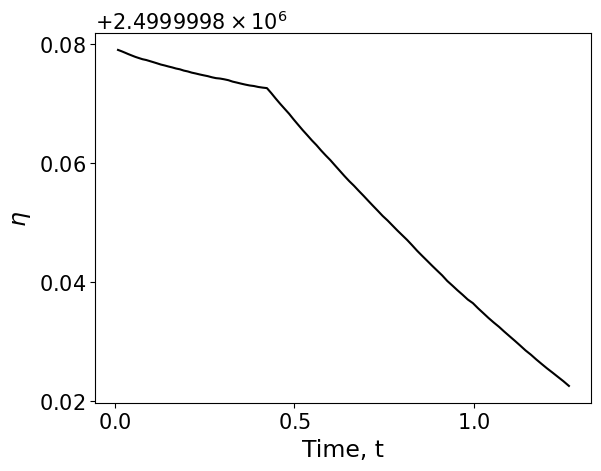}
\caption{Entropy, $\epsilon=0.001$}
\label{Fig:GV_eps0001_eta}
\end{subfigure}
\hspace{-0.2cm}
\begin{subfigure}[b]{0.32\textwidth}
\centering
\includegraphics[width=\textwidth]{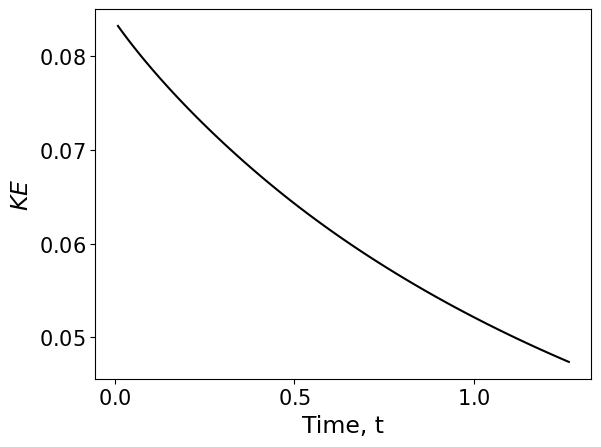}
\caption{KE, $\epsilon=0.001$}
\label{Fig:GV_eps0001_KE}
\end{subfigure}
\hspace{-0.2cm}
\begin{subfigure}[b]{0.32\textwidth}
\centering
\includegraphics[width=\textwidth]{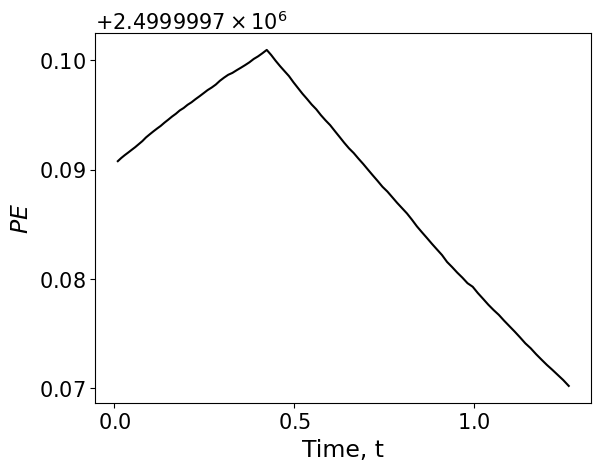}
\caption{PE, $\epsilon=0.001$}
\label{Fig:GV_eps0001_PE}
\end{subfigure}
\caption{\centering \textbf{Gresho vortex problem:} Entropy, KE and PE plots using space discretisation type 2 for $\epsilon=0.1,0.01,0.001$ on $50 \times 50$ grid.}
\label{Fig:GV}
\end{figure}

\begin{figure}[h!]
\centering
\begin{subfigure}[b]{0.32\textwidth}
\centering
\includegraphics[width=\textwidth]{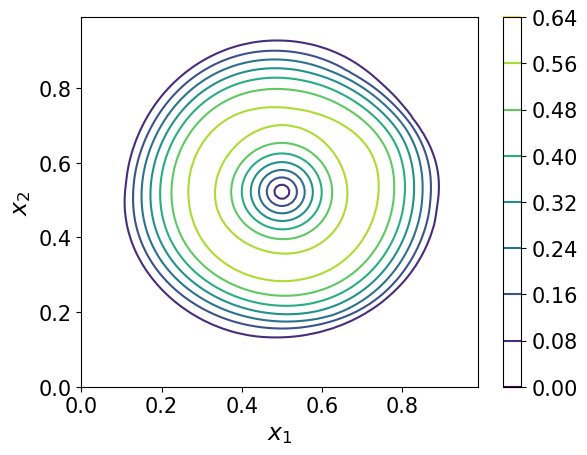}
\caption{Ma ratio, $\epsilon=0.1$}
\label{Fig:GV_eps01_Maratio}
\end{subfigure}
\hspace{-0.2cm}
\begin{subfigure}[b]{0.32\textwidth}
\centering
\includegraphics[width=\textwidth]{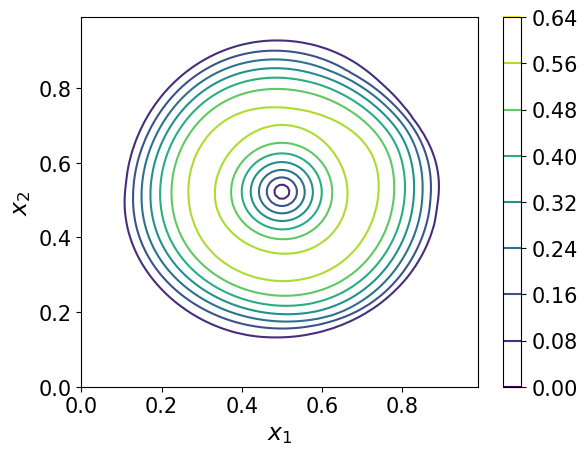}
\caption{Ma ratio, $\epsilon=0.01$}
\label{Fig:GV_eps001_Maratio}
\end{subfigure}
\hspace{-0.2cm}
\begin{subfigure}[b]{0.32\textwidth}
\centering
\includegraphics[width=\textwidth]{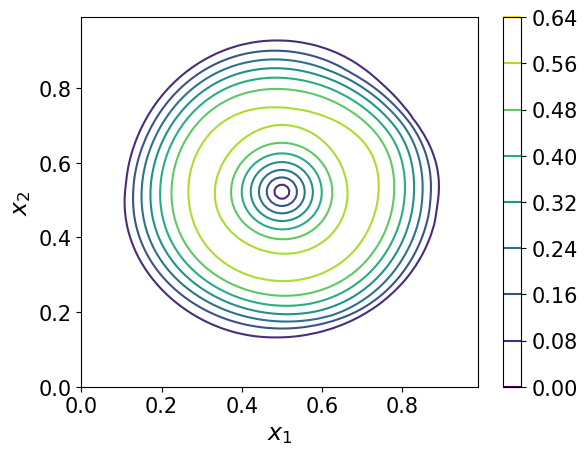}
\caption{Ma ratio, $\epsilon=0.001$}
\label{Fig:GV_eps0001_Maratio}
\end{subfigure}
\caption{\centering \textbf{Gresho vortex problem:} Ma ratio plots using space discretisation type 2 for $\epsilon=0.1,0.01,0.001$ on $100 \times 100$ grid.}
\label{Fig:GV_Maratio}
\end{figure}

This problem is taken from \cite{Gresho1990,Greshoprblm_Rieper2011}. A vortex of radius $R=0.4$ centered at $(x_{1_0},x_{2_0})=(0.5,0.5)$ is considered at initial time $t=0$. The initial background state is considered as: $\rho_0=1$, $\mathbf{u}_0=\left( u_{1_0},0 \right)^T$ with $u_{1_0}=0.1$, $p_0=1$ and hence $a_0=\sqrt{\frac{\gamma p_0}{\rho_0}}=\sqrt{\gamma}$. \\
The line velocity of the vortex is
\begin{equation*}
    u_{\theta}(r) = \left\{ \begin{matrix}
        2\frac{r}{R} & \text{if } 0\leq r < \frac{R}{2} \\
        2 \left( 1- \frac{r}{R} \right) & \text{if } \frac{R}{2} \leq r < R \\
        0 & \text{if } r\geq R
    \end{matrix} \right.,
\end{equation*}
and the velocity components in Cartesian coordinates are
\begin{equation*}
    u_1(x_1,x_2)=u_{1_0} - \frac{x_2 - x_{2_0}}{r} u_{\theta}(r), \ u_2(x_1,x_2)= \frac{x_1 - x_{1_0}}{r} u_{\theta}(r).
\end{equation*}
Here, $r=\sqrt{(x_1-0.5)^2+(x_2-0.5)^2}$. Upon balancing the pressure gradient and centrifugal force $\left( \textit{i.e., } \rho_0\frac{u_{\theta}^2}{r}= \frac{1}{\epsilon^2}\frac{\partial p}{\partial r} \right)$, pressure is derived as:
\begin{equation*}
    p(r)=p_0 + \epsilon^2 \left\{ \begin{matrix}
        2\frac{r^2}{R^2} +2-\log 16 & \text{if } 0\leq r < \frac{R}{2} \\
        2\frac{r^2}{R^2} - 8 \frac{r}{R} +4 \log \left( \frac{r}{R} \right)+ 6 & \text{if } \frac{R}{2} \leq r < R \\
        0 & \text{if } r\geq R
    \end{matrix} \right..
\end{equation*}
We assume adiabatic compression $p=\rho^{\gamma}$ with $\gamma=1.4$, and use asymptotic expansion: $p=p_0 + \epsilon^2 p_2$, $\rho=\rho_0 + \epsilon^2 \rho_2$. Comparing $\frac{\partial p}{\partial \epsilon}=2\epsilon p_2$ and $\frac{\partial p}{\partial \epsilon}= \frac{\partial p}{\partial \rho} \frac{\partial \rho}{\epsilon} = \frac{\gamma p }{\rho} 2\epsilon \rho_2$, we obtain $\rho_2 = \frac{p_2}{\gamma}$ by noting that $p=1,\rho=1$ up to the leading order. $\rho_2$ can then be used to evaluate $\rho$ as $\rho=\rho_0 + \epsilon^2 \rho_2$. \\
Periodic boundary conditions are imposed in both directions. The problem is simulated using $ARS(1,1,1)$ IMEX time discretisation and type 2 space discretisation for $\epsilon=0.1,0.01,0.001$. The following quantities are observed:
\begin{gather}
     \eta = 1/2 \ \rho (u_1^2 + u_2^2) + (1/\epsilon^2) p/(\gamma-1), \ \text{PE } = (1/\epsilon^2) p/(\gamma-1), \\ \text{KE } = \frac{1}{2}\left((u_1 - u_{1_0})^2 + u_2^2 \right) \text{ where } u_{1_0} \text{ is the the background velocity} \\
    \text{Mach number ratio, } \text{M}_{\text{ratio}} = \left(\sqrt{\frac{(u_1 - u_{1_0})^2 + u_2^2}{\gamma p/ \rho}} \right). 
\end{gather} 
Figure \ref{Fig:GV} shows the evolution of entropy ($\eta$), KE and PE over time upto $T=R\pi$. A CFL value of $C=0.5$ is used for $\epsilon=0.1,0.01,0.001$. As expected, the entropy and kinetic energy are decaying in time, although the potential energy separately does not decay. Figure \ref{Fig:GV_Maratio} shows the contours of $\text{M}_{\text{ratio}}$ at $T=R\pi$. \\
\Cref{tab:GV-EOC_rho,tab:GV-EOC_u1,tab:GV-EOC_u2} show the convergence rates of $\rho,u_1$ and $u_2$ for $\epsilon=0.1,0.01$ and $0.001$. It is observed that the required convergence rate is obtained for all the values of $\epsilon$.

\begin{table}[h!]
    \centering
    \renewcommand{\arraystretch}{1.3} % Increase row height for readability
    \setlength{\tabcolsep}{10pt}      % Adjust column spacing
    
    \begin{tabular}{|c|c|c|c|c|c|c|}
        \hline
        \multirow{2}*{\textbf{N}} & \multicolumn{2}{c|}{$\mathbf{\epsilon=10^{-1}}$} & 
                                     \multicolumn{2}{c|}{$\mathbf{\epsilon=10^{-2}}$} & 
                                     \multicolumn{2}{c|}{$\mathbf{\epsilon=10^{-3}}$} \\ 
        \cline{2-7}
        & $\mathbf{||\rho \textbf{ error}||_{L_2}}$ & \textbf{EOC} 
        & $\mathbf{||\rho \textbf{ error}||_{L_2}}$ & \textbf{EOC} 
        & $\mathbf{||\rho \textbf{ error}||_{L_2}}$ & \textbf{EOC} \\ 
        \hline
        10  & 0.00084239  & -      & 8.42982872 $\times 10^{-6}$  & -      & 8.45384875 $\times 10^{-8}$ & -      \\  
        20  & 0.00058693    & 0.4836 & 5.87421133$\times 10^{-6}$  & 0.4834 & 5.90367676 $\times 10^{-8}$ & 0.4805 \\  
        25  & 0.00050582    & 0.6367 & 5.06282663$\times 10^{-6}$  & 0.6363 & 5.10971760 $\times 10^{-8}$ & 0.6182 \\  
        50  & 0.00023863    & 1.0525 & 2.38934881$\times 10^{-6}$  & 1.0520 & 2.56732036 $\times 10^{-8}$ & 0.9643 \\  
        \hline
    \end{tabular}
    
    \caption{\centering \textbf{Gresho vortex problem:} Convergence rates of $L_2$ error in $\rho$ using $ARS(1,1,1)$ coupled with type 2 discretisation.} 
    \label{tab:GV-EOC_rho}
\end{table}

\begin{table}[h!]
    \centering
    \renewcommand{\arraystretch}{1.3} % Increase row height for readability
    \setlength{\tabcolsep}{10pt}      % Adjust column spacing
    
    \begin{tabular}{|c|c|c|c|c|c|c|}
        \hline
        \multirow{2}*{\textbf{N}} & \multicolumn{2}{c|}{$\mathbf{\epsilon=10^{-1}}$} & 
                                     \multicolumn{2}{c|}{$\mathbf{\epsilon=10^{-2}}$} & 
                                     \multicolumn{2}{c|}{$\mathbf{\epsilon=10^{-3}}$} \\ 
        \cline{2-7}
        & $\mathbf{||u_1 \textbf{ error}||_{L_2}}$ & \textbf{EOC} 
        & $\mathbf{||u_1 \textbf{ error}||_{L_2}}$ & \textbf{EOC} 
        & $\mathbf{||u_1 \textbf{ error}||_{L_2}}$ & \textbf{EOC} \\ 
        \hline
        10  & 0.17944234  & -      & 0.179428512  & -      & 0.179428376 & -      \\  
        20  & 0.105626    & 0.7092 & 0.105625341  & 0.7091 & 0.105625335 & 0.7091 \\  
        25  & 0.08712198    & 0.8244 & 0.0871283562  & 0.8241 & 0.0871284204 & 0.8241 \\  
        50  & 0.0367448    & 1.2095 & 0.0367560804  & 1.2092 & 0.0367561927 & 1.2092 \\  
        \hline
    \end{tabular}
    
    \caption{\centering \textbf{Gresho vortex problem:} Convergence rates of $L_2$ error in $u_1$ using $ARS(1,1,1)$ coupled with type 2 discretisation.} 
    \label{tab:GV-EOC_u1}
\end{table}

\begin{table}[h!]
    \centering
    \renewcommand{\arraystretch}{1.3} % Increase row height for readability
    \setlength{\tabcolsep}{10pt}      % Adjust column spacing
    
    \begin{tabular}{|c|c|c|c|c|c|c|}
        \hline
        \multirow{2}*{\textbf{N}} & \multicolumn{2}{c|}{$\mathbf{\epsilon=10^{-1}}$} & 
                                     \multicolumn{2}{c|}{$\mathbf{\epsilon=10^{-2}}$} & 
                                     \multicolumn{2}{c|}{$\mathbf{\epsilon=10^{-3}}$} \\ 
        \cline{2-7}
        & $\mathbf{||u_2 \textbf{ error}||_{L_2}}$ & \textbf{EOC} 
        & $\mathbf{||u_2 \textbf{ error}||_{L_2}}$ & \textbf{EOC} 
        & $\mathbf{||u_2 \textbf{ error}||_{L_2}}$ & \textbf{EOC} \\ 
        \hline
        10  & 0.1882019  & -      & 0.188182499  & -      & 0.188182305 & -      \\  
        20  & 0.11598214    & 0.64784505 & 0.115974954  & 0.6478 & 0.115974882 & 0.6478 \\  
        25  & 0.0951939    & 0.84549546 & 0.0951936708  & 0.8452 & 0.0951936685 & 0.8452 \\  
        50  & 0.04058436    & 1.1944146 & 0.0405921055  & 1.1941 & 0.0405921819e & 1.194  \\  
        \hline
    \end{tabular}
    
    \caption{\centering \textbf{Gresho vortex problem:} Convergence rates of $L_2$ error in $u_2$ using $ARS(1,1,1)$ coupled with type 2 discretisation.} 
    \label{tab:GV-EOC_u2}
\end{table}

\subsection{Travelling vortex problem}
\begin{figure}[h!]
\centering
\begin{subfigure}[b]{0.32\textwidth}
\centering
\includegraphics[width=\textwidth]{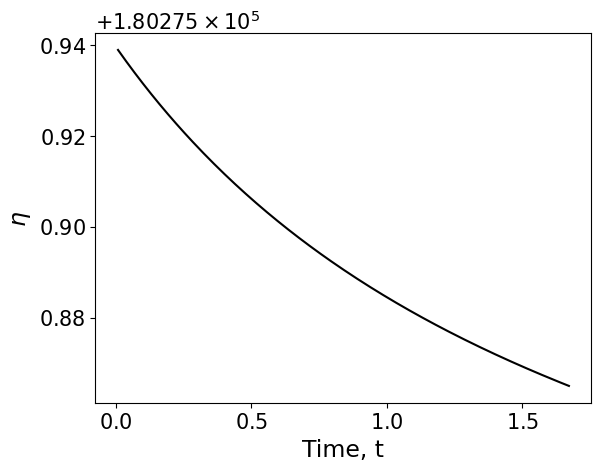}
\caption{Entropy, $\epsilon=0.1$}
\label{Fig:TV_eps01_eta}
\end{subfigure}
\hspace{-0.2cm}
\begin{subfigure}[b]{0.32\textwidth}
\centering
\includegraphics[width=\textwidth]{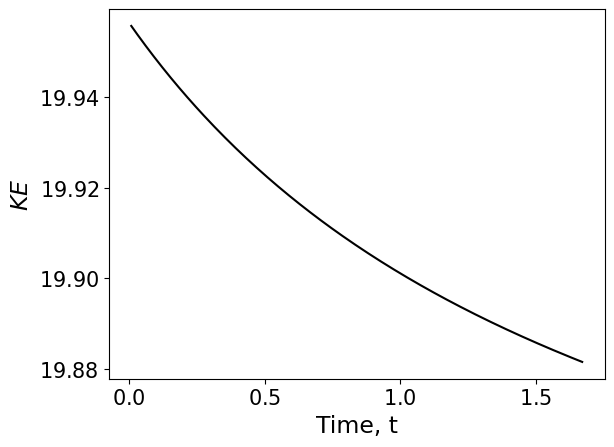}
\caption{KE, $\epsilon=0.1$}
\label{Fig:TV_eps01_KE}
\end{subfigure}
\hspace{-0.2cm}
\begin{subfigure}[b]{0.31\textwidth}
\centering
\includegraphics[width=\textwidth]{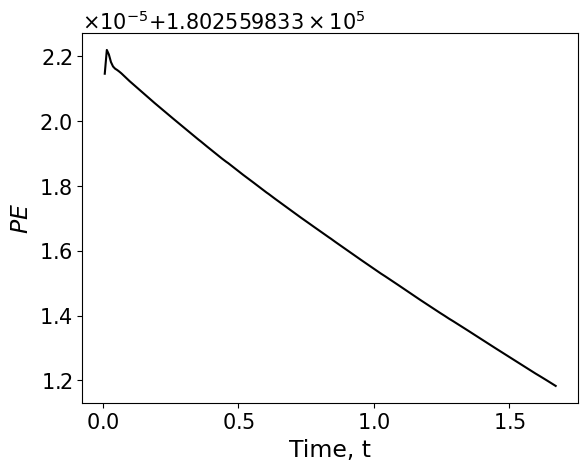}
\caption{PE, $\epsilon=0.1$}
\label{Fig:TV_eps01_PE}
\end{subfigure}
\hspace{-0.2cm}
\vfill
\begin{subfigure}[b]{0.32\textwidth}
\centering
\includegraphics[width=\textwidth]{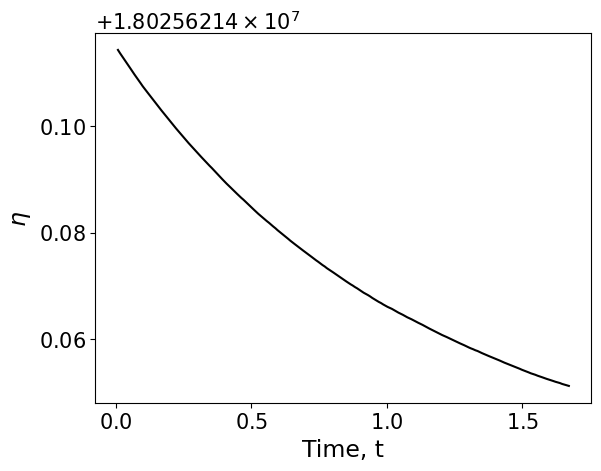}
\caption{Entropy, $\epsilon=0.01$}
\label{Fig:TV_eps001_eta}
\end{subfigure}
\hspace{-0.2cm}
\begin{subfigure}[b]{0.32\textwidth}
\centering
\includegraphics[width=\textwidth]{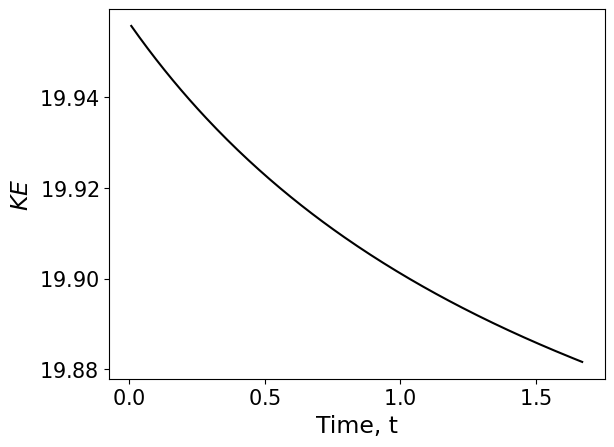}
\caption{KE, $\epsilon=0.01$}
\label{Fig:TV_eps001_KE}
\end{subfigure}
\hspace{-0.2cm}
\begin{subfigure}[b]{0.32\textwidth}
\centering
\includegraphics[width=\textwidth]{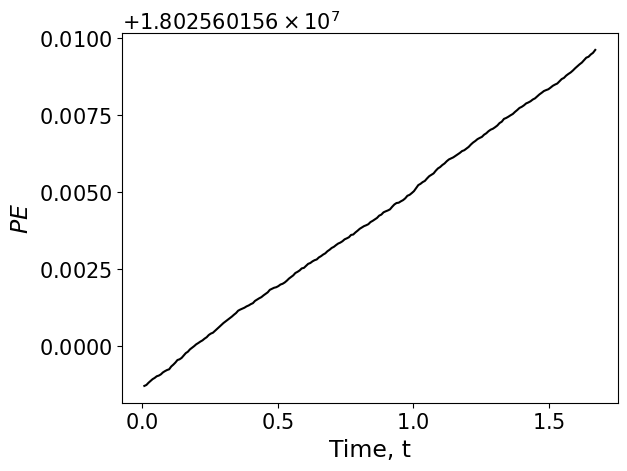}
\caption{PE, $\epsilon=0.01$}
\label{Fig:TV_eps001_PE}
\end{subfigure}
\hspace{-0.2cm}
\vfill
\begin{subfigure}[b]{0.32\textwidth}
\centering
\includegraphics[width=\textwidth]{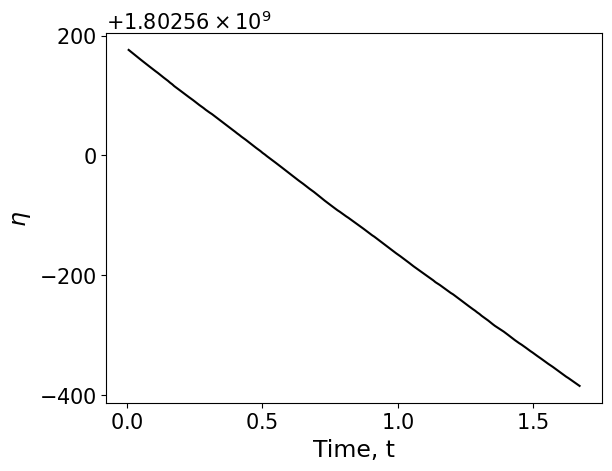}
\caption{Entropy, $\epsilon=0.001$}
\label{Fig:TV_eps0001_eta}
\end{subfigure}
\hspace{-0.2cm}
\begin{subfigure}[b]{0.32\textwidth}
\centering
\includegraphics[width=\textwidth]{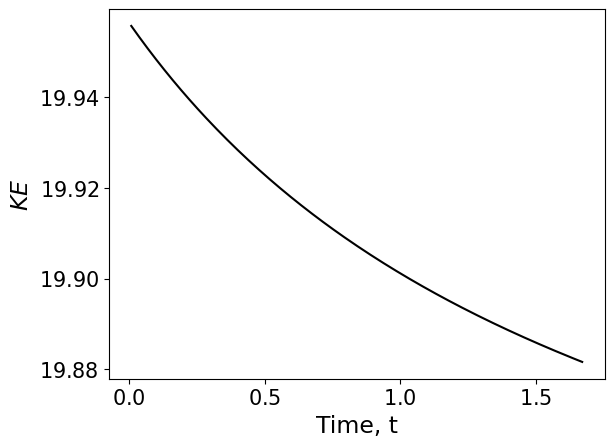}
\caption{KE, $\epsilon=0.001$}
\label{Fig:TV_eps0001_KE}
\end{subfigure}
\hspace{-0.2cm}
\begin{subfigure}[b]{0.32\textwidth}
\centering
\includegraphics[width=\textwidth]{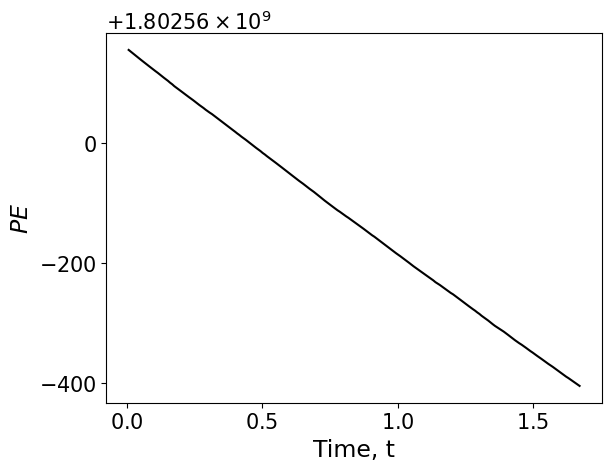}
\caption{PE, $\epsilon=0.001$}
\label{Fig:TV_eps0001_PE}
\end{subfigure}
\vfill
\begin{subfigure}[b]{0.31\textwidth}
\centering
\includegraphics[width=\textwidth]{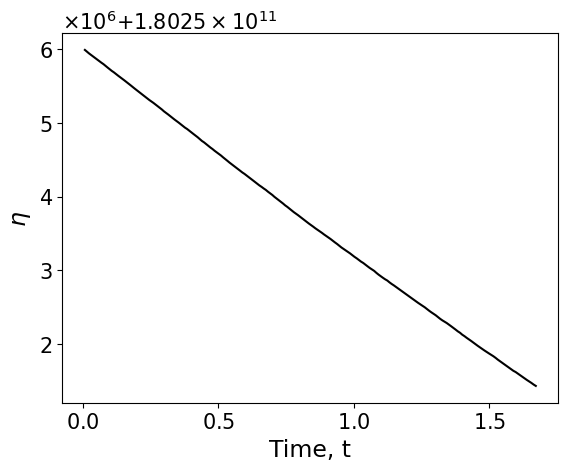}
\caption{Entropy, $\epsilon=0.0001$}
\label{Fig:TV_eps00001_eta}
\end{subfigure}
\hspace{-0.2cm}
\begin{subfigure}[b]{0.33\textwidth}
\centering
\includegraphics[width=\textwidth]{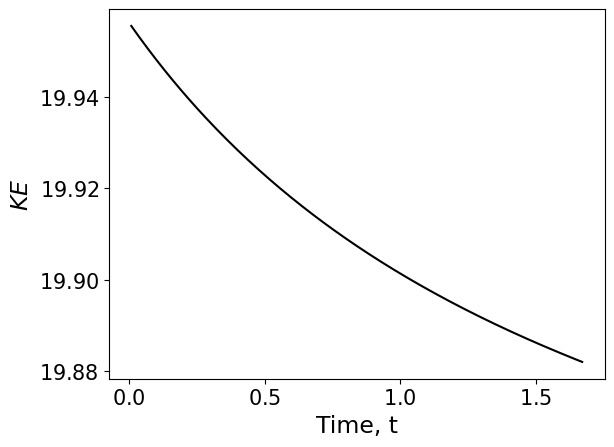}
\caption{KE, $\epsilon=0.0001$}
\label{Fig:TV_eps00001_KE}
\end{subfigure}
\hspace{-0.2cm}
\begin{subfigure}[b]{0.31\textwidth}
\centering
\includegraphics[width=\textwidth]{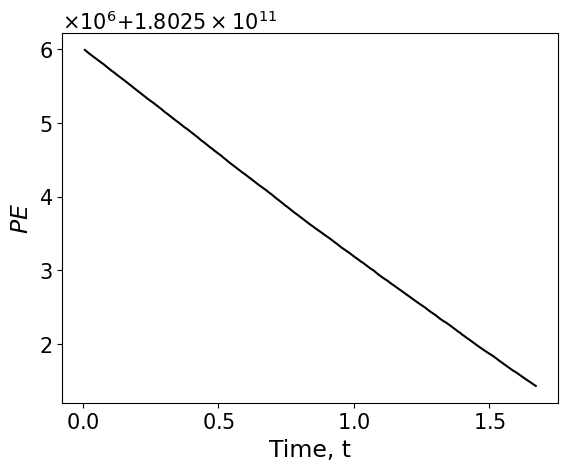}
\caption{PE, $\epsilon=0.0001$}
\label{Fig:TV_eps00001_PE}
\end{subfigure}
\caption{\centering \textbf{Travelling vortex problem:} Entropy, KE and PE plots using space discretisation type 2 for $\epsilon=0.1,0.01,0.001$ on $100 \times 100$ grid.}
\label{Fig:TV}
\end{figure}

\begin{figure}[h!]
\centering
\begin{subfigure}[b]{0.32\textwidth}
\centering
\includegraphics[width=\textwidth]{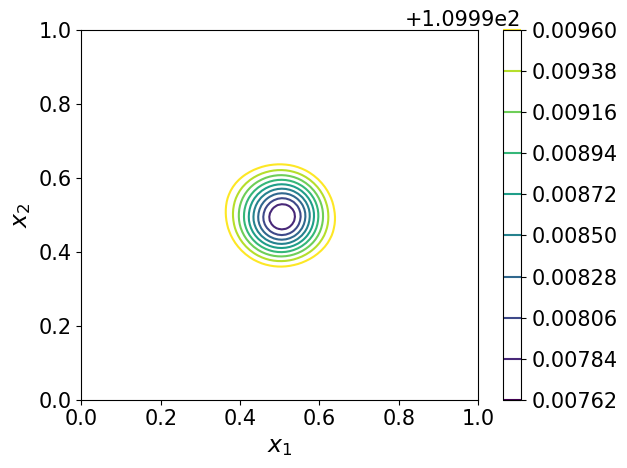}
\caption{Density, $\epsilon=10^{-1}$}
\label{Fig:TV_eps1_Density}
\end{subfigure}
\hspace{-0.2cm}
\begin{subfigure}[b]{0.32\textwidth}
\centering
\includegraphics[width=\textwidth]{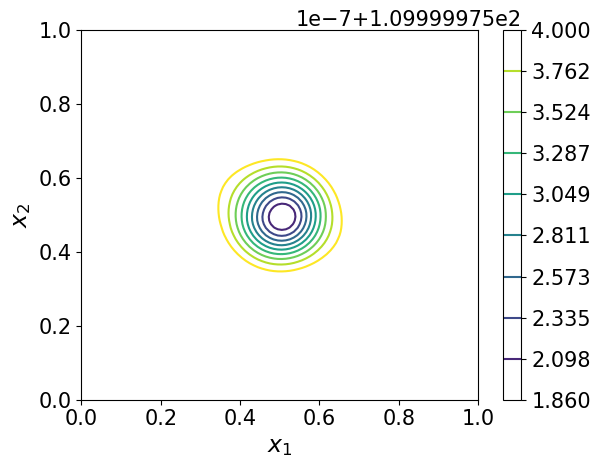}
\caption{Density, $\epsilon=10^{-3}$}
\label{Fig:TV_eps3_Density}
\end{subfigure}
\hspace{-0.2cm}
\begin{subfigure}[b]{0.32\textwidth}
\centering
\includegraphics[width=\textwidth]{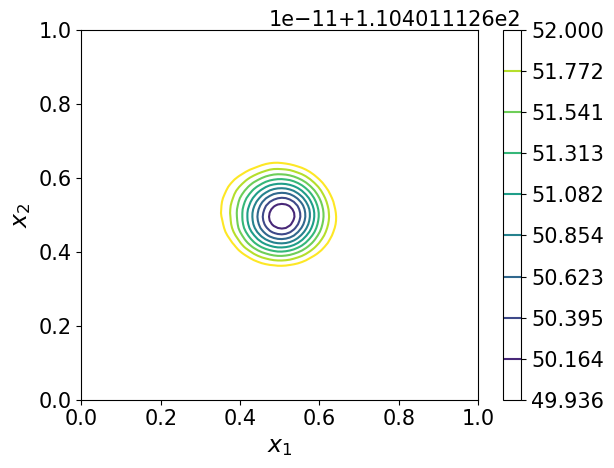}
\caption{Density, $\epsilon=10^{-5}$}
\label{Fig:TV_eps5_Density}
\end{subfigure}
\caption{\centering \textbf{Travelling vortex problem:} Density contours using space discretisation type 2 for $\epsilon=10^{-1},10^{-3},10^{-5}$ on $100 \times 100$ grid.}
\label{Fig:TV_Density}
\end{figure}

The travelling vortex problem taken from \cite{Basic_AP_Lukacova3} is considered here for the isentropic Euler system with $p=\rho^{1.4}$. The initial condition for the problem is given by,
\begin{eqnarray}
\rho^0 \left( x_1,x_2\right)= 110 + \epsilon^2\left( \frac{1.5}{4\pi}\right)^2 D\left(x_1,x_2\right) \left( k\left(r\right)-k\left( \pi\right)\right) \\
u_1^0 \left( x_1,x_2\right)=  0.6+ 1.5 \left( 1+\cos \left( r \left(x_1,x_2\right)\right)\right) D\left(x_1,x_2\right) \left( 0.5-x_2\right)  \\
u_2^0 \left( x_1,x_2\right)=  0+ 1.5 \left( 1+\cos \left( r \left(x_1,x_2\right)\right)\right) D\left(x_1,x_2\right) \left( x_1-0.5\right) 
\end{eqnarray}
with 
\begin{eqnarray}
k\left(q\right)=2\cos\left(q\right)+2q \ \sin\left(q\right) + \frac{1}{8}\cos\left(2q\right) + \frac{1}{4}q \ \sin\left(2q\right) + \frac{3}{4}q^2 \\
r \left( x_1,x_2\right) = 4 \pi \left(  \left( x_1-0.5\right)^2 + \left( x_2-0.5\right)^2 \right)^{\frac{1}{2}} \\
D\left(x_1,x_2\right) = \left\{ \begin{matrix}  1 & \text{if } r\left( x_1,x_2\right)<\pi \\ 0 & \text{otherwise}\end{matrix} \right.
\end{eqnarray}
The boundary is periodic, and the domain of the problem is $[0,1]\times[0,1]$. Figure \ref{Fig:TV} shows the evolution of entropy ($\eta$), KE and PE over time upto $T=1/0.6$ (which is the time period for the vortex to return to its initial position) for different values of $\epsilon$. A CFL value of $C=0.6$ is used for for all the values of $\epsilon$ that are considered. As expected, the entropy and kinetic energy are decaying in time, although the potential energy separately does not decay. Figure \ref{Fig:TV_Density} shows the contours of density at $T=1/0.6$. \\
\Cref{tab:TV-EOC_u1a,tab:TV-EOC_u1b,tab:TV-EOC_u2a,tab:TV-EOC_u2b} show the convergence rates of $u_1$ and $u_2$ for $\epsilon=10^{-1},10^{-2},10^{-3},10^{-4},10^{-5},10^{-6}$. It is observed that the required convergence rate is obtained for all the values of $\epsilon$.

\begin{table}[h!]
    \centering
    \renewcommand{\arraystretch}{1.3} 
    \setlength{\tabcolsep}{10pt}     
    
    \begin{tabular}{|c|c|c|c|c|c|c|}
        \hline
        \multirow{2}*{\textbf{N}} & \multicolumn{2}{c|}{$\mathbf{\epsilon=10^{-1}}$} & 
                                     \multicolumn{2}{c|}{$\mathbf{\epsilon=10^{-2}}$} & 
                                     \multicolumn{2}{c|}{$\mathbf{\epsilon=10^{-3}}$} \\ 
        \cline{2-7}
        & $\mathbf{||u_1 \textbf{ error}||_{L_2}}$ & \textbf{EOC} 
        & $\mathbf{||u_1 \textbf{ error}||_{L_2}}$ & \textbf{EOC} 
        & $\mathbf{||u_1 \textbf{ error}||_{L_2}}$ & \textbf{EOC} \\ 
        \hline
        10  & 2.53 $\times 10^{-2}$  & -      & 2.53 $\times 10^{-2}$  & -      & 2.53 $\times 10^{-2}$ & -  \\  
        20  & 1.48 $\times 10^{-2}$    & 0.71 & 1.48$\times 10^{-2}$  & 0.71 & 1.48 $\times 10^{-2}$ & 0.71  \\  
        25  & 1.24 $\times 10^{-2}$     & 0.75 & 1.24 $\times 10^{-2}$  & 0.75 & 1.24 $\times 10^{-2}$ & 0.75  \\  
        50  & 5.71 $\times 10^{-3}$    & 1.09 & 5.71 $\times 10^{-3}$  & 1.09 & 5.71 $\times 10^{-3}$ & 1.09 \\  
        \hline
    \end{tabular}
    
    \caption{\centering \textbf{Travelling vortex problem:} Convergence rates of $L_2$ error in $u_1$ using $ARS(1,1,1)$ coupled with type 2 discretisation.} 
    \label{tab:TV-EOC_u1a}
\end{table}

\begin{table}[h!]
    \centering
    \renewcommand{\arraystretch}{1.3} 
    \setlength{\tabcolsep}{10pt}     
    
    \begin{tabular}{|c|c|c|c|c|c|c|}
        \hline
        \multirow{2}*{\textbf{N}} & \multicolumn{2}{c|}{$\mathbf{\epsilon=10^{-4}}$} & 
                                     \multicolumn{2}{c|}{$\mathbf{\epsilon=10^{-5}}$} & 
                                     \multicolumn{2}{c|}{$\mathbf{\epsilon=10^{-6}}$} \\ 
        \cline{2-7}
        & $\mathbf{||u_1 \textbf{ error}||_{L_2}}$ & \textbf{EOC} 
        & $\mathbf{||u_1 \textbf{ error}||_{L_2}}$ & \textbf{EOC} 
        & $\mathbf{||u_1 \textbf{ error}||_{L_2}}$ & \textbf{EOC} \\ 
        \hline
        10  & 2.53 $\times 10^{-2}$  & -            & 2.53 $\times 10^{-2}$   & -      & 8.74 $\times 10^{-2}$ & -      \\  
        20  & 1.48 $\times 10^{-2}$  & 0.71         & 1.48 $\times 10^{-2}$  & 0.71 & 5.27 $\times 10^{-2}$ & 0.68 \\  
        25  & 1.24 $\times 10^{-2}$ & 0.75         & 1.24 $\times 10^{-2}$ & 0.75 & 4.57 $\times 10^{-2}$ & 0.61 \\  
        50  & 5.71 $\times 10^{-3}$  & 1.09         & 5.71 $\times 10^{-3}$  & 1.09 & 2.06 $\times 10^{-2}$ & 1.12 \\  
        \hline
    \end{tabular}
    
    \caption{\centering \textbf{Travelling vortex problem:} Convergence rates of $L_2$ error in $u_1$ using $ARS(1,1,1)$ coupled with type 2 discretisation.} 
    \label{tab:TV-EOC_u1b}
\end{table}

\begin{table}[h!]
    \centering
    \renewcommand{\arraystretch}{1.3} 
    \setlength{\tabcolsep}{10pt}     
    
    \begin{tabular}{|c|c|c|c|c|c|c|}
        \hline
        \multirow{2}*{\textbf{N}} & \multicolumn{2}{c|}{$\mathbf{\epsilon=10^{-1}}$} & 
                                     \multicolumn{2}{c|}{$\mathbf{\epsilon=10^{-2}}$} & 
                                     \multicolumn{2}{c|}{$\mathbf{\epsilon=10^{-3}}$} \\ 
        \cline{2-7}
        & $\mathbf{||u_2 \textbf{ error}||_{L_2}}$ & \textbf{EOC} 
        & $\mathbf{||u_2 \textbf{ error}||_{L_2}}$ & \textbf{EOC} 
        & $\mathbf{||u_2 \textbf{ error}||_{L_2}}$ & \textbf{EOC} \\ 
        \hline
        10  & 2.32 $\times 10^{-2}$  & -      & 2.32 $\times 10^{-2}$  & -      & 2.32 $\times 10^{-2}$ & -  \\  
        20  & 1.71 $\times 10^{-2}$    & 0.41 & 1.71 $\times 10^{-2}$  & 0.41 & 1.71 $\times 10^{-2}$ & 0.41  \\  
        25  & 1.51 $\times 10^{-2}$     & 0.55 & 1.51 $\times 10^{-2}$  & 0.55 & 1.51 $\times 10^{-2}$ & 0.55  \\  
        50  & 7.69 $\times 10^{-3}$    & 0.94 & 7.69 $\times 10^{-3}$  & 0.94 & 7.69 $\times 10^{-3}$ & 0.94 \\  
        \hline
    \end{tabular}
    
    \caption{\centering \textbf{Travelling vortex problem:} Convergence rates of $L_2$ error in $u_2$ using $ARS(1,1,1)$ coupled with type 2 discretisation.} 
    \label{tab:TV-EOC_u2a}
\end{table}

\begin{table}[h!]
    \centering
    \renewcommand{\arraystretch}{1.3} 
    \setlength{\tabcolsep}{10pt}     
    
    \begin{tabular}{|c|c|c|c|c|c|c|}
        \hline
        \multirow{2}*{\textbf{N}} & \multicolumn{2}{c|}{$\mathbf{\epsilon=10^{-4}}$} & 
                                     \multicolumn{2}{c|}{$\mathbf{\epsilon=10^{-5}}$} & 
                                     \multicolumn{2}{c|}{$\mathbf{\epsilon=10^{-6}}$} \\ 
        \cline{2-7}
        & $\mathbf{||u_2 \textbf{ error}||_{L_2}}$ & \textbf{EOC} 
        & $\mathbf{||u_2 \textbf{ error}||_{L_2}}$ & \textbf{EOC} 
        & $\mathbf{||u_2 \textbf{ error}||_{L_2}}$ & \textbf{EOC} \\ 
        \hline
        10  & 2.32 $\times 10^{-2}$  & -            & 2.31 $\times 10^{-2}$   & -      & 2.19 $\times 10^{-2}$ & -      \\  
        20  & 1.71 $\times 10^{-2}$  & 0.41         & 1.71 $\times 10^{-2}$  & 0.41 & 1.93 $\times 10^{-2}$ & 0.17 \\  
        25  & 1.51 $\times 10^{-2}$ & 0.55         & 1.50 $\times 10^{-2}$ & 0.55 & 1.88 $\times 10^{-2}$ & 0.10 \\  
        50  & 7.69 $\times 10^{-3}$  & 0.94         & 7.76 $\times 10^{-3}$  & 0.92 & 1.32 $\times 10^{-2}$ & 0.49 \\  
        \hline
    \end{tabular}
    
    \caption{\centering \textbf{Travelling vortex problem:} Convergence rates of $L_2$ error in $u_2$ using $ARS(1,1,1)$ coupled with type 2 discretisation.} 
    \label{tab:TV-EOC_u2b}
\end{table}

\section{Summary and Conclusions}
\label{AP_ES_Sec:Con}
In this paper, we have derived and analysed entropy stable methods in the low Mach number regime. First, an entropy inequality corresponding to convex entropy function depending on Mach number $\epsilon$ is derived for the baratropic Euler system. Then, in order to accurately and efficiently approximate low Mach number regimes of the barotropic Euler system, the IMEX-type discretisation in time is combined with a suitable entropy stable discretisation in space yielding asymptotic preserving methods. We have proposed three strategies for dicretisation in space and analysed them with respect to accuracy and entropy stability. Numerical experiments confirm entropy stable approximation in the low Mach number regime. 
\begin{comment}
 The entropy inequality corresponding to a convex entropy function depending on Mach number $\epsilon$ has been derived for the baratropic Euler system. Further, numerical schemes satisfying such an entropy stability for different values of $\epsilon$ have been developed by using IMEX-AP time discretisation and three space discretisation strategies. The entropy stability and asymptotic preserving properties of the methods have been validated numerically. 
\end{comment}

\begin{appendix}
\section{Appendix: Butcher tableau}
\label{App: AP_ES_Butcher tableau}
The first order type CK-ARS double Butcher tableau (known as ARS$(1,1,1)$) is:
\begin{equation}
\begin{tabular}{p{0.25cm}|p{0.25cm}p{0.25cm}}
\centering $0$ & $0$ & $0$ \cr 
\centering $1$ & $1$ & $0$ \cr 
\hline
 & $1$ & $0$
\end{tabular} \quad \quad
\begin{tabular}{p{0.25cm}|p{0.25cm}p{0.25cm}}
\centering $0$ & $0$ & $0$ \cr 
\centering $1$ & $0$ & $1$ \cr 
\hline
 & $0$ & $1$
\end{tabular} 
\end{equation}
The following is the 2-stage second order accurate Butcher tableau ARS$(2,2,2)$:  
\begin{equation*}
\begin{tabular}{p{0.25cm}|p{0.4cm}p{0.7cm}p{0.5cm}}
\centering $0$ & $0$ & $0$ & $0$ \cr
\centering $\gamma$ & $\gamma$ & $0$ & $0$ \cr
\centering $1$ & $\delta$ & $1-\delta$ & $0$ \cr
\hline
 & $\delta$ & $1-\delta$ & $0$
\end{tabular}
 \quad \quad
\begin{tabular}{p{0.25cm}|p{0.4cm}p{0.7cm}p{0.5cm}}
\centering $0$ & $0$ & $0$ & $0$ \cr
\centering $\gamma$ & $0$ & $\gamma$ & $0$ \cr
\centering $1$ & $0$ & $1-\gamma$ & $\gamma$ \cr
\hline
 & $0$ & $1-\gamma$ & $\gamma$
\end{tabular}
\end{equation*}
Here, $\gamma=1-\frac{1}{\sqrt{2}}$ and $\delta=1-\frac{1}{2\gamma}$.
\end{appendix}

\section*{Declaration}
The authors declare that they have no conflict of interest. Data sharing is not applicable to this article as no datasets were generated or analysed during the current study. Computational codes are available upon request.

\section*{Acknowledgements}
The authors thank Hendrik Ranocha (Mainz) and Arpit Babbar (Mainz) for fruitful discussions on this topic. 

\bibliographystyle{abbrv}
\bibliography{references} 

\end{document}